\newif\ifpictures
\picturestrue 

\documentclass{conm-p-l}
\title{First Steps in Tropical Geometry}

\author{J\"urgen Richter-Gebert}
\address{J\"urgen Richter-Gebert: Zentrum Mathematik, Technische Universit\"at M\"unchen,
Boltzmannstr.\ 3, D-85747 Garching bei M\"unchen, Germany}

\author{Bernd Sturmfels}\thanks{
Bernd Sturmfels  was partially supported by NSF grant DMS-0200729 and a John-von-Neumann Professorship during the
summer semester 2003 at Technische Universit\"at  M\"unchen. 
}
\address{Bernd Sturmfels: Department of Mathematics, University of California at Berkeley,
  Berkeley, CA 94720, USA}

\author{Thorsten Theobald}
\address{Thorsten Theobald: Zentrum Mathematik, Technische Universit\"at M\"unchen,
Boltzmannstr.\ 3, D-85747 Garching bei M\"unchen, Germany}

\date{}
\usepackage[dvips]{graphicx}

\theoremstyle{plain}
\newtheorem{thm}{Theorem}[section]
\newtheorem{lemma}[thm]{Lemma}
\newtheorem{prop}[thm]{Proposition}
\newtheorem{cor}[thm]{Corollary}

\newtheorem{defn}[thm]{Definition}
\newtheorem{ex}[thm]{Example}

\newtheorem{rmk}[thm]{Remark}

\newcommand{\zz}{\mathbb{Z}}
\newcommand{\nn}{\mathbb{N}}
\newcommand{\pp}{\mathbb{P}}
\newcommand{\tp}{\mathbb{TP}}
\newcommand{\qq}{\mathbb{Q}}
\newcommand{\rr}{\mathbb{R}}

\newcommand{\cc}{\mathbb{C}}

\begin{document}

\begin{abstract}
Tropical algebraic geometry is the
geometry of the tropical semi\-ring $(\rr,{\rm min},+)$.
Its objects are polyhedral cell complexes
which behave like complex algebraic varieties. 
We give an  introduction to this theory,  
with an emphasis on plane curves and linear spaces.
New results include a complete description of the
families of quadrics through four points in the
tropical projective plane and a counterexample
to the incidence version of Pappus' Theorem.
\end{abstract}

\maketitle

\section{Introduction}
Idempotent semirings arise in
a variety of contexts in applied mathematics,
including control theory, optimization
and mathematical physics (\cite{cgq-99,CGQ,Pin}).
 An important such semiring
is the \emph{min-plus algebra} or 
\emph{tropical semiring} $(\rr,\oplus,\odot)$. The underlying
set $ \rr $ is the set of real numbers, sometimes augmented by $+ \infty $.
The arithmetic operations of \emph{tropical addition} $\oplus$ and 
\emph{tropical multiplication} $\odot$ are 
$$ x \,\oplus \, y \,\,\, := \,\,\, {\rm min} \{x,y \} \qquad \hbox{and} \qquad
 a \,\odot \, b \,\,\, := \,\,\, a + b . $$
 The tropical semiring is idempotent in the sense that
 $\, a \,\oplus \, a \,\oplus \cdots \oplus a \,= \, a $.
 While linear algebra and matrix theory over idempotent semirings
 are well-developed and have had numerous successes in applications,
  the corresponding analytic geometry has received less attention until
  quite recently (see \cite{CGQ} and the references therein). 

 The $n$-dimensional real vector space $\rr^n$ is a module over the
  tropical semiring $(\rr,\oplus,\odot)$, with the
  operations of coordinatewise tropical addition
  $$ (a_1, \ldots, a_n) \,\oplus \, (b_1, \ldots, b_n) \quad = \quad 
  \bigl( {\rm min} \{a_1,b_1 \}, \ldots,  {\rm min} \{a_n,b_n \} \bigr). $$
  and tropical scalar multiplication (which is ``scalar addition'' classically).  $$ \lambda \, \odot \,  (a_1, a_2,\ldots, a_n)\quad = \quad 
  \bigl(\lambda + a_1,\lambda + a_2,\ldots,\lambda + a_n \bigr). $$
  Here are two suggestions of how one might define a 
  \emph{tropical linear space}.

  \medskip \noindent {\bf Suggestion 1. } A
tropical linear space $L$ is a subset of  $\rr^n$ which consists
of all solutions $(x_1,x_2,\ldots,x_n)$ to a finite system of
tropical linear equations
$$ a_1 \odot x_1 \,\oplus  \,\cdots  \,\oplus  \,a_n \odot x_n \,\,\, = \,\,\, 
b_1 \odot x_1  \,\oplus \, \cdots  \,\oplus \, b_n  \odot x_n .$$

  \medskip \noindent {\bf Suggestion 2. } 
  A tropical linear space $L$ in $\rr^n$ consists of
all tropical linear combinations 
$\, \lambda \odot a \,\oplus \,
\mu \odot b \,\oplus \,\cdots \,\oplus \,
\nu \odot c \,$ of a fixed finite
 subset $\{a,b,\ldots,c\} \subset  \rr^n $.
\smallskip

In both cases, the set $L$ is closed under tropical scalar multiplication,
$\,L \, = \, L + \rr (1,1,\ldots,1)$. We therefore identify
$L$ with its image in the \emph{tropical projective space}
$$ \tp^{n-1} \quad = \quad \rr^n / \rr (1,1,\ldots,1). $$
Let us consider the case of lines
in the tropical projective plane ($n=3$). According to
Suggestion 1, a line in $\tp^2$ would be the solution set of one linear equation
$$ a \odot x \,\oplus \, b \odot y \,\oplus \, c \odot z 
\quad = \quad
 a' \odot x \,\oplus \, b' \odot y \,\oplus \, c' \odot z . $$
 Figure~\ref{fi:linetype1} shows that such lines are one-dimensional in
 most cases, but can be two-dimensional. There is a total of
  twelve combinatorial types; see
  \cite[Figure 5]{CGQ}.

\ifpictures
\begin{figure}[h]
\vspace*{0cm}

\[
  \begin{array}{c@{\qquad}c}
  \includegraphics[height=4cm]{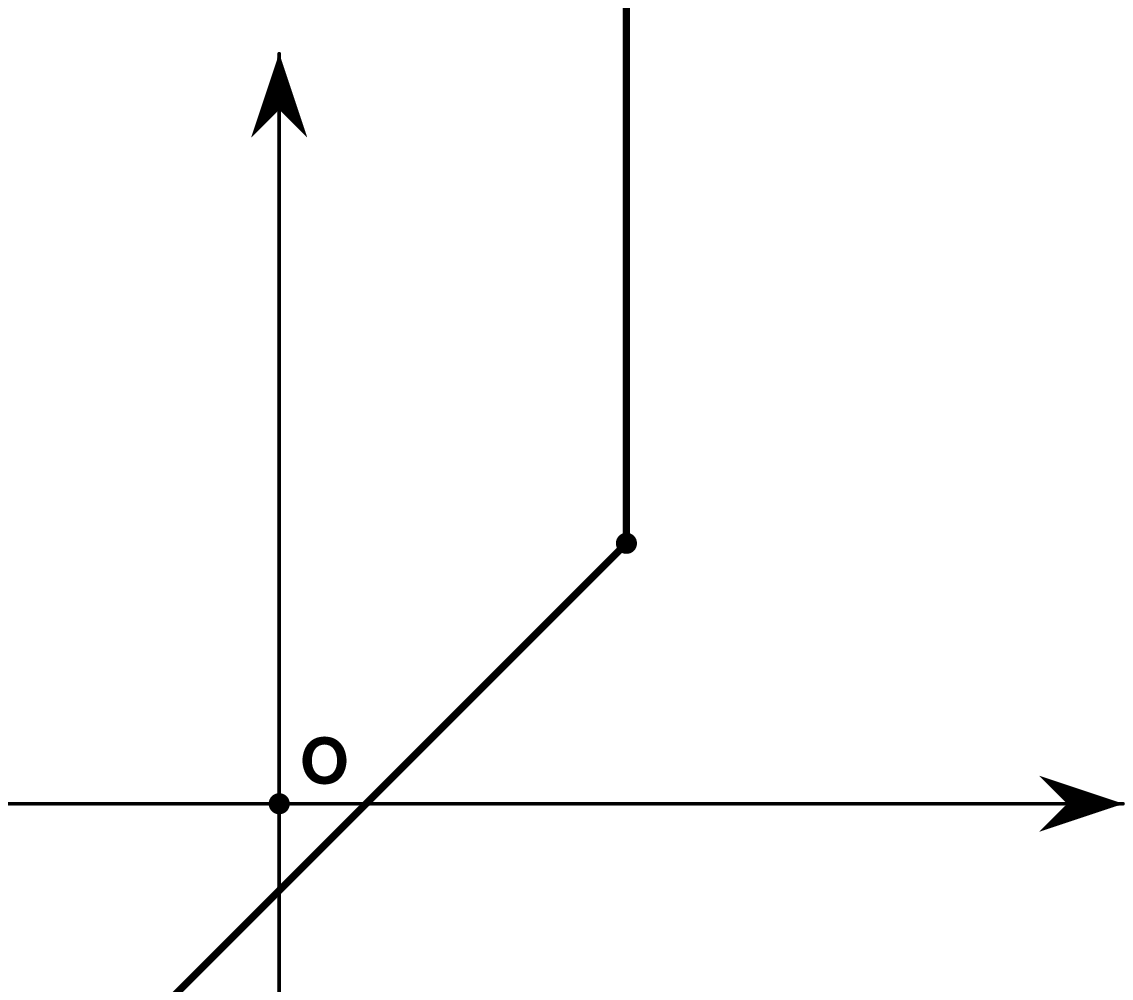} &
  \includegraphics[height=4cm]{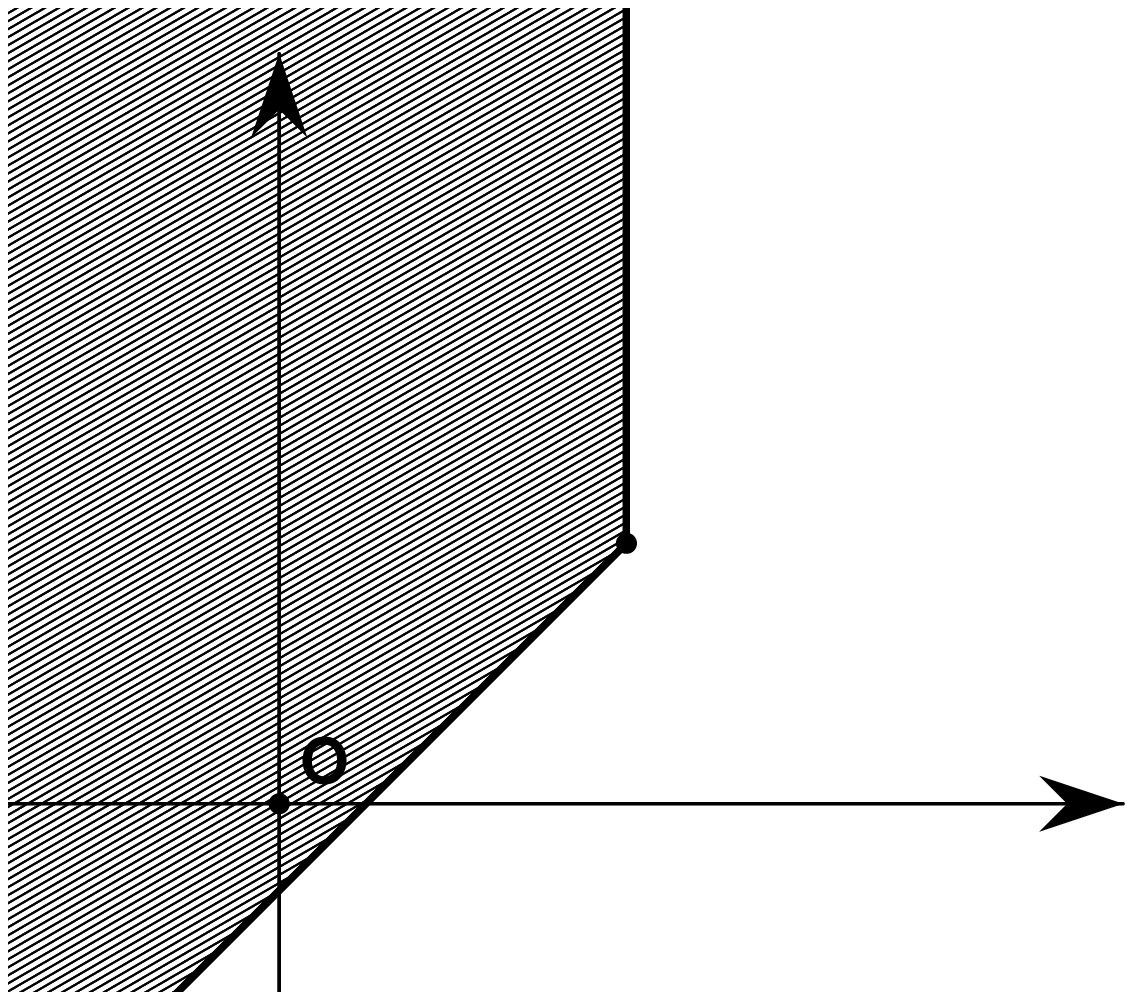} \\
  x \,= \,2 \odot x \,\oplus \, 1 \odot y \,\oplus\, 4 \odot z &
  x \,= \, x \,\oplus\, 1 \odot y \,\oplus\, 4 \odot z
  \end{array}
\]

\vspace*{0cm}

\caption{Two lines according to Suggestion 1.
 In all our pictures of the tropical projective plane $\tp^2$
we normalize with $z=0$.}
\label{fi:linetype1}
\end{figure}
\fi  

  Suggestion 2 implies that a line in $\tp^2$ is the span of  two points $a$ and $b$.
  This is the set of the following points in $\tp^2$ as the scalars
  $\lambda$ and $\mu$ range over $ \rr$:
    $$
  \lambda \odot a \oplus \mu \odot b  \quad = \quad 
  \bigl( {\rm min} \{ \lambda + a_1, \mu + b_1 \} , \,
   {\rm min} \{ \lambda + a_2, \mu + b_2 \} , \,
   {\rm min} \{ \lambda + a_3, \mu + b_3 \} \bigr). $$
  Such ``lines'' are pairs of segments connecting the two points $a$ and $b$. See 
  Figure~\ref{fi:linetype2}.

\ifpictures
\begin{figure}[h]
\vspace*{0cm}

\[
  \begin{array}{c@{\quad}c@{\quad}c}
  \includegraphics[height=3.4cm]{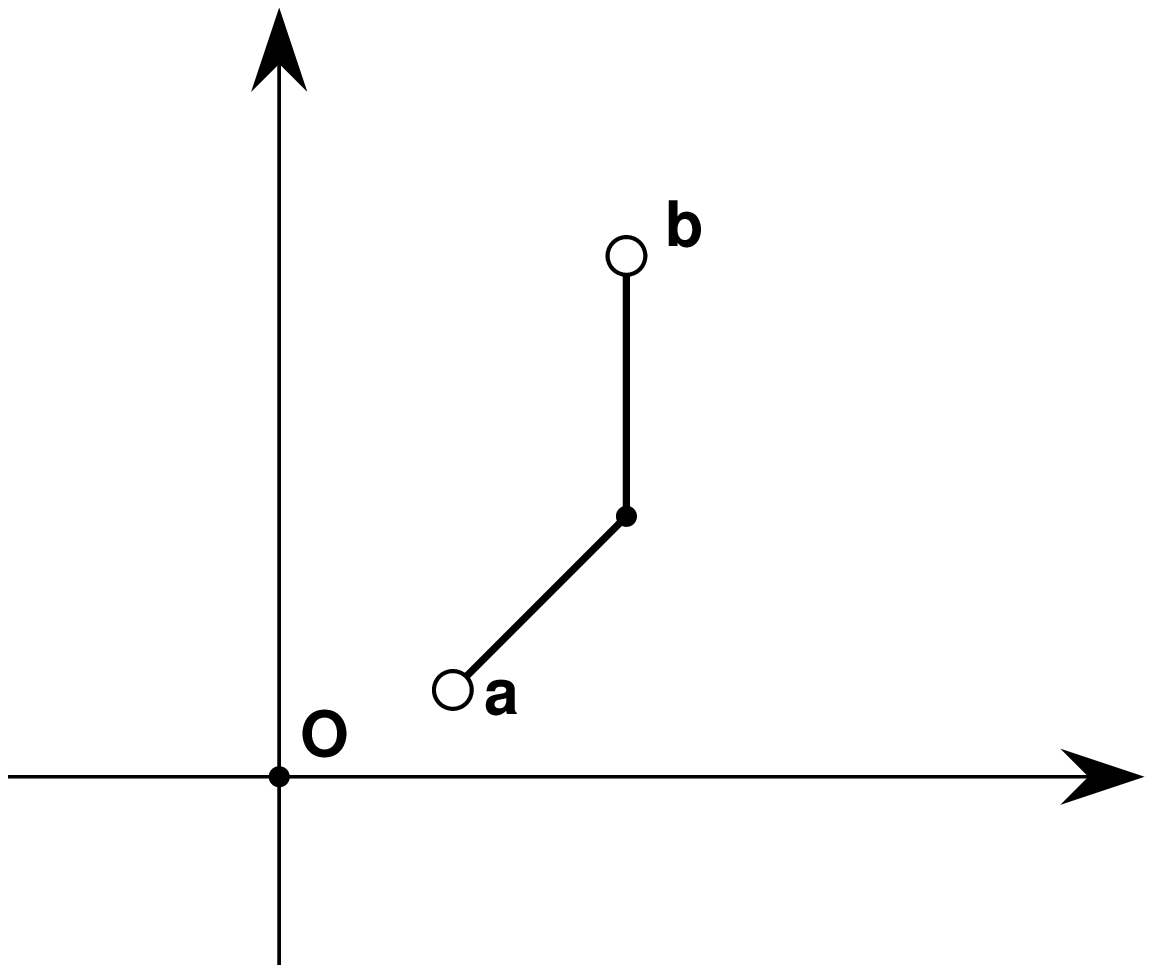} &
  \includegraphics[height=3.4cm]{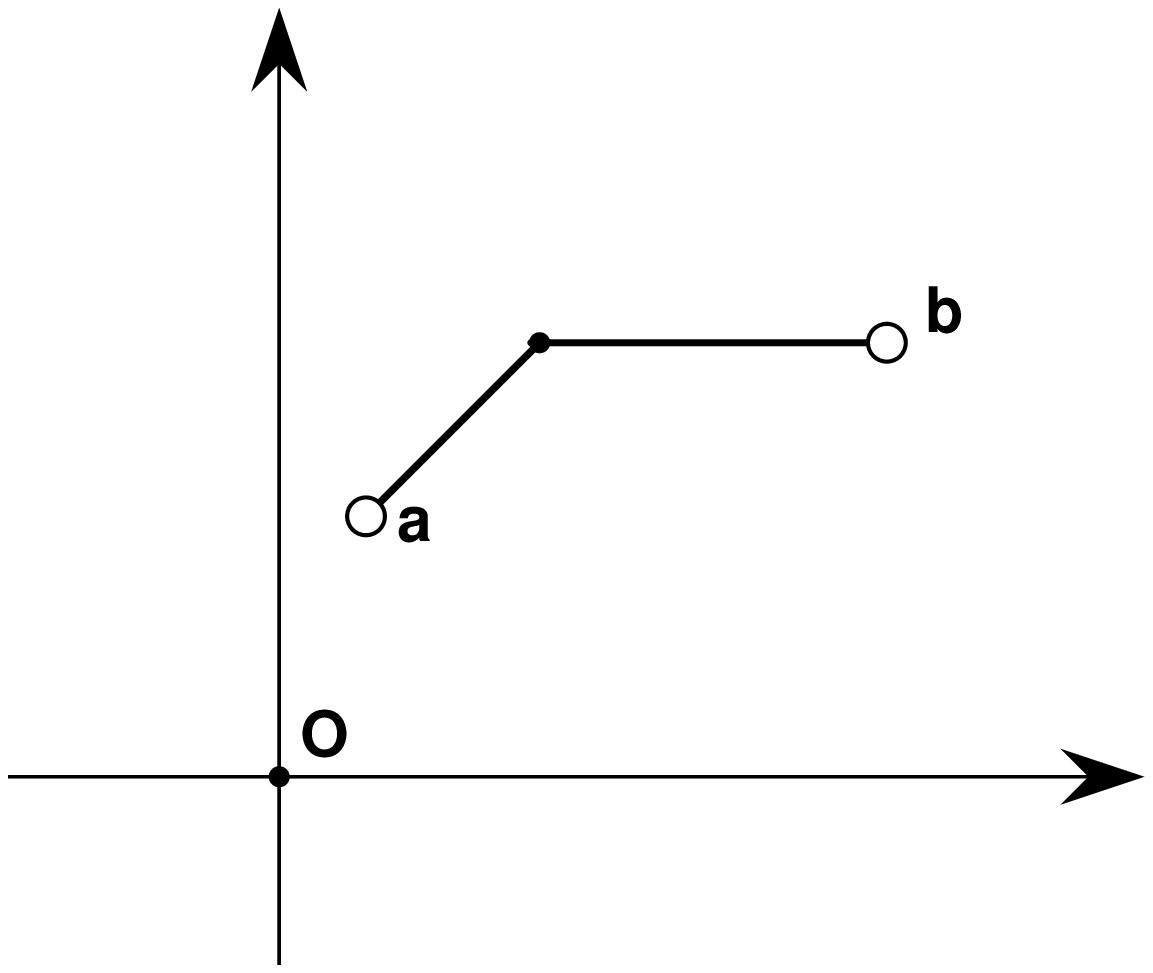} &
  \includegraphics[height=3.4cm]{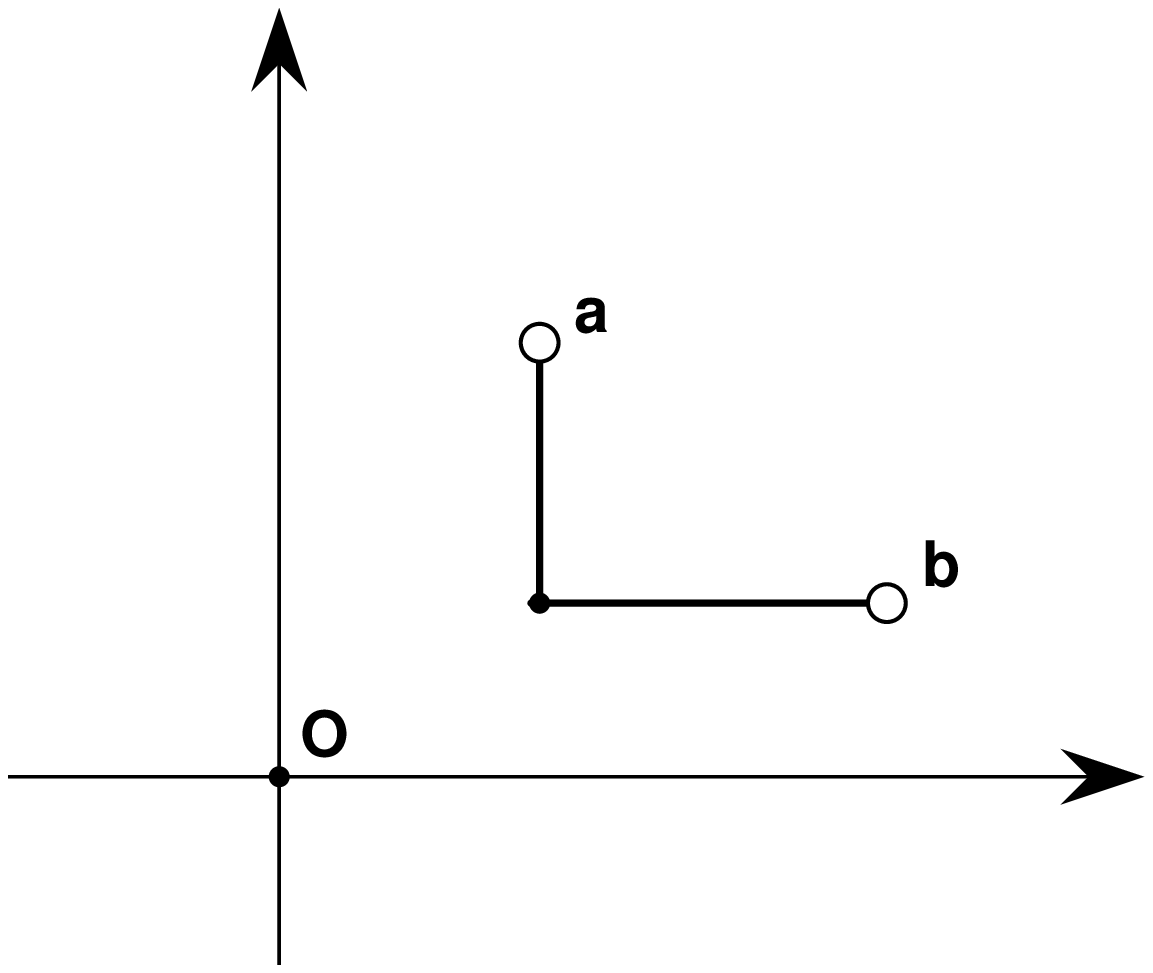}
  \end{array}
\]

\vspace*{0cm}

\caption{The three combinatorial types of lines in Suggestion 2.} 
\label{fi:linetype2}
\end{figure}
\fi

As shown in Figure~\ref{fi:span3points},
the span of three points $a$, $b$ and $c$ in $\tp^2$ is usually a two-dimensional figure.
Such figures are called \emph{tropical triangles}.

\ifpictures
\begin{figure}[h]
\vspace*{0cm}

\[
  \begin{array}{c@{\qquad}c@{\quad}c}
  \includegraphics[height=3cm]{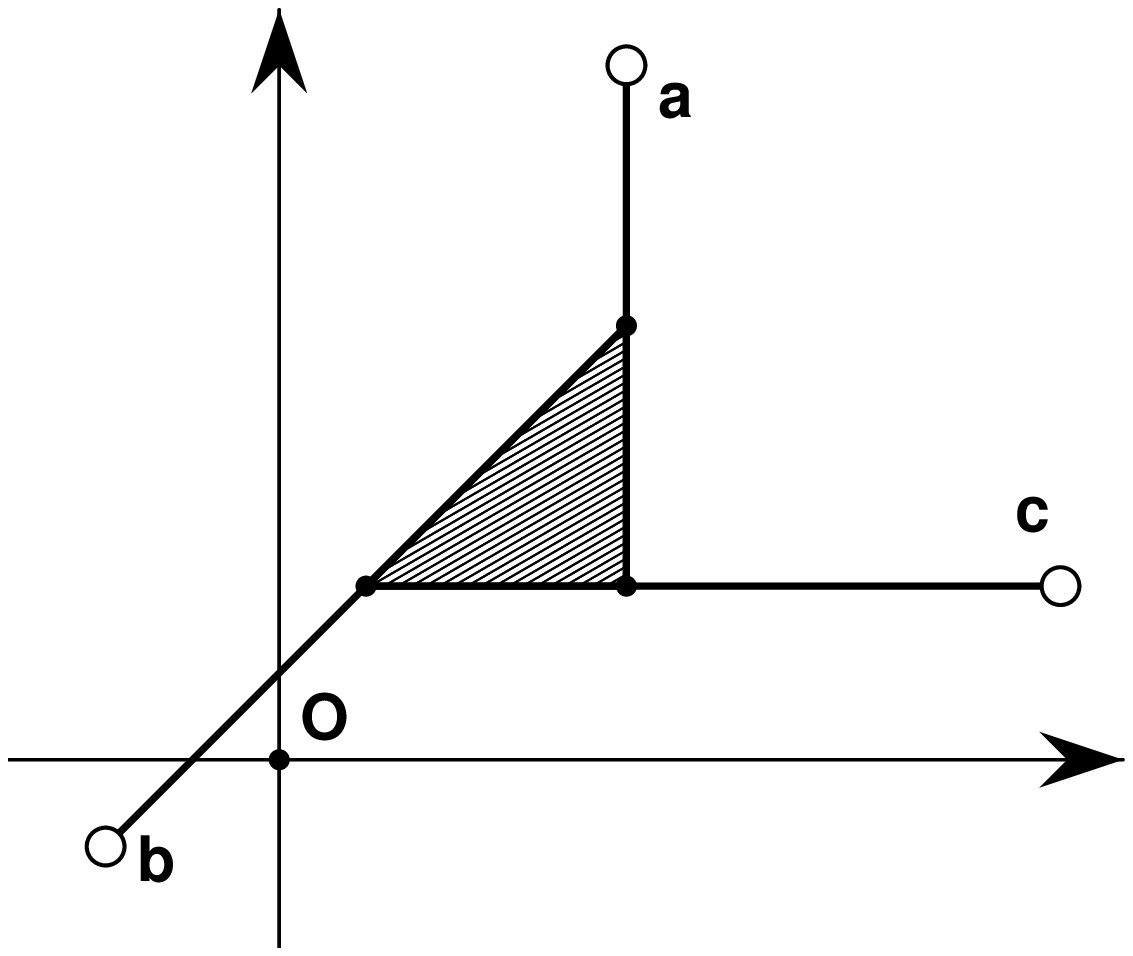} &
  \includegraphics[height=3cm]{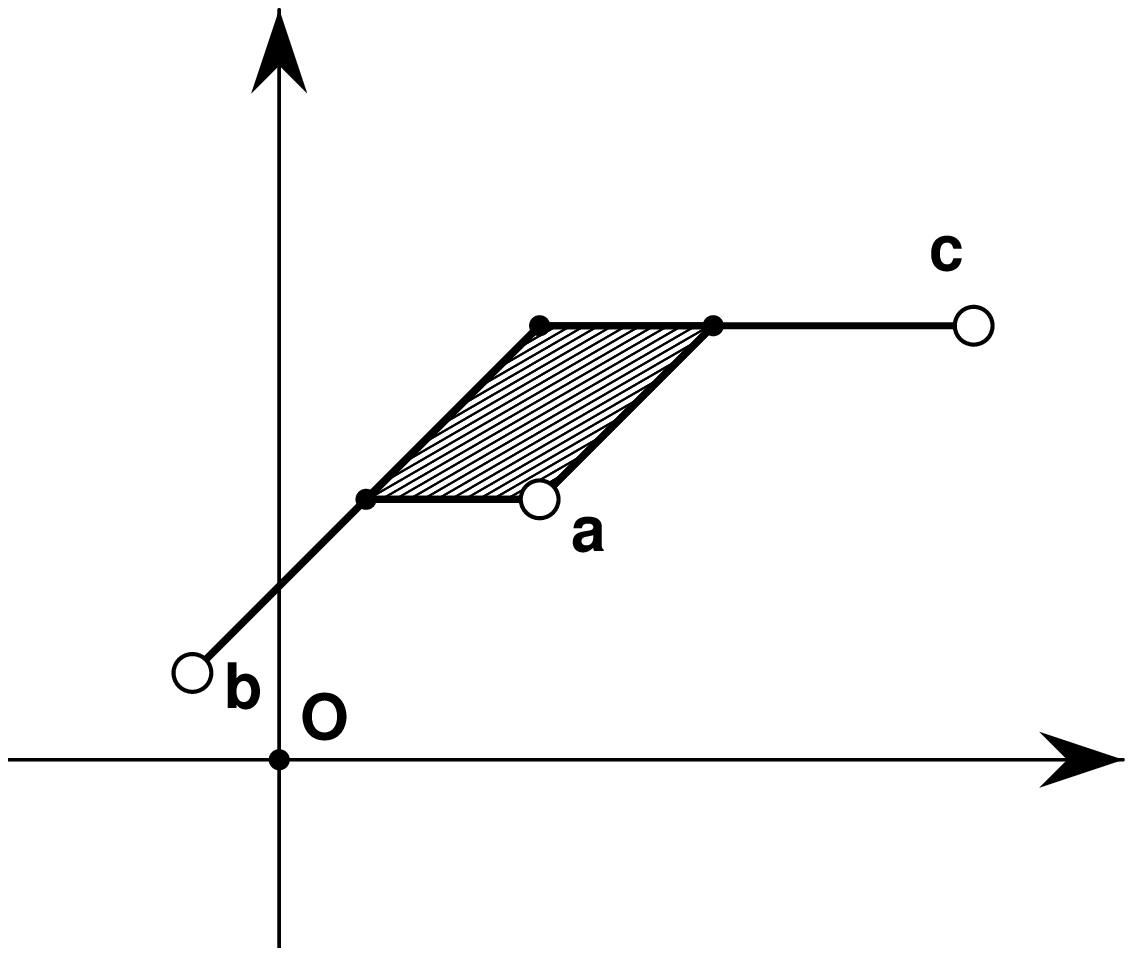} &
  \includegraphics[height=3cm]{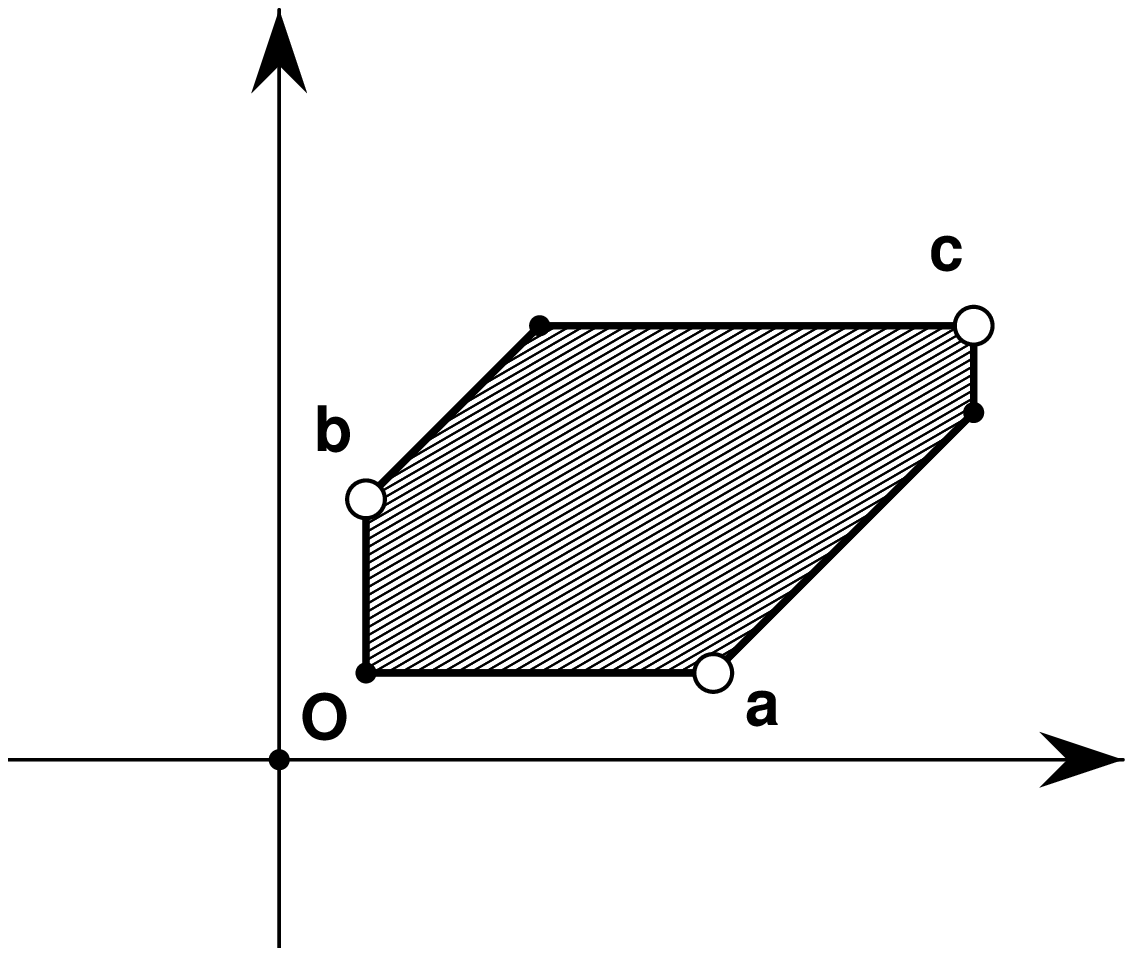} \\
  \end{array}
\]

\vspace*{0cm}

\caption{Different combinatorial types of the span of three points.}
\label{fi:span3points}
\end{figure}
\fi

It is our opinion that both of the suggested definitions of
linear spaces are incorrect. Suggestion 2 gives lines
that are too small. They are just tropical segments as in
Figure~\ref{fi:linetype2}. Also if we attempt to get the entire
plane $\tp^2$ by the construction in Suggestion 2, then we
end up only with tropical triangles as in Figure~\ref{fi:span3points}.
The lines arising from Suggestion 1 are bigger, but they are
sometimes too big. Certainly, no line should be two-dimensional.
We wish to argue that both suggestions, in spite of their
algebraic appeal, are not satisfactory for geometry. What we want is this:

  \medskip \noindent {\bf Requirement.} Lines are one-dimensional objects.
  Two general lines in the plane should meet in one point.
  Two general points in the plane should lie on one line.

  \medskip

Our  third and final definition of tropical linear space, to be
presented in the next section, will meet this requirement.
We hope that the reader will agree
that the resulting figures are
not too big and not too small but just right.

This paper is organized as follows. In Section 2 we present
the general definition of tropical algebraic varieties. In Section 3
we explain how certain tropical varieties, such as curves in the
plane, can be constructed by elementary polyhedral means.
B\'ezout's Theorem is discussed in Section 4.
In Section 5 we study linear systems of
equations in tropical geometry. The tropical
Cramer's rule is shown to compute the 
stable solutions of these systems. We apply this to
construct the unique tropical conic through any five given points
in the plane $\tp^2$.
In Section 6, we construct the pencil of conics through any four points.
We show that of the $105$ trivalent trees with six nodes representing such lines
\cite[\S 3]{speyer-sturmfels-2003} precisely $14$ trees arise
from quadruples of points in $\tp^2$.
Section 7 addresses the validity of incidence theorems in tropical
geometry. We show that Pappus' Theorem is false in general,
but we conjecture that  a certain constructive version 
of Pappus' Theorem is valid. We also report on first steps in
implementing tropical geometry in the software {\tt Cinderella} \cite{CINDERELLA}.

Tropical algebraic geometry  is an emerging field of mathematics, and
different researchers have used different names for tropical varieties:
logarithmic limit sets, Bergman fans, Bieri-Groves sets, and non-archimedean amoebas.
All of these notions are essentially the same.
Recent references include 
\cite{eklw-2003,iks-2003,mikhalkin-2002,shustin-2003,SSPE,tillmann-2003}.
For the relationship to Maslov dequantization see \cite{viro-2000}.

\section{Algebraic definition of tropical varieties}

In our algebraic definition of tropical varieties, we start from a lifting to
the field of algebraic functions in one variable. Similar liftings have 
already been
used in the max-plus literature in the context of Cramer's rule
and eigenvalue problems  (see~\cite{olsder-roos-88,schutter-demoor-97}). 
Our own version of Cramer's rule will be given in Section~\ref{se:cramer}.

The order of a rational function in one complex variable $t$
is the order of its zero or pole at the origin. It is computed as the
 smallest exponent in the numerator polynomial minus the
smallest exponent in the denominator polynomial. This definition
of order extends uniquely to the algebraic closure $\,K =
\overline{\cc(t)} \,$ of the field $\cc(t)$ of rational functions. Namely,
any non-zero algebraic function $p(t) \in K$ can be locally expressed as a  \emph{Puiseux series}
 $$ p(t) \quad = \quad c_1 t^{q_1} \,+\, c_2 t^{q_2}  \,+ \, c_3 t^{q_3} \,+ \, \cdots. $$
 Here $c_1,c_2,\ldots$ are non-zero complex numbers and
 $q_1 < q_2  < \cdots \, $ are rational numbers
 with bounded denominators. 
 The \emph{order} of  $p(t)$ is the exponent $q_1$.
  The order of an $n$-tuple of algebraic functions
 is  the $n$-tuple of their orders. This gives a map
 \begin{equation}
 \label{degreemap}
  {\rm order} \,\,: \,\, (K \backslash \{ 0 \})^n \,\, \rightarrow \,\,
 \qq^n \,\,\, \subset  \,\,\, \rr^n .
 \end{equation}
 Let $I$ be any ideal in the Laurent polynomial ring
 $\,K[x_1^{\pm 1}, \ldots, x_n^{\pm 1} ] $ and consider its
 affine variety $\, V(I) \subset  (K \backslash \{ 0 \})^n \,$
 over the algebraically closed field $K$.  The image of
 $V(I)$ under the map (\ref{degreemap}) is a subset
 of $\qq^n$. We take its topological closure. The resulting
 subset of $\rr^n$ is the tropical variety $\mathcal{T}(I)$.

 \begin{defn} \rm A \emph{tropical algebraic variety}
 is any subset of $\rr^n$ of the form
 $$ \mathcal{T}(I) \quad = \quad \overline{ {\rm order}(V(I))} ,$$
 where $I$ is an ideal in the ring of Laurent polynomials
 in $n$ unknowns with coefficients in the field $K$
 of algebraic functions in one complex variable $t$.
  \end{defn}

 An ideal $I \subset K[x_1^{\pm 1}, \ldots, x_n^{\pm 1} ] $ 
 is \emph{homogeneous} if all monomials $\,x_1^{i_1} \cdots x_n^{i_n}\,$ appearing
 in a given generator of $I$
 have the same total degree $i_1 + \cdots + i_n$. 
 Such a homogeneous ideal $I  $ 
 defines a variety $V(I)$ in
 projective space $\pp^{n-1}_K$ minus the coordinate hyperplanes $x_i = 0$.
 Its image under the order map (\ref{degreemap}) becomes a subset of tropical
 projective space $\, \tp^{n-1} \, = \, \rr^n / \rr (1,1,\ldots,1) $.

 \begin{defn} \rm A \emph{tropical projective variety}
 is a subset of $\tp^{n-1}$ of the form
 $$ \mathcal{T}(I) \quad = \quad \overline{ {\rm order}(V(I))} / \rr (1,1,\ldots,1) $$
 where $I$ is a homogeneous ideal in the Laurent polynomial ring
 $\,K[x_1^{\pm 1}, \ldots, x_n^{\pm 1}]$.
  \end{defn}

We are now prepared to give the correct definition of tropical linear space.

 \begin{defn} \rm A \emph{tropical linear space}  is a subset of
 tropical projective space $\tp^{n-1}$ of the form $\, \mathcal{T}(I) \,$ where the ideal
 $I$ is generated by linear forms
 $$ p_1(t) \cdot x_1 \,+ \,p_2(t) \cdot x_2 \,+ \,\cdots \,+ \, p_n(t) \cdot x_n $$
 whose coefficients $p_i(t)$ are algebraic functions in one complex variable $t$.
 \end{defn}

Before discussing the geometry of tropical varieties in general,
let us first see that this definition satisfies the requirement expressed in
the introduction.

\ifpictures
\begin{figure}[h]
\vspace*{0cm}

\[
  \begin{array}{c@{\qquad}c}
  \includegraphics[width=5cm]{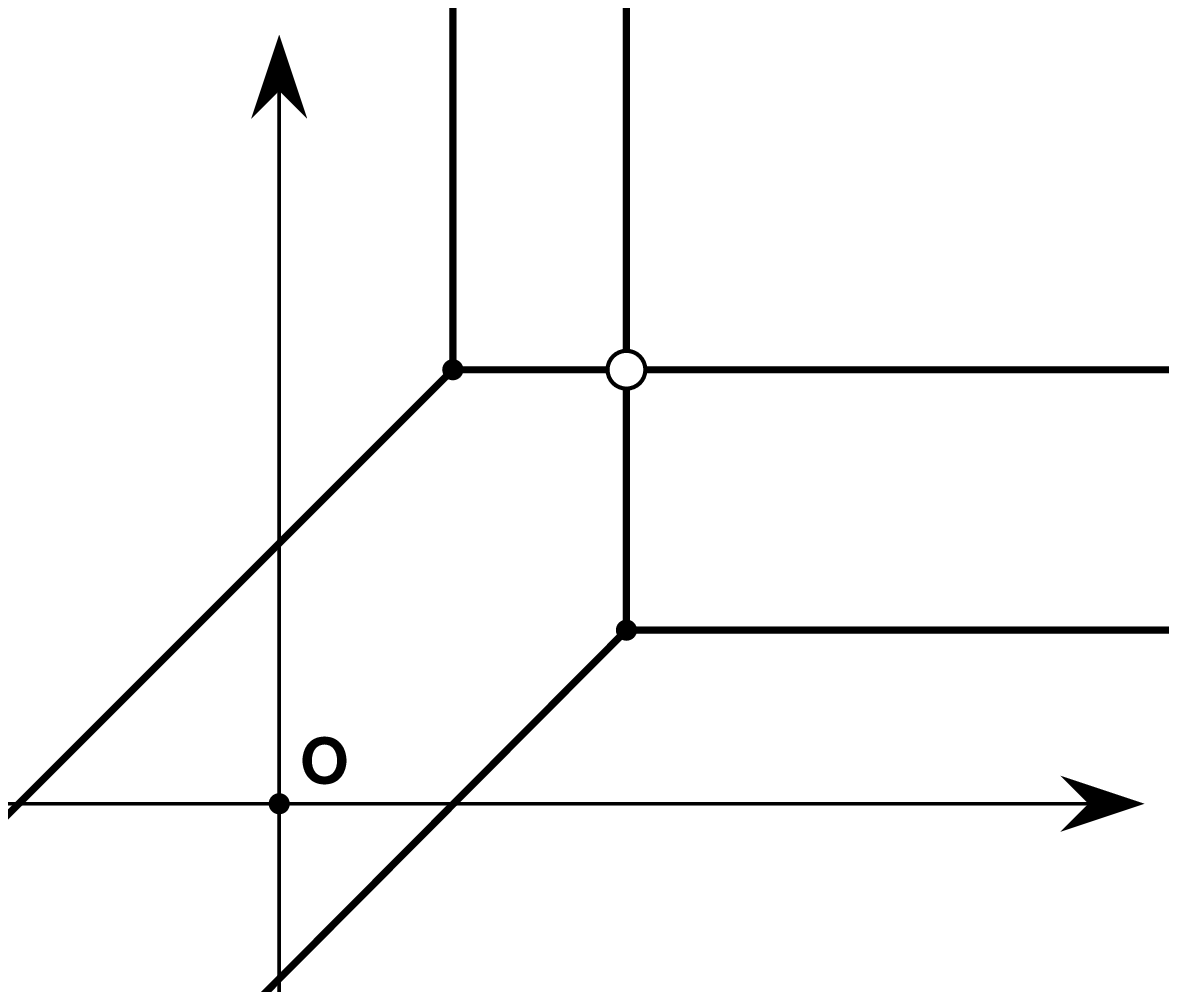} &
  \includegraphics[width=5cm]{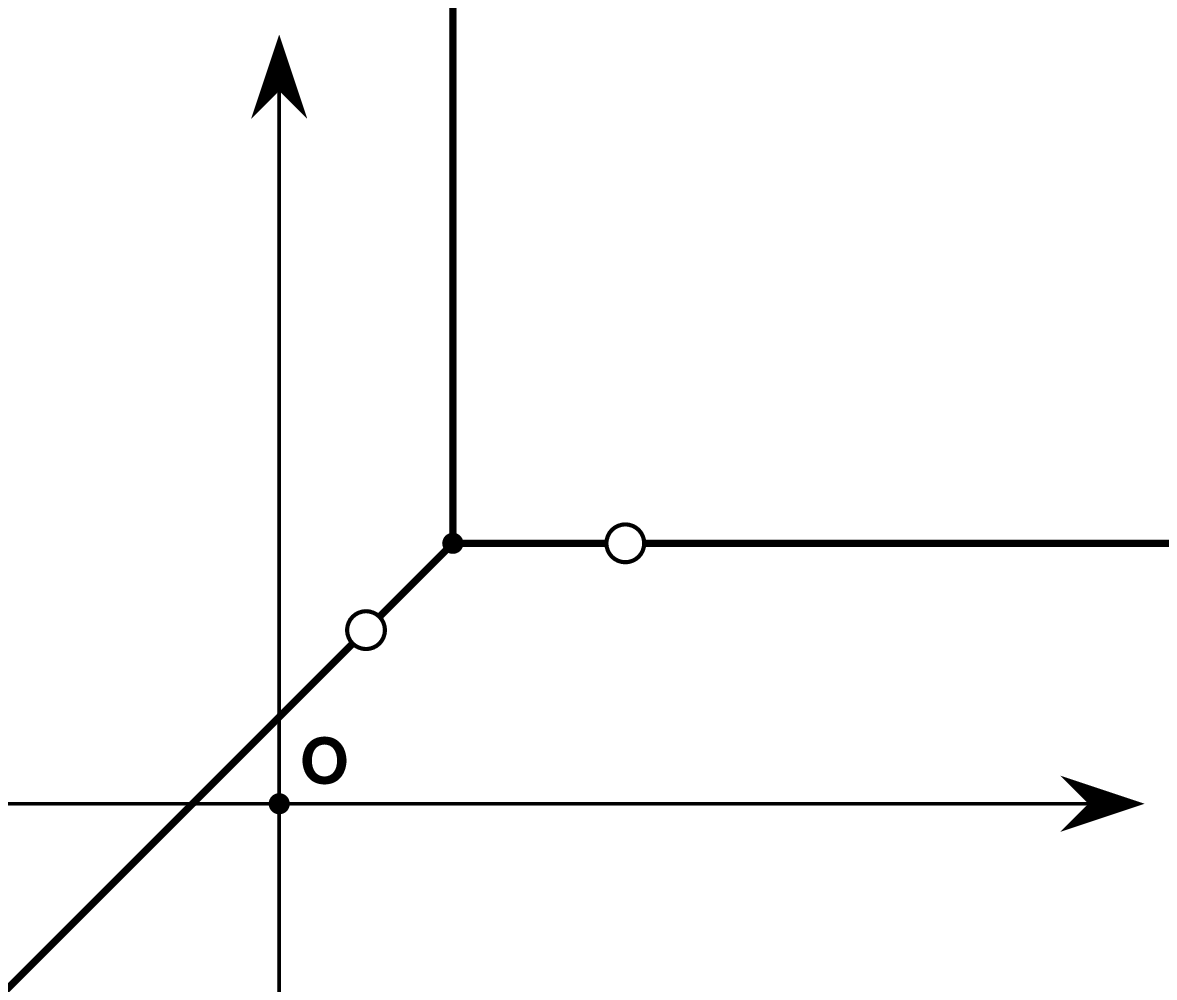} 
    \end{array}
\]

\vspace*{0cm}

\caption{Lines in the tropical plane satisfy our Requirement.}
\label{fi:tropicallines}
\end{figure}
\fi

\begin{ex} \label{ex:linetp2}
\rm \emph{Lines in the tropical plane} are defined by principal ideals
$$ I \quad = \quad \langle \,
p_1(t) \cdot x_1 \,+ \,p_2(t) \cdot x_2 \,+ \,p_3(t) \cdot x_3 \,\rangle. $$
If we abbreviate $\,a_i \,= \, {\rm order}\bigl( p_i(t) \bigr) \,$ then the line equals
\begin{eqnarray}
\mathcal{T}(I) \ = \
\bigl\{ \,(w_1,w_2,w_3) \in \tp^2 \,\,: \,
a_1+ w_1 & =  & a_2 + w_2 \leq  a_3 + w_3  \,\,\,\, \hbox{or} \nonumber \\ 
\qquad \qquad
a_1+ w_1 & = & a_3 + w_3 \leq  a_2 + w_2  \,\,\hbox{  or  } \,\, \label{eq:condlinetp2} \\
a_2 + w_2 & =  & a_3 + w_3 \leq  a_1 + w_1  \bigr\}. \nonumber
\end{eqnarray}
Thus $ \mathcal{T}(I)$ consists of three half rays emanating 
from $(-a_1,-a_2,-a_3)$ in the three
coordinate directions. 
Any two general lines meet in a unique point in $\tp^2$
and any two general points in $\tp^2$ lie on a unique line.
This is shown in Figure~\ref{fi:tropicallines}.
\end{ex}

\begin{ex} \rm \emph{Planes in the tropical $3$-space} are defined by principal ideals
$$ I \quad = \quad \langle \,
p_1(t) \cdot x_1 \,+ \,p_2(t) \cdot x_2 \,+ \,p_3(t) \cdot x_3 \,+ \,p_4(t) \cdot x_4 \,
\rangle \, . $$
Set $\,a_i =  {\rm order}\bigl( p_i(t) \bigr) \,$ as before.
Then $\mathcal{T}(I)$ is the union of six two-dimensional cones
emanating from the  point $M = (-a_1,-a_2,-a_3,-a_4)$.
See Figure~\ref{fi:plane3d} a). The intersection of two general
tropical planes is a tropical line. A line lying on a tropical plane is depicted in
 Figure~\ref{fi:plane3d} b). For a detailed algebraic discussion 
 of lines in $\tp^3$ see in Example \ref{ex:linestp3} below.

\ifpictures
\begin{figure}[h]
\vspace*{0cm}

\[
  \begin{array}{c@{\qquad}c}
  \includegraphics[height=5.5cm]{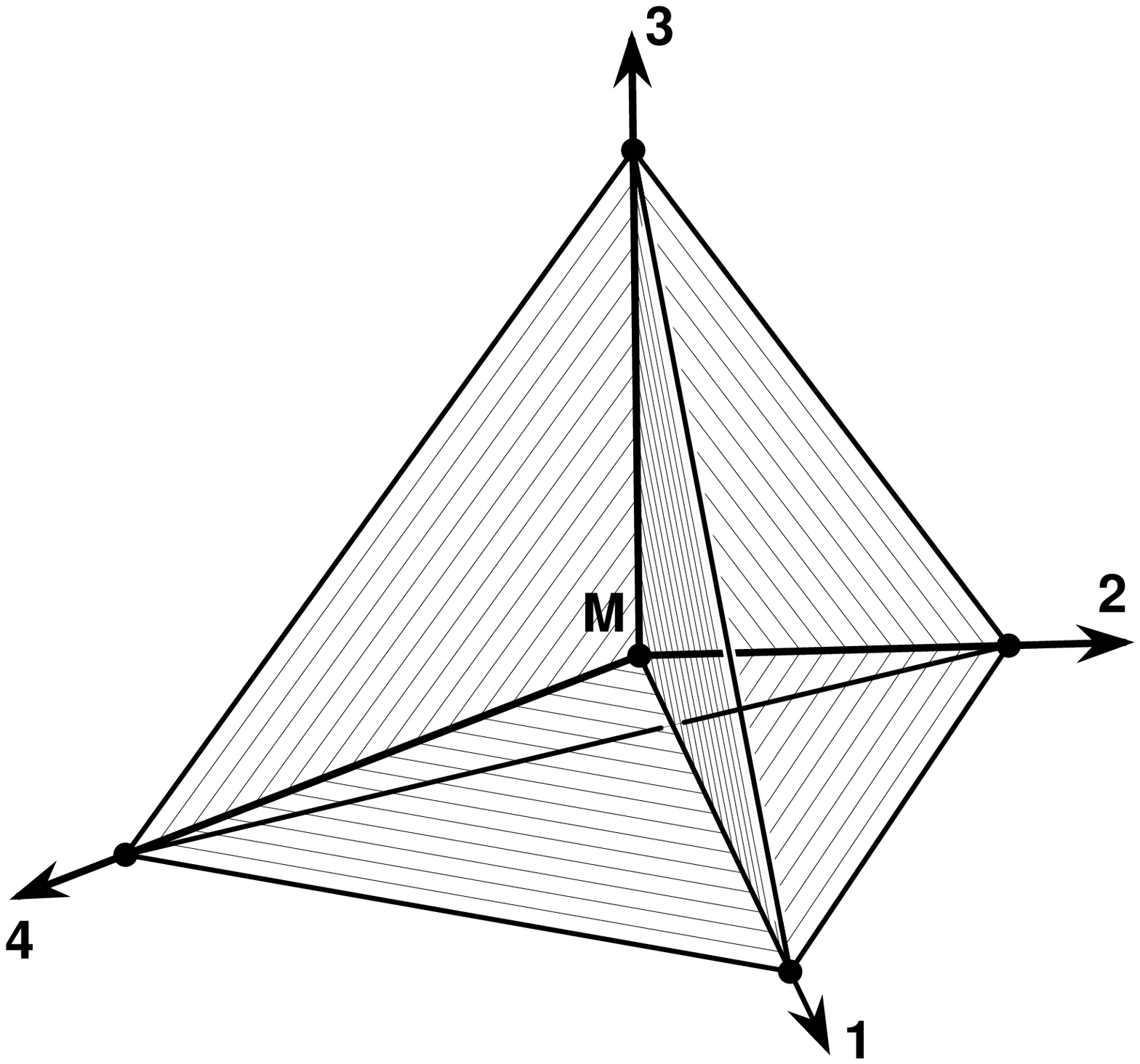} &
  \includegraphics[height=5.5cm]{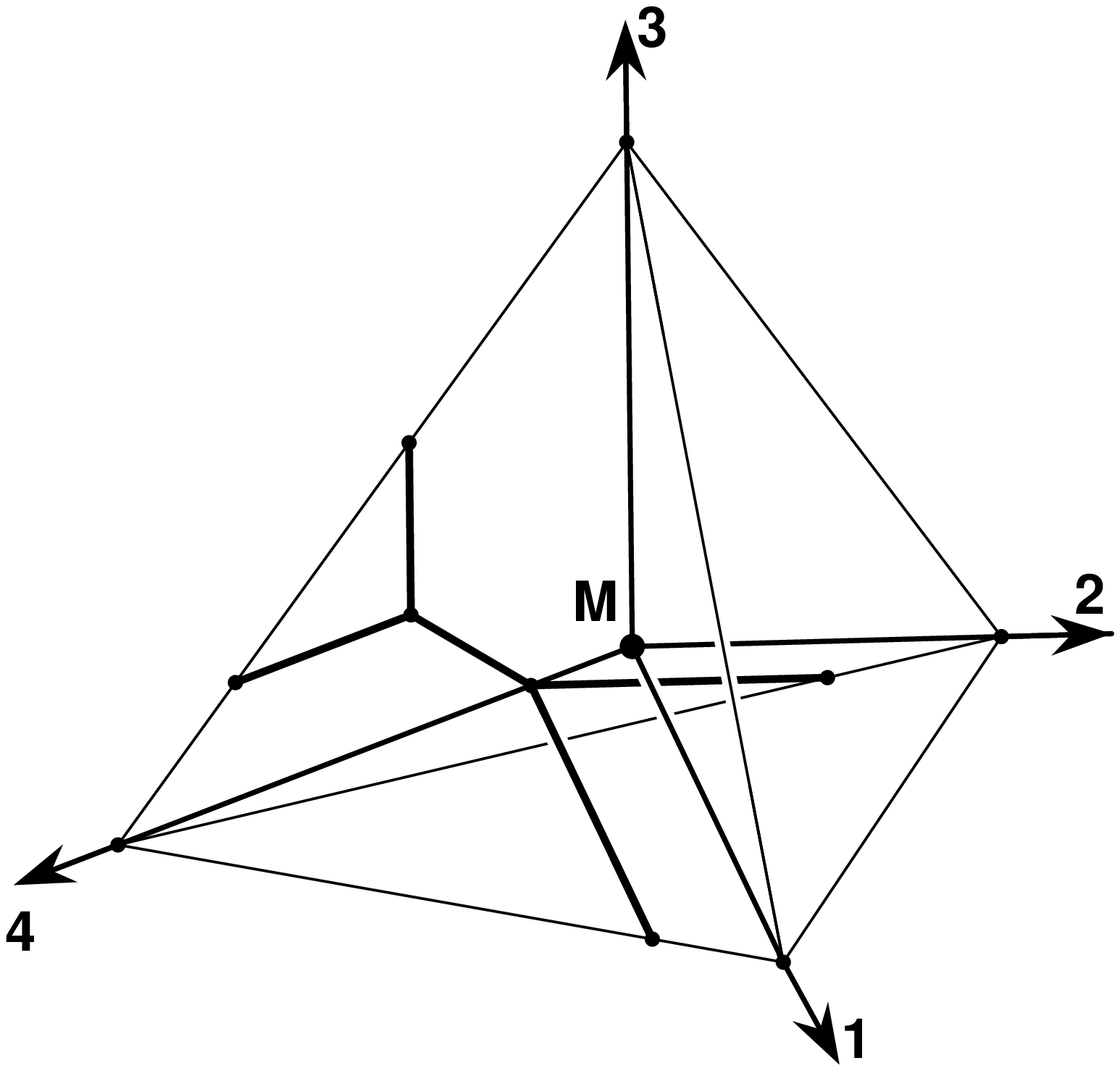} \\
  \text{a)} & \text{b)} 
  \end{array}
\]

\vspace*{0cm}

\caption{A tropical plane and a tropical line in $\tp^3$.}
\label{fi:plane3d}
\end{figure}
\fi

Our results in Section 5 imply that
any three general points in $\tp^3$ lie on a unique tropical plane.
And, of course, three general planes meet in a unique point.
\end{ex}

\medskip

Returning to the general discussion, we now present a method
for computing an arbitrary tropical variety $\mathcal{T}(I)$.
We can assume that $I$ is generated by homogeneous polynomials in
$K[x_1,\ldots,x_n]$, and, for the purpose of our discussion, we shall
regard $I$ as an ideal in this polynomial ring rather than the
Laurent polynomial ring. The input consists of an arbitrary 
generating set of the ideal $I$, and the algorithm basically
amounts to computing the \emph{Gr\"obner fan} of $I$,
as described in \cite[\S 3]{sturmfels-b96}.

 We fix a  weight vector $w \in \rr^n$. The
 \emph{weight} of the variable $x_i$ is $w_i$. The weight of a term
$\,p(t) \cdot x_1^{\alpha_1} \cdots x_n^{\alpha_n}\,$ is the real number 
$\,{\rm order}(p(t)) + \alpha_1 w_1 + \ldots + \alpha_n w_n$.
Consider a polynomial $f \in K[x_1,\ldots,x_n]$. It is a sum of
\emph{terms} $\,p(t) \cdot x_1^{\alpha_1} \cdots x_n^{\alpha_n}$.
Let $\,\overline{w}\,$ be the smallest weight among all terms in $f$.
The \emph{initial form} of $f$ equals
\[ {\rm in}_w(f) \quad = \quad 
  \sum c_{\alpha_1, \ldots, \alpha_n} \cdot 
   x_1^{\alpha_1} \cdots x_n^{\alpha_n} \, \]
where the sum ranges over all terms $\,p(t) \cdot x_1^{\alpha_1} \cdots x_n^{\alpha_n}\,$
in $f$ whose $w$-weight coincides 
with $\overline{w}$ and where 
$\,c_{\alpha_1, \ldots, \alpha_n} \in \cc \,$ denotes
the coefficient of $\,t^{\overline{w}-
 \alpha_1 w_1 - \ldots - \alpha_n w_n} \,$ 
in the Puiseux series $p(t)$. We set ${\rm in}_w(0) = 0$.
The \emph{initial ideal}  $\,{\rm in}_w(I)\,$ is defined
as the ideal generated by all initial forms
$\, {\rm in}_w(f) \,$ as $f$ runs over $I$.
For a fixed ideal $I$,  there are only finitely many initial ideals, and they
can be computed using the algorithms in \cite[\S 3]{sturmfels-b96}.
This implies the following result; 
see also \cite[\S 9]{SSPE} and \cite{speyer-sturmfels-2003}.

\begin{thm} \label{ugbthm}
Every ideal $I$ has a finite subset $\,\mathcal{G}\,$
with the following properties:
\begin{enumerate}
\item If $w \in  \mathcal{T}(I)\,$  then
$\,\{ {\rm in}_w (g) \,:\, g \in \mathcal{G}\,\bigr\}\,$
generates the initial ideal $\,{\rm in}_w(I)$.
\item If $w \not\in \mathcal{T}(I)\,$  then
$\,\{ {\rm in}_w (g) \,:\, g \in \mathcal{G}\,\bigr\}\,$
contains a monomial.
\end{enumerate}
\end{thm}

The finite set $\mathcal{G}$ in this theorem is
said to be a \emph{tropical basis} of the ideal $I$.
If $I$ is generated by a single polynomial $f$,
then the singleton $\{f \}$ is a tropical basis of $I$.
If $I$ is generated by linear forms, then the set of all \emph{circuits}
in $I$ is a tropical basis of $I$.
These are
the linear forms in $I$ whose set of variables is minimal
with respect to inclusion.
For any ideal $I$ we have
$\mathcal{T}(I) = \cap_{g \in \mathcal{G}} \mathcal{T}(\langle g \rangle)$.
We note that an earlier version of this
paper made the claim that every universal Gr\"obner basis
of $I$ is automatically a tropical basis, but this claim
is not true.

Theorem  \ref{ugbthm} implies that every tropical variety
$\mathcal{T}(I)$ is a \emph{polyhedral cell complex}, i.e., it
is a finite union of closed convex polyhedra in $\tp^{n-1}$
where the intersection of any two polyhedra is a common face.
Bieri and Groves \cite{bieri-groves-84} proved that
the dimension of this cell complex coincides
with the Krull dimension of the ring 
$K[x_1^{\pm 1}, \ldots, x_n^{\pm 1} ] /I$.
An alternative proof using Gr\"obner bases appears in \cite[\S 9]{SSPE}.

\begin{thm} \label{bierigroves}
If $V(I)$ is equidimensional of dimension $d$ then so is
$\mathcal{T}(I)$.
\end{thm}

In order to appreciate the role played by the tropical basis $\mathcal{G}$
as the representation of its ideal $I$, one needs to look at
varieties that are not hypersurfaces. The simplest example
is that of a line in the three-dimensional space $\tp^3$.

\begin{ex} \label{ex:linestp3}
\rm A \emph{line in three-space} is the tropical variety
$\,\mathcal{T}(I) \,$ of an ideal $I$ which is
 generated by a two-dimensional space of linear
forms in $K[x_1,x_2, x_3,x_4]$.
A tropical basis of such an ideal $I$ consists of four linear forms,
\begin{eqnarray*}
 U \quad =  \quad \bigl\{ &
p_{12}(t) \cdot x_2 \,+ \,p_{13}(t) \cdot x_3 \,+ \,p_{14}(t) \cdot x_4 , \,  \\ &
-p_{12}(t) \cdot x_1 \,+ \,p_{23}(t) \cdot x_3 \,+ \,p_{24}(t) \cdot x_4 , \, \\ &
-p_{13}(t) \cdot x_1 \,- \,p_{23}(t) \cdot x_2 \,+ \,p_{34}(t) \cdot x_4 , \, \\ &
\, -p_{14}(t) \cdot x_1 \, -\,p_{24}(t) \cdot x_2 \,- \,p_{34}(t) \cdot x_3 \, \bigr\},
\end{eqnarray*}
where the coefficients of the linear forms satisfy the \emph{Grassmann-Pl\"ucker relation}
\begin{equation}
\label{gpr} p_{12}(t) \cdot p_{34}(t) \,- \,
p_{13}(t) \cdot p_{24}(t) \,+ \,
p_{14}(t) \cdot p_{23}(t) \,\,\, = \,\,\, 0 . 
\end{equation}
We abbreviate $\, a_{ij} =  {\rm order}\bigl( p_{ij}(t) \bigr) $.
According to Theorem \ref{ugbthm}, the line $\mathcal{T}(I)$ is the
set of all points $w \in \tp^3$ which satisfy a Boolean combination
of linear inequalities:
\begin{eqnarray*}   & 
\bigl( \,\, a_{12} + x_2 = a_{13} + x_3 \leq a_{14} + x_4
        \quad \text{or}  \qquad \qquad \qquad\qquad \qquad \qquad \\ & \qquad \quad
a_{12} + x_2 = a_{14} + x_4 \leq a_{13} + x_3       \,\,\,\text{or} \,\,\,
a_{13} + x_3 = a_{14} + x_4 \leq a_{12} + x_2  \,\,\bigr)     \\
\text{and} & 
\bigl( \,\, 
a_{12} + x_1 = a_{23} + x_3 \leq a_{24} + x_4
        \quad \text{or}  \qquad \qquad \qquad\qquad \qquad \qquad \\ & \qquad \quad
a_{12} + x_1 = a_{24} + x_4 \leq a_{23} + x_3       \,\,\,\text{or} \,\,\,
a_{23} + x_3 = a_{24} + x_4 \leq a_{12} + x_1  \,\,\bigr)     \\
\text{and} & 
\bigl( \,\, 
a_{13} + x_1 = a_{23} + x_2 \leq a_{34} + x_4
        \quad \text{or}  \qquad \qquad \qquad\qquad \qquad \qquad \\ & \qquad \quad
a_{13} + x_1 = a_{34} + x_4 \leq a_{23} + x_2       \,\,\,\text{or} \,\,\,
a_{23} + x_2 = a_{34} + x_4 \leq a_{13} + x_1  \,\,\bigr)     \\
\text{and} & 
\bigl( \,\, 
a_{14} + x_1 = a_{24} + x_2 \leq a_{34} + x_3
        \quad \text{or}  \qquad \qquad \qquad\qquad \qquad \qquad \\ & \qquad \quad
a_{14} + x_1 = a_{34} + x_3 \leq a_{24} + x_2       \,\,\,\text{or} \,\,\,
a_{24} + x_2 = a_{34} + x_3 \leq a_{14} + x_1  \,\,\bigr)   .
\end{eqnarray*}
To resolve this Boolean combination, one distinguishes three cases
arising from (\ref{gpr}):
\[ \begin{array}{c@{\qquad}l}
\hbox{Case }{[12,34]:} \quad\quad & a_{14} + a_{23} \,= \, a_{13} + a_{24} \, \leq \,a_{12} + a_{34} , \\
\hbox{Case }{[13,24]:} \quad\quad & a_{14} + a_{23} \,= \, a_{12} + a_{34} \, \leq \,a_{13} + a_{24}  , \\
\hfill \break
\hbox{Case }{[14,23]:} \quad\quad & a_{13} + a_{24} \,= \, a_{12} + a_{34} \, \leq \,a_{14} + a_{23}  .
\end{array}
\]
In each case, the line $\mathcal{T}(I)$ consists of a line segment, with
two of the four coordinate rays emanating from each end point.
The two end points of the line segment are
\[
\begin{array}{l@{\quad}l@{\,}l}
\qquad \hbox{Case }{[12,34]}: \quad & ( a_{23} + a_{34}, \, a_{13} + a_{34}, \, a_{14}+a_{23}, \, a_{13} + a_{23}) & \text{ and } \\ 
& ( a_{13} + a_{24}, \, a_{13} + a_{14}, \, a_{12} + a_{14}, \, a_{12} + a_{13}) \, ,\\ [1ex]
\qquad \hbox{Case }{[13,24]}: \quad  & ( a_{23} + a_{34}, \, a_{13} + a_{34}, \, a_{13} + a_{24}, \, a_{13} + a_{23}) & \text{ and }\\
& ( a_{24} + a_{34}, \, a_{14} + a_{34}, \, a_{14} + a_{24}, \,  a_{12} + a_{34}) \, ,\\ [1ex]
\qquad \hbox{Case }{[14,23]}: \quad & ( a_{24} + a_{34}, \, a_{14} + a_{34}, \, a_{14} + a_{24}, \, a_{14} + a_{23}) & \text{ and } \\
& ( a_{23} + a_{34}, \, a_{13} + a_{34}, \, a_{12} + a_{34}, \, a_{13} + a_{23}) \,. \\ 
\end{array}
\]

\ifpictures
\begin{figure}[h]
\vspace*{0cm}

\[
  \begin{array}{c@{}c@{}c@{}c} 
  \includegraphics[width=4.2cm]{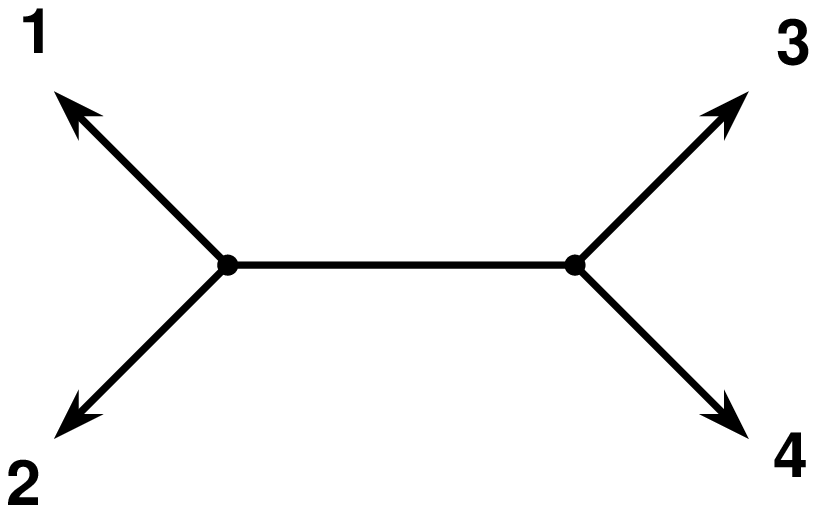} &
  \includegraphics[width=4.2cm]{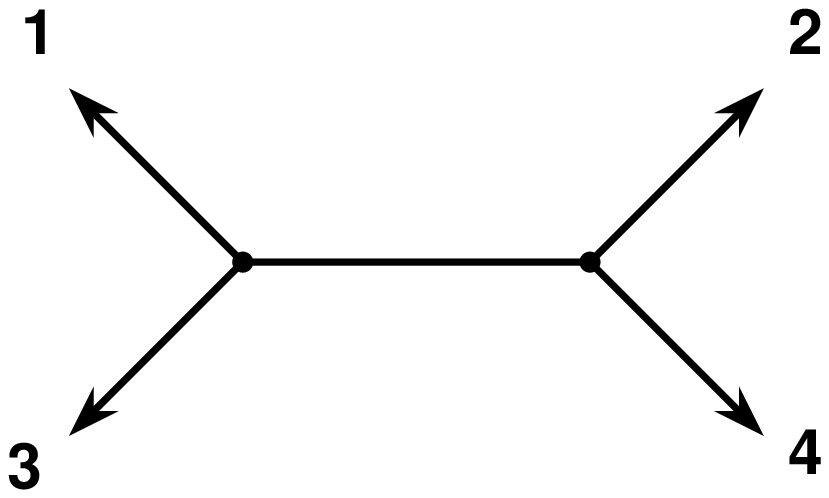} &
  \includegraphics[width=4.2cm]{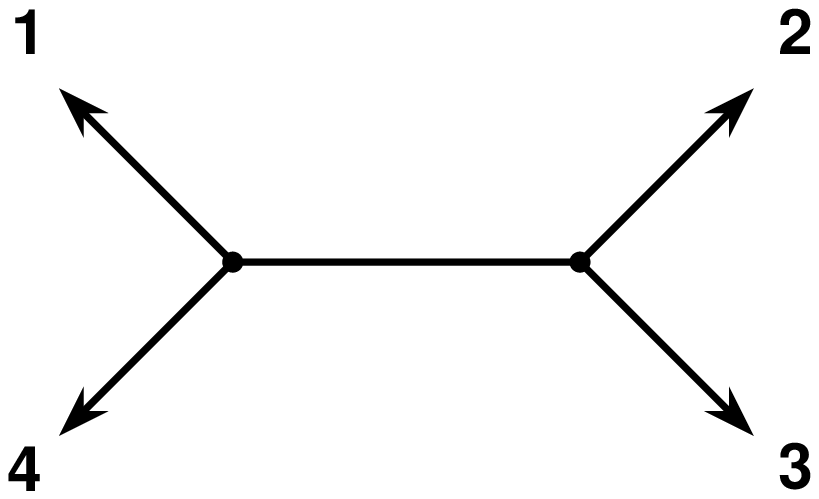} 
  \end{array}
\]

\vspace*{0cm}

\caption{The three types of tropical lines in $\tp^3$.}
\label{fi:linetypes}
\end{figure}
\fi

The three types of lines in $\tp^3$ are depicted in Figure~\ref{fi:linetypes}. Combinatorially,
they are the trivalent trees with four  labeled leaves. It was shown in
\cite{speyer-sturmfels-2003}
that lines in $\tp^{n-1}$ correspond to trivalent trees with $n$ labeled leaves.
See Example~\ref{linesinPn} below.
\end{ex}

\section{Polyhedral construction of tropical varieties}

After our excursion in the last section to polynomials over the field $K$,
let us now return to the tropical semiring
$(\rr,\oplus,\odot)$. Our aim is to derive an elementary description 
of  tropical varieties. A \emph{tropical monomial} is
an expression of the form
\begin{equation}
\label{tropmon}
 c \odot x_1^{a_1} \odot  \cdots \odot x_n^{a_n} ,
 \end{equation}
where the powers of the variables are computed tropically as well,
for instance, $x_1^3 = x_1 \odot x_1 \odot x_1$. 
The  tropical monomial
(\ref{tropmon}) represents the classical linear function
$$ \rr^n \rightarrow \rr \,\,, \,\,\,(x_1,\ldots,x_n) \mapsto  a_1 x_1 + \cdots + a_n x_n + c . $$
A \emph{tropical polynomial} is a finite tropical sum of tropical monomials,
\begin{equation}
\label{troppol}
F \quad = \quad 
 c_1 \odot x_1^{a_{11}} \odot  \cdots \odot x_n^{a_{1n}} \,\, \oplus \,\,
 \cdots  \,\, \oplus \,\,
  c_r \odot x_1^{a_{r1}} \odot  \cdots \odot x_n^{a_{rn}} .
 \end{equation}
\begin{rmk} The tropical polynomial $F$ is a piecewise
linear concave function, given as the minimum of
$r$ linear functions $\,  (x_1,\ldots,x_n) \mapsto  
a_{j1} x_1 + \cdots + a_{jn} x_n + c_j $.
\end{rmk}

We define the \emph{tropical hypersurface}
$\mathcal{T}(F)$ as the set of all points
$x = (x_1,\ldots,x_n) $ in $ \rr^n$ with the property that
$F$ is not linear at $x$. Equivalently,
$\mathcal{T}(F)$  is the set of points
$x$ at which the minimum is attained
by two or more of the linear functions.
Figure~\ref{fi:3dattainedtwice} shows 
 the graph of the
 piecewise-linear concave function $F$
 and the resulting curve $\mathcal{T}(F) \subset \rr^2$
 for a quadratic tropical polynomial $F(x,y)$.
 
\ifpictures
\begin{figure}[h]
\vspace*{0cm}

\[
  \includegraphics[width=6cm]{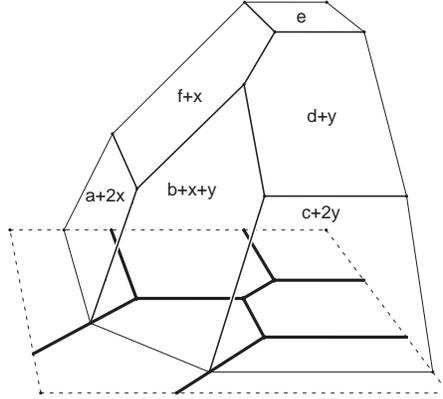}
\]

\vspace*{0cm}

\caption{The graph of a piecewise linear concave function on $\rr^2$.}

\label{fi:3dattainedtwice}
\end{figure}
\fi

\begin{lemma} \label{codim1} If $F$ is a tropical polynomial then
there exists a polynomial $f \in K[x_1,\ldots,x_n]$ such that
$ \,\mathcal{T}(F) \, = \,\mathcal{T}(\langle f \rangle)$, and vice versa.
\end{lemma}

\begin{proof}
If $F$ is the tropical polynomial (\ref{troppol}) then we define
\begin{equation}
\label{ordinarypol}
 f \quad = \quad \sum_{i=1}^r p_i(t) \cdot x_1^{a_{i1}} \cdots x_n^{a_{in}}, 
 \end{equation}
where $p_i(t) \in K $ is any Puiseux series of order $c_i$, for instance,
$p_i(t) = t^{c_i}$. Since $\mathcal{G} = \{f\}$ is a tropical basis for
its ideal, Theorem \ref{ugbthm} implies 
$ \,\mathcal{T}(F) \, = \,\mathcal{T}(\langle f \rangle)$.
Conversely, given any polynomial $f$ we can define a
tropical polynomial $F$ defining the same tropical 
hypersurface in $\rr^n$ 
by setting $c_i = {\rm order}(p_i(t))$.
\end{proof}

\begin{thm} \label{codim1thm}
Every  purely $(n-1)$-dimensional tropical variety in $\rr^n$ is
a tropical hypersurface and hence equals $\mathcal{T}(F)$ for some
tropical polynomial $F$.
\end{thm}

\begin{proof}
Let $X = \mathcal{T}(I)$ be a tropical variety of
pure dimension $n-1$ in $ \rr^n$. This means that
every maximal face of the polyhedral complex
is a convex polyhedron of dimension $n-1$.
If $P_1, \ldots, P_s$ are the minimal 
primes of the ideal $I$ then 
$$ X \quad = \quad 
 \mathcal{T}(P_1) \,\cup \,
 \mathcal{T}(P_2) \,\cup \,\cdots \,\cup
 \mathcal{T}(P_s) .$$
 Each $\mathcal{T}(P_i)$ is pure of codimension $1$,
 hence Theorem \ref{bierigroves} implies that
 $P_i$ is a codimension $1$ prime in the polynomial ring $K[x_1,\ldots,x_n]$.
  The prime ideal $P_i$
  is generated by a single irreducible polynomial,
  $\, P_i = \langle f_i \rangle$. If we set
  $f = f_1 f_2 \cdots f_s$ then $\,X = \mathcal{T}(\langle f \rangle)$,
  and Lemma \ref{codim1} gives the desired conclusion.
  \end{proof}

Theorem \ref{codim1thm} states that every tropical hypersurface 
in $\rr^n$ has an elementary construction as the locus where a
piecewise linear concave function fails to be linear.
Tropical hypersurfaces in $\tp^{n-1}$ arise in the same manner
from \emph{homogeneous tropical polynomials} (\ref{troppol}),
where $\,a_{i1} + \cdots + a_{in} \,$ is the same for all $i$.
We next describe this elementary construction for some curves in
$\tp^2$ and some surfaces in $\tp^3$.

\begin{ex} \label{ex:conictypes}
\rm \emph{Quadratic curves in the plane} are defined by 
tropical quadrics
$$ F \,\,\, = \,\,\,
 a_1 \odot x \odot x  \,\,\oplus \, \,
 a_2 \odot x \odot y   \,\,\oplus \,\,
 a_3 \odot y \odot y  \,\,\oplus \, \,
 a_4 \odot y \odot z  \,\,\oplus \, \,
 a_5 \odot z \odot z  \,\,\oplus \, \,
 a_6 \odot x \odot z  . $$
 The curve $\mathcal{T}(F)$ is a graph which has
 six unbounded edges and at most three bounded
 edges. The unbounded edges are pairs
 of parallel half rays in the three coordinate directions.
 The number of bounded edges depends on the $3 \times 3$-matrix 
  \begin{equation}
  \label{symmetric3by3}
   \begin{pmatrix} 
 a_1 & a_2 & a_6 \\
 a_2 & a_3 & a_4  \\
 a_6 & a_4 & a_5 \end{pmatrix} \, .
 \end{equation}
 We regard the  row vectors of this matrix
 as three points in $\tp^2$. If all three points are identical
  then $\mathcal{T}(F)$ is a tropical line counted with multiplicity two.
  If the three points lie on a  tropical line then
  $\mathcal{T}(F)$ is the union of two tropical lines.
    Here the number of bounded edges
 of  $\mathcal{T}(F)$ is two. In the general situation,
 the three points do not lie on a tropical line. Up to symmetry,
  there are five such general cases:

\smallskip \noindent
\emph{Case a:} $\mathcal{T}(F)$ looks like a tropical line of multiplicity two
  (depicted in Figure~\ref{fi:nonproperconics} a)). This happens if and only if
\[
   2 a_2 \, \geq \,a_1 + a_3 \quad \hbox{and} \quad
   2 a_4 \, \geq \, a_3 + a_5 \quad \hbox{and} \quad
   2 a_6 \geq \, a_1 + a_5 \, .
\]

 \smallskip \noindent
 \emph{Case b:} $\mathcal{T}(F)$ has two double half rays: 
There are three symmetric possibilities. 
The one in Figure~\ref{fi:nonproperconics}~b)
occurs if and only if
\[
   2 a_2 \, \geq \,a_1 + a_3 \quad \hbox{and} \quad
   2 a_4 \, \geq \, a_3 + a_5 \quad \hbox{and} \quad
   2 a_6 \, < \, a_1 + a_5 \, .
\]

 \smallskip \noindent
 \emph{Case c:}  $\mathcal{T}(F)$ has one double half ray:
The double half ray is emanating in the $y$-direction 
if and only if
\[
   2 a_2 \, < \,a_1 + a_3 \quad \hbox{and} \quad
   2 a_4 \, < \, a_3 + a_5 \quad \hbox{and} \quad
   2 a_6 \, \geq \, a_1 + a_5 \, .
\]
Figure~\ref{fi:nonproperconics}~c) depicts the two 
combinatorial types for this situation. They are distinguished
by whether $\,2 a_2 + a_5 - a_1 - 2 a_4\,$ is negative or positive.

\ifpictures
\begin{figure}[h]
\vspace*{0cm}

\[
  \begin{array}{c@{\;}c@{\;}c@{}c} 
  \includegraphics[width=2.8cm]{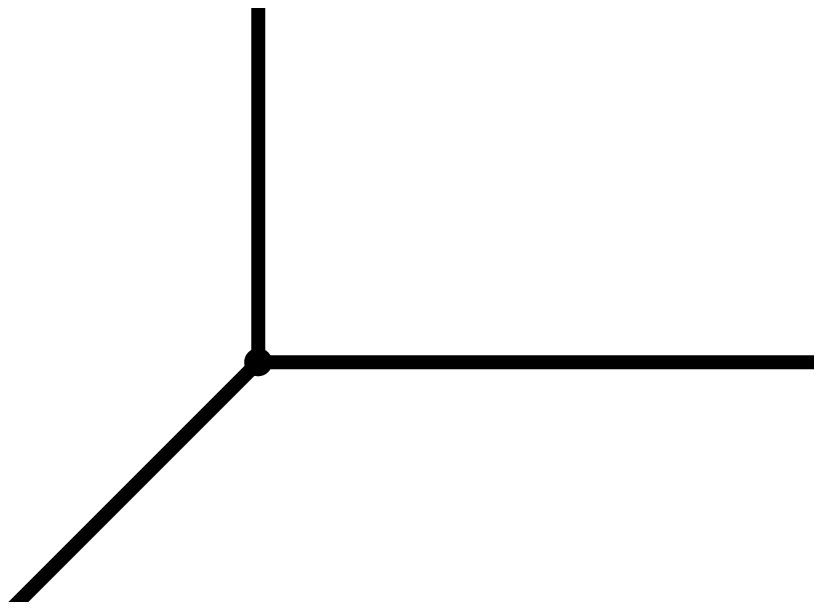} &
  \includegraphics[width=2.8cm]{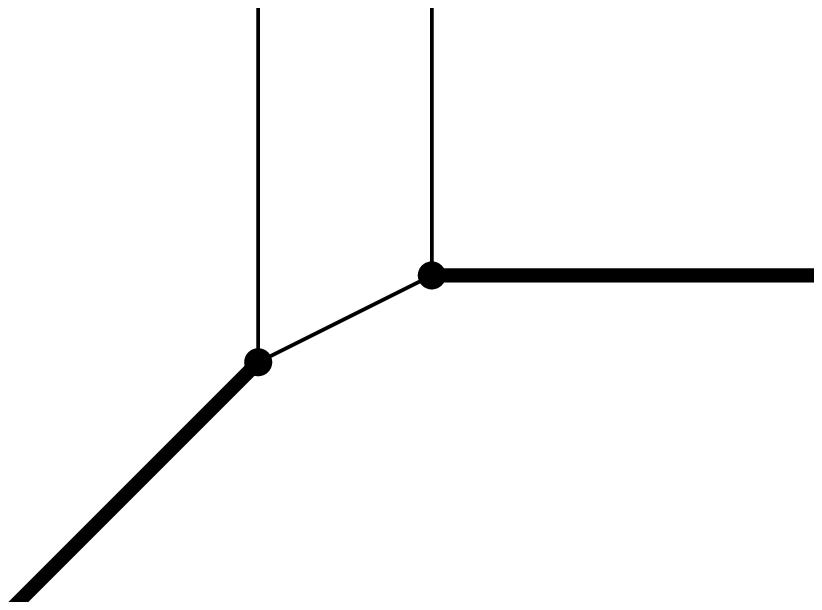} &
  \includegraphics[width=2.8cm]{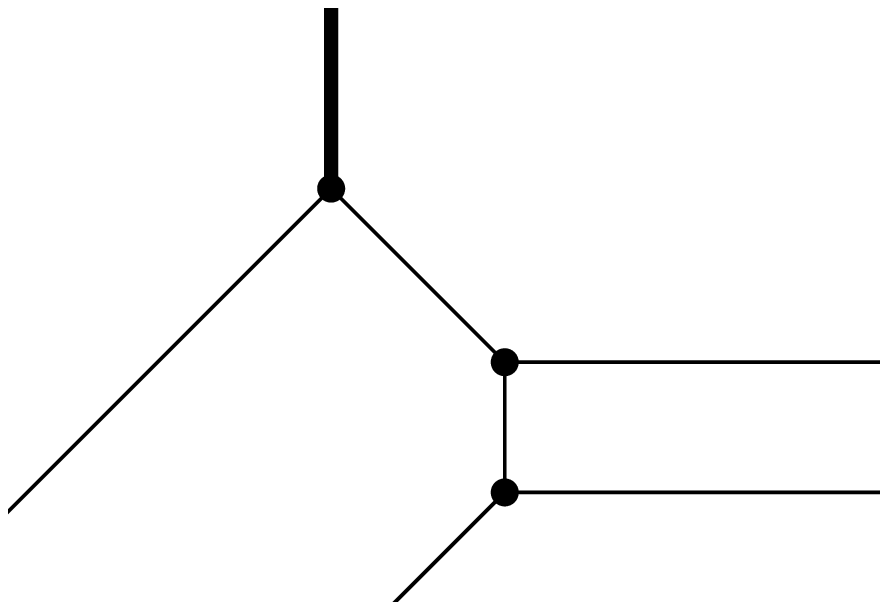} &
  \includegraphics[width=2.8cm]{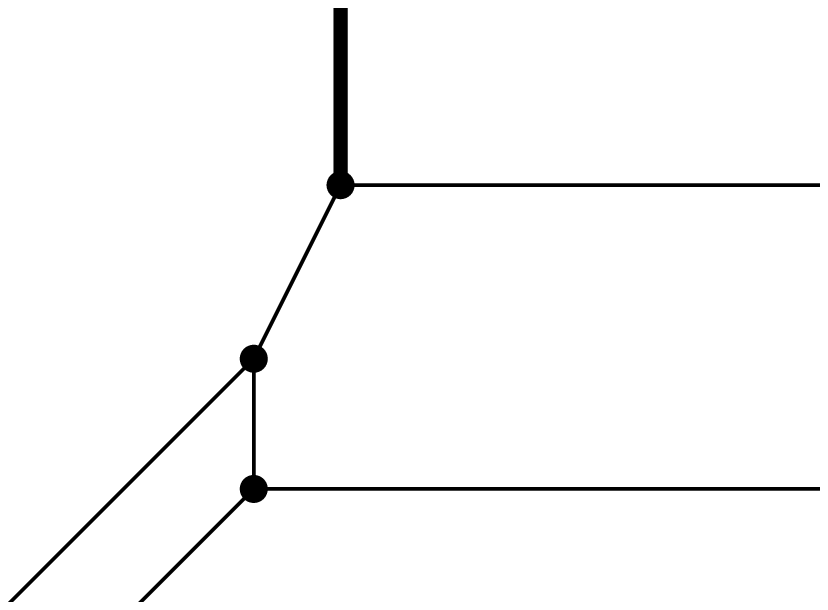} \\
  \text{a)} & \text{b)} & \multicolumn{2}{c}{\text{c)}}  
  \end{array}
\]

\vspace*{0cm}

\caption{Types of non-proper tropical conics in $\tp^2$.}
\label{fi:nonproperconics}
\end{figure}
\fi

 \smallskip \noindent
 \emph{Case d:}  $\mathcal{T}(F)$ has one vertex not on any half ray.
 This happens if and only if
\[
   a_2 + a_4 \, < \,a_3 + a_6 \quad \hbox{and} \quad
   a_2 + a_6 \, < \, a_1 + a_4 \quad \hbox{and} \quad
   a_4 + a_6  \, < \, a_2 + a_5 \, .
\]
If one of these inequalities becomes an equation, then $\mathcal{T}(F)$
is a union of two lines.

 \smallskip \noindent
\emph{Case e:}  $\mathcal{T}(F)$ has four vertices and each of them lies on some half ray.
Algebraically,
\begin{eqnarray*}
  & &  2 a_2 \, < \,a_1 + a_3 \quad \hbox{and} \quad
   2 a_4 \, < \, a_3 + a_5 \quad \hbox{and} \quad
   2 a_6 < \, a_1 + a_5 \\
  & \text{ and } & 
   (a_2 + a_4 \, > \,a_3 + a_6 \quad \hbox{or} \quad
   a_2 + a_6 \, > \, a_1 + a_4 \quad \hbox{or} \quad
   a_4 + a_6  \, > \, a_2 + a_5) \, .
\end{eqnarray*}

\ifpictures
\begin{figure}[h]
\vspace*{0cm}

\[
  \begin{array}{c@{\quad}c@{\quad}c@{\quad}c} 
  \multicolumn{2}{c}{\includegraphics[width=4.8cm]{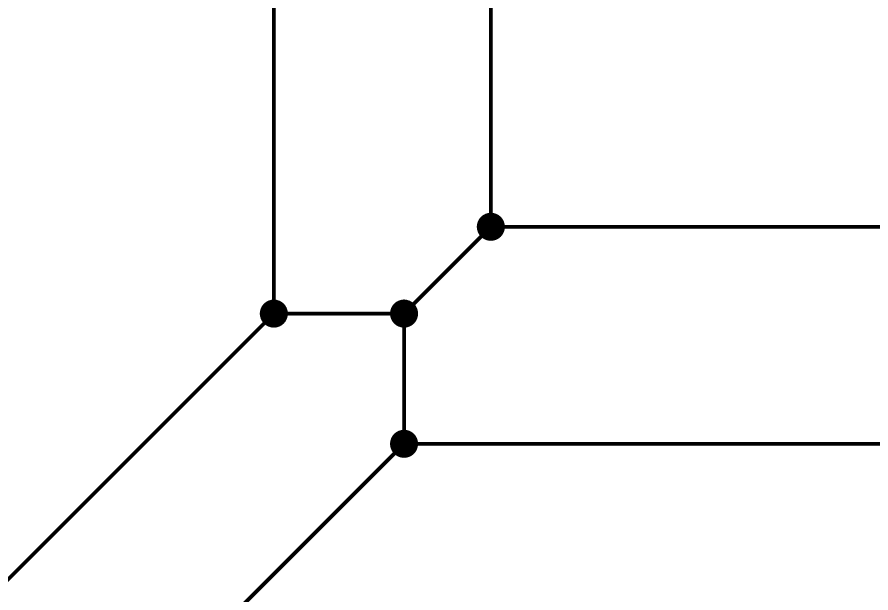}} &
  \multicolumn{2}{c}{\includegraphics[width=4.8cm]{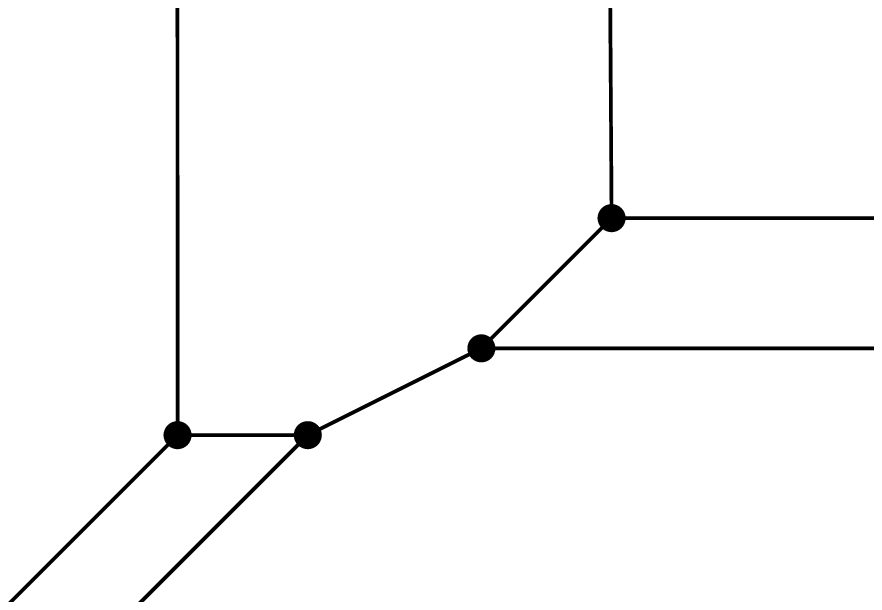}} \\  
  \multicolumn{2}{c@{\qquad}}{\text{d)}} & \multicolumn{2}{c}{\text{e)}}   \end{array}
\]

\vspace*{0cm}

\caption{Types of proper tropical conics in $\tp^2$.}
\label{fi:properconics}
\end{figure}
\fi

 \bigskip

 The curves in cases d) and e) are called \emph{proper conics}.
 They are shown in Figure~\ref{fi:properconics}.
 The set of proper conics forms a polyhedral cone.
 Its closure in  $\tp^5$ is called
 the \emph{cone of proper conics}. This cone is
 defined by the three inequalities
 \begin{equation}
\label{properconics}
 2 a_2 \, \leq \,a_1 + a_3 \quad \hbox{and} \quad
 2 a_4 \, \leq \, a_3 + a_5 \quad \hbox{and} \quad
 2 a_6 \leq \, a_1 + a_5 \, .
\end{equation}
We will see later that proper conics
play a special role in interpolation. \qed
\end{ex}

We next describe arbitrary curves in the tropical plane $\tp^2$.
Let $\mathcal{A}$ be a subset of 
$\,\{ (i,j,k) \in \nn_0^3 \, : \, i+j+k = d\}$ for some $d$ and consider
a tropical polynomial 
\begin{equation}
\label{troppoly}
  F(x,y,z) \quad = \ \sum_{(i,j,k) \in \mathcal{A}} a_{ijk} \odot x^i \odot y^j \odot z^k \, ,
  \qquad \text{where $a_{ijk} \in \rr $.}
\end{equation}
Then  $\mathcal{T}(F)$ is a tropical curve in the tropical projective plane $\tp^2$,
and Theorem~\ref{codim1thm} implies that every tropical curve in $\tp^2$
has the form $\mathcal{T}(F)$ for some $F$.

Here is an algorithm for drawing the curve 
 $\mathcal{T}(F)$ in the plane. The input to this algorithm is the \emph{support} 
 $\mathcal{A}$ and the list of coefficients $a_{ijk}$.
For any pair of points 
$(i',j',k'), (i'',j'',k'') \in \mathcal{A}$, consider the system of
linear inequalities
\begin{eqnarray*}
  & a_{i'j'k'} + i'x + j'y + k'z \ = \ 
  a_{i''j''k''} + i''x + j'' y + k''z \\
  & \qquad \qquad \qquad \qquad \qquad \qquad 
   \le  a_{ijk} + ix + jy + kz \qquad 
   \text{for } (i,j,k) \in \mathcal{A} \, .
\end{eqnarray*}
The solution set to this system is either empty or a point or a segment
or a ray in $\tp^2$.
The tropical curve $\mathcal{T}(F)$ is 
the union of these segments and rays.

It appears as if the running time of this procedure is
quadratic in the cardinality of $\mathcal{A}$, as we 
are considering arbitrary pairs of points
$(i',j',k')$ and $(i'',j'',k'')$ in $\mathcal{A}$.
However, most of these pairs can be ruled out a priori.
The following refined algorithm runs in time $O(m \log m)$ where
$m$ is the cardinality of $\mathcal{A}$. Compute the convex
hull of the points $(i,j,k,a_{ijk})$. This is a three-dimensional
polytope. The lower faces of this
 polytope project bijectively onto
the convex hull of $\mathcal{A}$ under deleting the last coordinate.
This defines a \emph{regular subdivision} $\Delta$ of $\mathcal{A}$.
A pair of vertices $(i',j',k')$ and $(i'',j'',k'')$ needs to be considered
if and only if they form an edge in the regular subdivision $\Delta$.
The segments of $\mathcal{T}(F)$ arise from the
interior edges of $\Delta$, and the rays of 
$\mathcal{T}(F)$ arise from the boundary edges of $\Delta$. This shows:

\begin{prop}
The tropical curve $\mathcal{T}(F)$ is an embedded graph in
$\tp^2$ which is dual to the regular subdivision $\Delta$ of
the support $\mathcal{A}$ of the tropical polynomial $F$.
Corresponding edges of $\mathcal{T}(F)$ and $\Delta$ are perpendicular.
\end{prop}

If the coefficients of $F$ in (\ref{troppoly}) are sufficiently generic then the
subdivision $\Delta$ is a regular triangulation. This means that the 
curve $\mathcal{T}(F)$ is a trivalent graph.

\ifpictures
\begin{figure}[h]
\vspace*{0cm}

\[\begin{array}{c@{\qquad}c}
  \includegraphics[width=5cm]{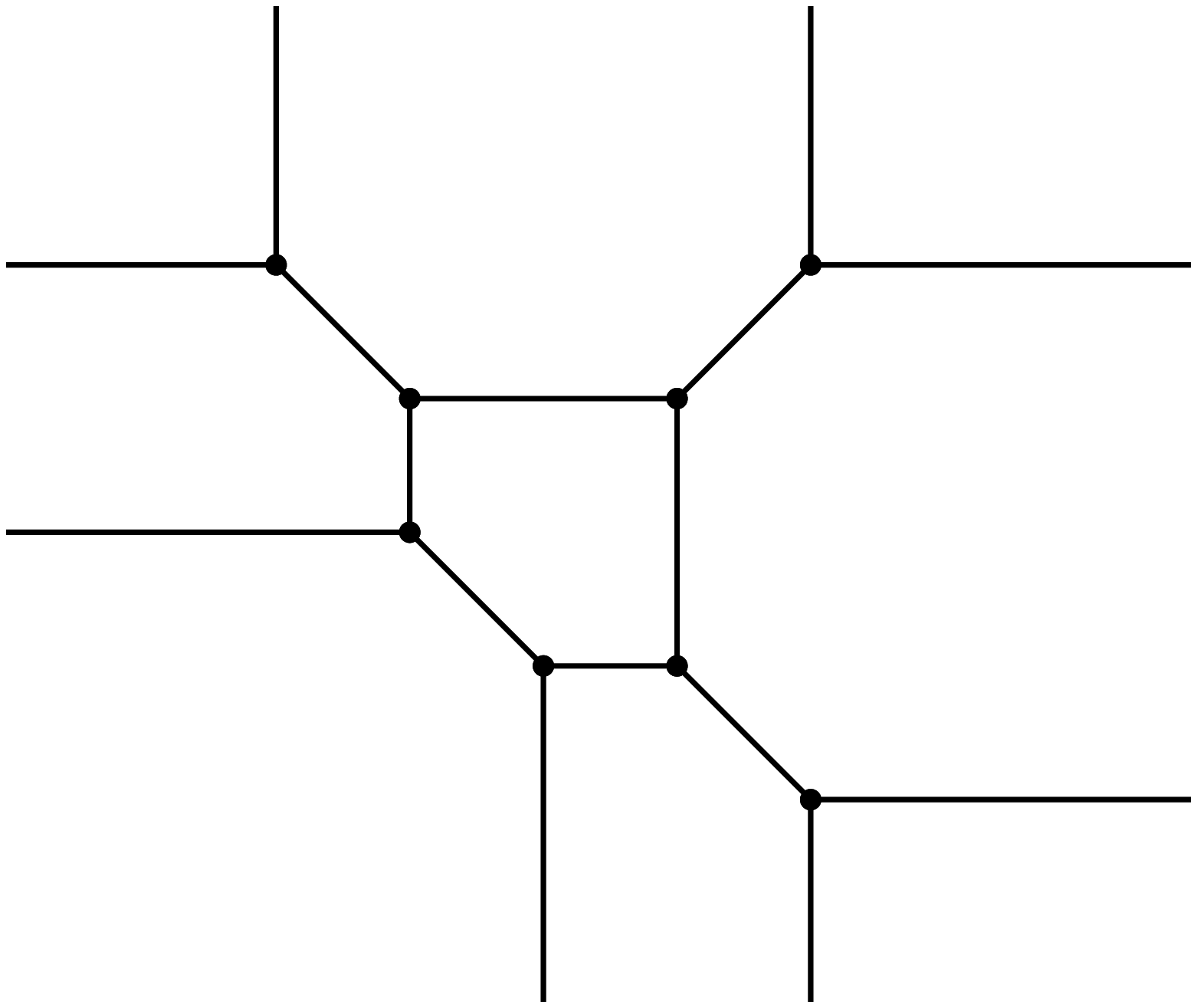}&
  \includegraphics[width=5cm]{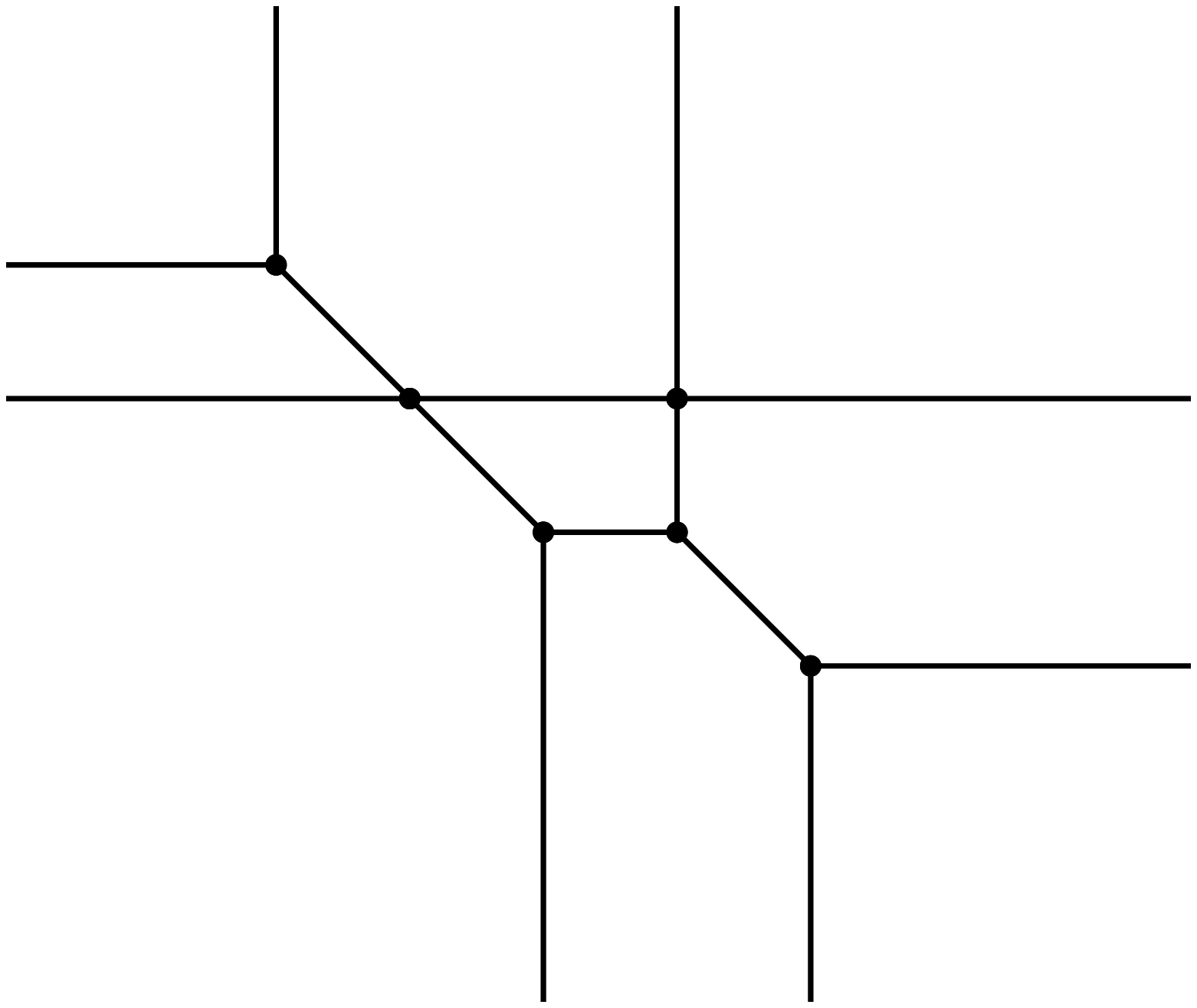}\\[.5cm]
  \includegraphics[width=2.5cm]{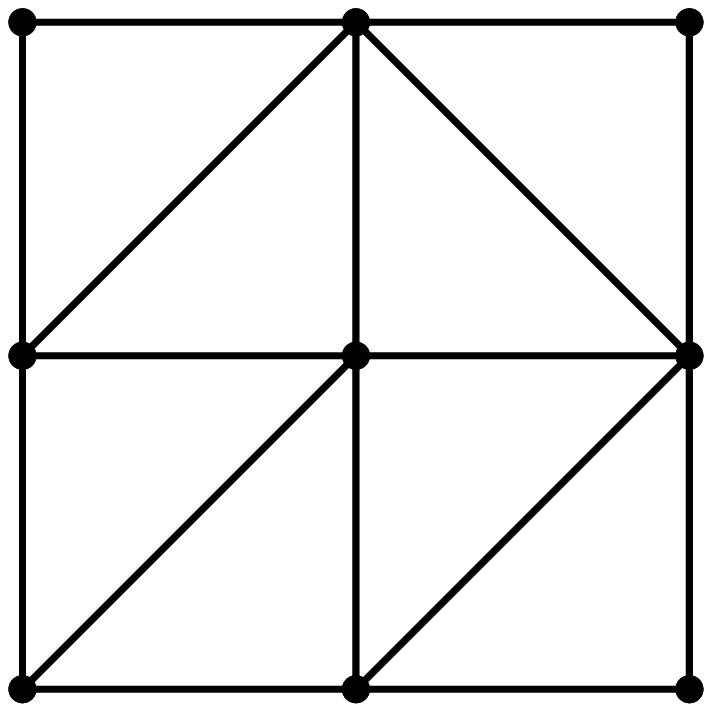}&
  \includegraphics[width=2.5cm]{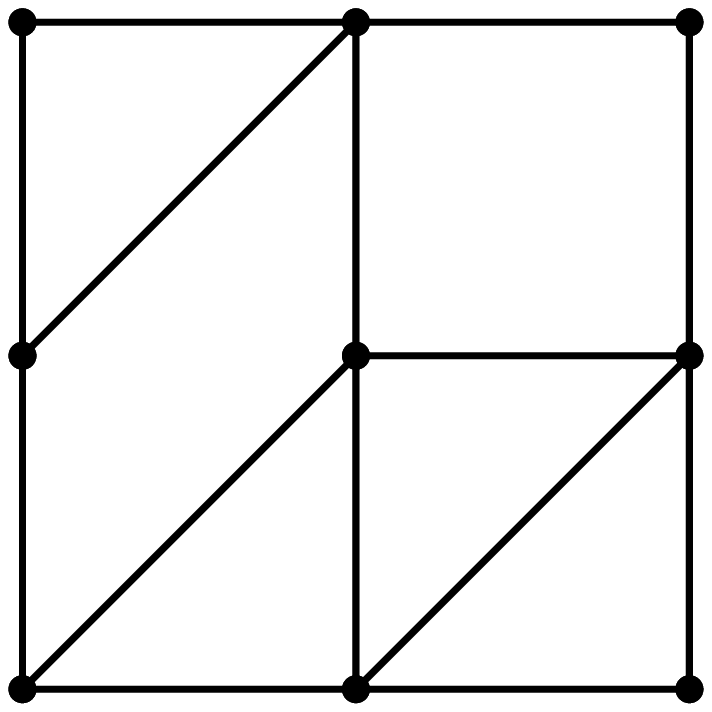}\\
\end{array}
\]

\vspace*{0cm}

\caption{Tropical biquadratic curves.}
\label{fi:elliptic}
\end{figure}
\fi

Figure \ref{fi:elliptic} shows two tropical curves whose
support is a square with side length two.
In both cases, the corresponding subdivision $\Delta$ is drawn
below the curve~$\mathcal{T}(F)$.
These are tropical versions of biquadratic curves in $\mathbb{P}^1 \times \mathbb{P}^1$, so they
represent families of elliptic curves over $K$.
The unique cycle in 
 $\mathcal{T}(F)$ arises from the interior vertex of $\Delta$.
 The same construction yields tropical Calabi-Yau hypersurfaces
 in all dimensions.

All the edges in a tropical curve $\mathcal{T}(F)$ have a natural
\emph{multiplicity}, which is the lattice length of the corresponding
edge in $\Delta$. In our algorithm,
the multiplicity is computed as the greatest
common divisor of $\,i' - i''$, $j' - j''$ and $k' - k''$. 
Let $p$ be any vertex of the tropical curve $\mathcal{T}(F)$,
 let $v_1,v_2,\ldots,v_r$ be the primitive lattice
vectors in the directions of the edges emanating from $p$,
and let $m_1,m_2,\ldots,m_r$ be the multiplicities of these
edges. Then the following \emph{equilibrium condition} holds:
\begin{equation}
\label{equilibrium}
 m_1 \cdot v_1 \, + \,
 m_2 \cdot v_2 \, + \,\cdots \,+ \,
 m_r \cdot v_r \quad = \quad 0 .
 \end{equation}
The validity of this identity can be seen by considering
the convex $r$-gon dual to $p$ in the subdivision $\Delta$.
The edges of this $r$-gon are obtained from the vectors
$m_i \cdot v_i$ by a $90$ degree rotation. But, clearly,
the edges of a convex polygon sum to zero.

Our next theorem states that this equilibrium condition 
actually characterizes tropical curves in $\tp^2$.  
This remarkable fact provides an alternative definition of
tropical curves. A subset $\Gamma$ of $\tp^2$ is a 
\emph{rational graph} if $\Gamma$ is a finite
union of rays and segments  whose endpoints and directions
have coordinates in the rational numbers $\qq$, and
each ray or segment has a positive integral multiplicity.
A rational graph $\Gamma$ is said to be \emph{balanced} 
if the condition (\ref{equilibrium}) holds at 
each vertex $p$ of $\Gamma$. 

\begin{thm}
The tropical curves in $\tp^2$ are the balanced rational
graphs.
\end{thm}

\begin{proof}
We have shown that every tropical curve is a balanced rational 
graph. For the converse, suppose that
$\Gamma$ is any rational graph which is balanced. Considering
the connected components of $\Gamma$ separately, we can assume
that $\Gamma$ is connected.
By a theorem of Crapo and Whiteley (\cite{crapo-whiteley-82},
see also \cite[Section~13.1]{richter-gebert-b96}), the balanced graph
$\Gamma$ is the projection of
the lower edges of a convex polytope. The rationality of $\Gamma$ 
ensures that we can choose this polytope to be rational. Its
defining inequalities have the form 
$ i x + j y + k z \geq a_{ijk}$ for some real numbers $a_{ijk}$.
Now define a tropical 
polynomial $F$ as in (\ref{troppol}). Then 
our algorithm implies that $\Gamma $ equals $ \mathcal{T}(F)$.
Theorem \ref{codim1thm} completes the proof.
\end{proof}

\smallskip 
Our polyhedral construction of curves can be generalized to hypersurfaces 
 in tropical projective space $\tp^{n-1}$. This raises the question whether
the class of all tropical varieties has a similar characterization.
If so, then perhaps the algebraic introduction in Section 2
was irrelevant? \ We wish to argue that this is most certainly
not the case. Here is the crucial definition.
A subset of $\rr^n$ or $\tp^{n-1}$ is 
a \emph{tropical prevariety} if it is the intersection of finitely
many tropical hypersurfaces~$\mathcal{T}(F)$.

\begin{lemma}
Every tropical variety is a tropical
prevariety, but not conversely.
\end{lemma}

\begin{proof}
By Theorem \ref{ugbthm} and Theorem \ref{codim1thm},
every tropical variety is a finite intersection of tropical hypersurfaces,
arising from the polynomials in a tropical basis.
But there are many examples of tropical prevarieties which are
not tropical varieties. Consider the tropical lines
$\mathcal{T}(L)$ and $\mathcal{T}(L')$ where 
$$
 L \,\, = \,\, 0 \odot x \,\, \oplus \,\, 0 \odot y \,\, \oplus \,\, 0 \odot z 
\quad \hbox{and} \quad
 L' \,\, = \,\, 0 \odot x \,\, \oplus \,\, 0 \odot y \,\, \oplus \,\, 1 \odot z \, . 
 $$
Then $\,\mathcal{T}(L) \, \cap \,\mathcal{T}(L')\,$ is the half ray
consisting of all points $(a,a,0)$ with $a \leq 0$. Such 
a half ray is not a tropical variety.
\end{proof}

When performing constructions in tropical algebraic geometry, 
it appears to be crucial that we work with tropical
varieties and not just with tropical prevarieties.
The distinction between these two categories is
very subtle, with the algebraic notion of a
tropical basis providing the key.
In order for a tropical prevariety  to be a tropical variety,
the defining tropical hypersurfaces must obey some strong
combinatorial consistency conditions, namely, those present among
the  Newton polytopes of some polynomials which form
 a tropical basis.

Synthetic constructions of families of tropical varieties that are
not hypersurfaces require great care. The simplest 
case is that of lines in projective space  $\tp^{n-1}$.

\begin{ex}
\label{linesinPn} \rm
A \emph{line} in $\tp^{n-1}$ is an embedded tree whose edges
are either bounded line segments or unbounded half rays, 
subject to the following three rules:
\begin{enumerate}
 \item The directions of all edges in the tree are spanned by integer vectors.
  \item There are precisely $n$ unbounded half rays. Their directions
  are the $n$ standard coordinate directions $e_1, \ldots,e_n$ in $\tp^{n-1}$.
    \item If $u_1,u_2,\ldots,u_r$ are the primitive integer vectors in the
    directions of all outgoing edges at any fixed vertex of the tree then $u_1 + u_2 + \cdots + u_r = 0$.
        \end{enumerate}      
\end{ex}

The correctness of this description follows from the
results on tropical Grassmannians in \cite{speyer-sturmfels-2003}.
We refer to this article for details on tropical linear spaces.

\section{B\'ezout's Theorem\label{se:bezout}}

In classical projective geometry, B\'ezout's Theorem states
that the number of intersection points of two general
curves in the complex projective plane is
the product of the degrees of the curves.
In this section we prove the same theorem
for tropical geometry.  The first step is to clarify
what we mean by a curve of degree $d$.

A tropical polynomial $F$ as in (\ref{troppoly}) is said to 
be a  \emph{tropical polynomial of degree $d$} if its support
$\mathcal{A}$ is {\bf equal to} the set $\{ (i,j,k) \in \nn_0^3 \, : \, i+j+k = d\}$.
Here the coefficients $a_{ijk}$ can be any real numbers,
including $0$. Changing a coefficient $a_{ijk}$ to $0$ does not alter
the support of a polynomial. After all, $0$ is the neutral element
for multiplication $\odot$ and not for addition $\oplus$.
 Deleting a term from the polynomial $F$
and thereby shrinking its support
corresponds to changing $a_{ijk}$  to $+ \infty$.
If $F$ is a tropical polynomial of degree $d$ then
we call $\mathcal{T}(F)$ a \emph{tropical curve of degree $d$}.

\begin{ex} \rm Let $d = 2$ and consider the following tropical polynomials:
\begin{eqnarray*}
 F_1 & = & 3 x^2 \,\oplus \, 5 x y  \,\oplus \, 7 y^2 \,\oplus \, 11 xz
 \,\oplus \, 13 yz \,\oplus \, 17 z^2, \\
  F_2 & = & 3 x^2 \,\oplus \, 5 x y  \,\oplus \, 7 y^2 \,\oplus \, 11 xz
 \,\oplus \, 13 yz \,\oplus \, 0  z^2 ,\\
F_3 &= & 0 x^2 \,\oplus \, 0 x y  \,\oplus \, 0 y^2 \,\oplus \, 0 xz
 \,\oplus \, 0 yz \,\oplus \, 0  z^2 ,\\
F_4 &= & 3 x^2 \,\oplus \, 5 x y  \,\oplus \, 7 y^2 \,\oplus \, 11 xz
 \,\oplus \, 13 yz \,\oplus \, (+  \infty)  z^2, \\
 F_5 &= & 3 x^2 \,\oplus \, 5 x y  \,\oplus \, 7 y^2 \,\oplus \, 11 xz
 \,\oplus \, 13 yz .
\end{eqnarray*}
$\mathcal{T}(F_1)$,
$\mathcal{T}(F_2)$ and
$\mathcal{T}(F_3)$ are tropical curves of degree $ 2$.
$\mathcal{T}(F_4) = \mathcal{T}(F_5)$ is a tropical curve, but it does 
not have a degree $d $. Can you draw these curves ? \qed
\end{ex}

In order to state  B\'ezout's Theorem, we need to define
 intersection multiplicities for two balanced rational graphs in $\tp^2$.
 Consider two intersecting line segments with rational slopes,
 where the segments have multiplicities $m_1$ and $m_2$
 and where the primitive direction vectors are $(u_1,u_2,u_3) , (v_1,v_2,v_3) \in \zz^3/ \zz (1,1,1)$.
Since the line segments are not parallel, the following determinant is
nonzero:
\[
  \det \left( \begin{array}{ccc}
    u_1 & u_2 & u_3 \\
    v_1 & v_2 & v_3 \\
     1 & 1 & 1
  \end{array} \right)
\] 
The \emph{(tropical) multiplicity} of the intersection point is
defined as the absolute value of this determinant times $m_1$ times $m_2$.

\begin{thm} \label{bezout}
Consider two tropical curves $C$ and $D$ of degrees $c$ and $d$
in the tropical projective plane $\tp^2$.
If the two curves intersect in finitely many points
then the number of intersection points, counting
multiplicities, is equal to $\,c \cdot d$.
\end{thm}

\ifpictures
\begin{figure}[h]
\vspace*{-.2cm}

\[
  \includegraphics[width=4.5 cm]{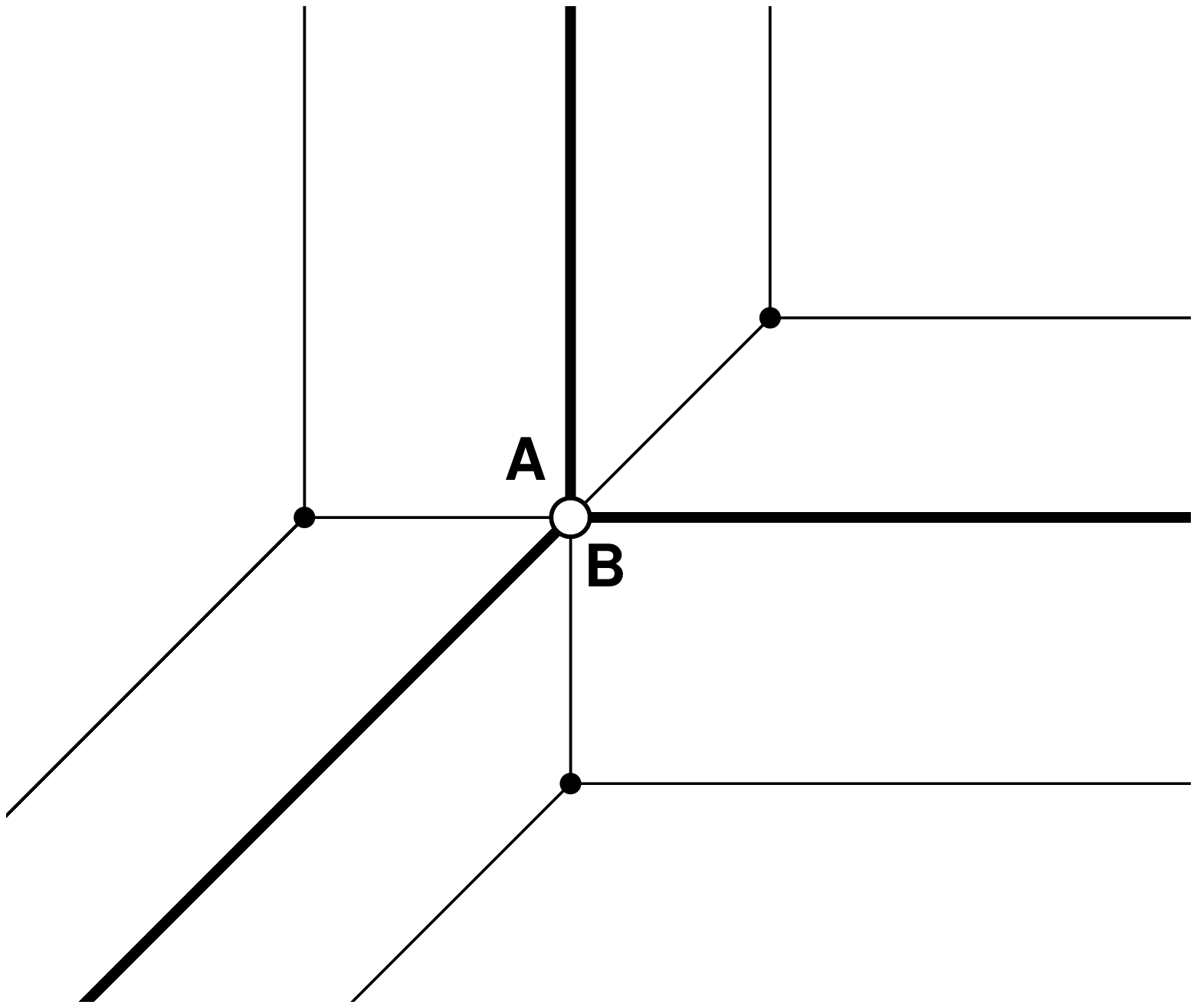}\qquad\qquad
  \includegraphics[width=4.5 cm]{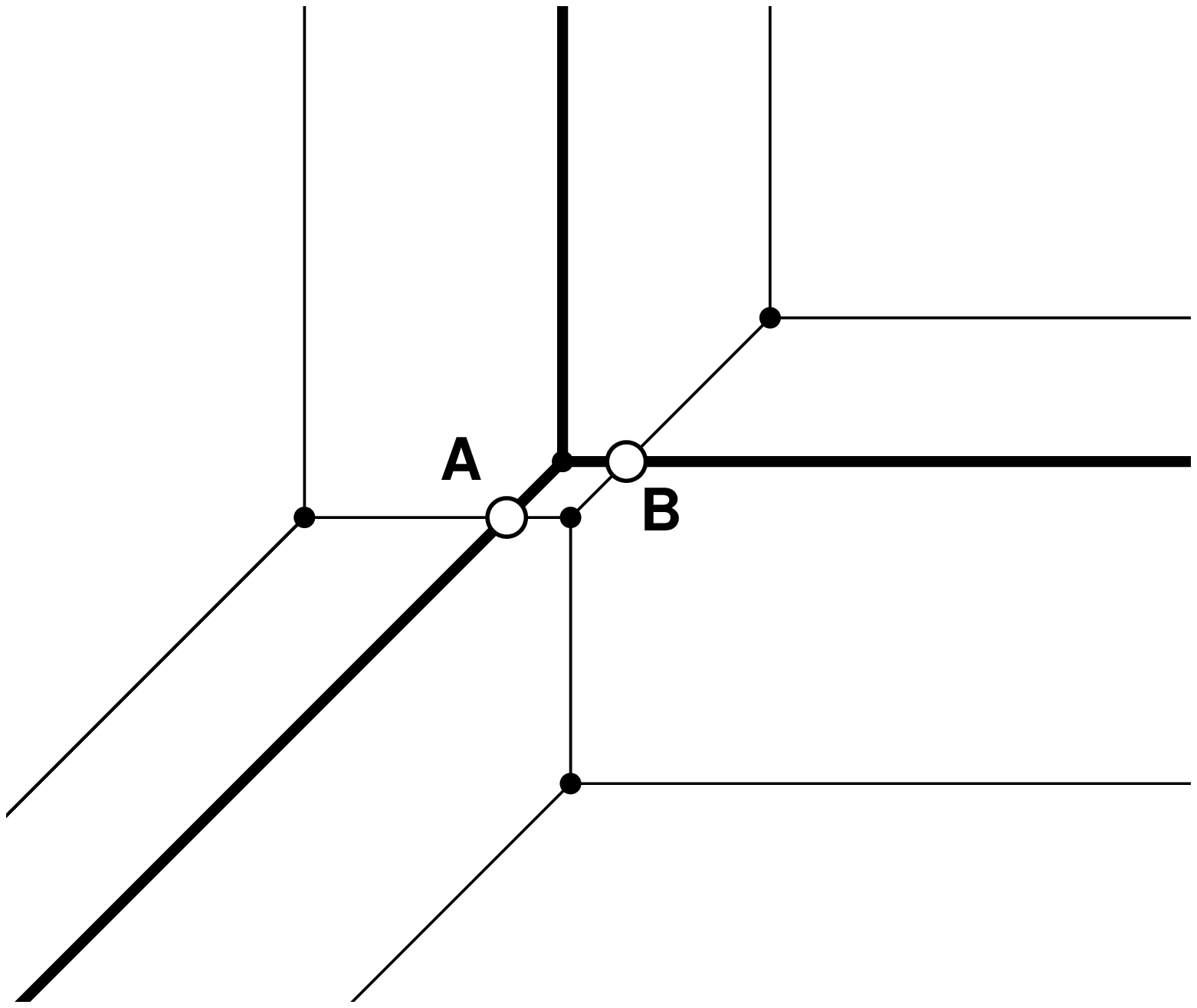}
\]

\vspace*{-.2cm}

\caption{Non-transversal intersection of a line and a conic.}
\label{fi:nontransversal}
\end{figure}
\fi

We say that the curves $C$ and $D$ intersect \emph{transversally}
if each intersection point lies in the relative interior
of an edge of $C$ and in the relative interior of an edge of $D$.
Theorem \ref{bezout} is now properly stated for the
case of transversal intersections. Figure~\ref{fi:nontransversal}
shows a non-transversal intersection of a 
 tropical conic with a tropical line. In the left picture a slight perturbation of the situation is
shown. It shows that the point of intersection {\it really} comes from two points of intersection and has to be counted 
with the multiplicity that is the sum of the two points in the nearby situation.
We will first give the proof of B\'{e}zout's Theorem for the transversal
case, and subsequently we will discuss the case of non-transversal 
intersections.

\begin{proof}
The statement holds for curves in special position for which
all intersection points occur among the half rays of
the first curve in $x$-direction and the half rays of the second curve in
$y$-direction. Such a position is shown in Figure~\ref{fi:trans}.

\ifpictures
\begin{figure}[h]
\vspace*{0cm}

\[
  \includegraphics[width=4.8 cm]{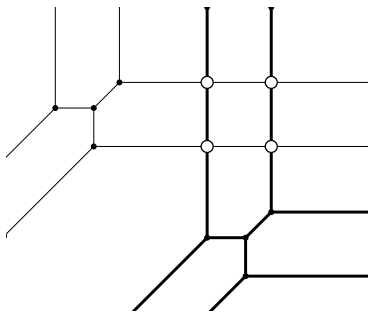}
\]

\vspace*{0cm}

\caption{Two conics intersect in four points.}
\label{fi:trans}
\end{figure}
\fi

The following homotopy moves
any instance of two transversally intersecting curves to such a special situation.
We fix the first curve $C$
and we translate the second curve $D$ with constant velocity along
a sufficiently general piecewise linear path. Let $D_t$ denote the curve $D$ 
at time $t \geq 0$.
We can assume that for no value of $t$ a vertex of $C$ coincides
with a vertex of $D_t$ and that for all but finitely many values of $t$
the two curves $C$ and $D_t$ intersect transversally.
Suppose these special values of $t$ are the time stamps 
$t_1 < t_2 < \cdots < t_r $. For any value of $t$ strictly between
two successive time stamps $t_i$ and $t_{i+1}$, the number of
intersection points in $\,C \,\cap \,D_t \,$ remains unchanged,
and so does the multiplicity of each intersection point.
We claim that the total intersection number also remains unchanged
across a time stamp $t_i$. 

Let $P$ be the set of branching points of $C$ which are also contained
in $D_{t_i}$ and the set of branching points of $D_{t_i}$ which are
also contained in $C$. Since $P$ is finite it suffices to show
the invariance of intersection multiplicity for any point $p \in P$.
Either $p$ is a vertex of $C$ and lies in the relative interior 
of a segment of $D_{t_i}$, 
or $p$ is a vertex of $D_{t_i}$ and lies in the relative 
interior of a line segment of $C$.
The two cases are symmetric, so we may assume that $p$ is a vertex
of $D_{t_i}$ and lies in the relative interior of a segment $S$ of $C$.
Let $\ell$ be the line underlying $S$ and $u$ be the weighted outgoing
direction vector of $p$ along $\ell$. Further let
$v^{(1)}, \ldots, v^{(k)}$ and $w^{(1)}, \ldots, w^{(l)}$ be the 
weighted direction vectors of the outgoing edges of $p$ into 
the two open half planes defined by $\ell$. 
At an infinitesimal time $t$ before and after $t_i$ the 
total intersection multiplicities at the neighborhoods of $p$ are
\[
  m' = 
  \sum_{i=1}^k \left|
  \det \left( \begin{array}{ccc}
    u_1 & u_2 & u_3 \\
    v^{(i)}_1 & v^{(i)}_2 & v^{(i)}_3 \\
    1 & 1 & 1
  \end{array} \right) \right| \; \text{ and } \;
  m'' = 
  \sum_{j=1}^l
  \left| \det \left( \begin{array}{ccc}
    u_1 & u_2 & u_3 \\
    w^{(j)}_1 & w^{(j)}_2 & w^{(j)}_3 \\
    1 & 1 & 1
  \end{array} \right) \right| .
\] 
Since within each of the two sums the determinants have the same
sign, equality of $m'$ and $m''$ follows immediately from the equilibrium 
condition at $p$.

In case of a non-transversal intersection, the intersection
multiplicity is the (well-defined) multiplicity of any perturbation
 in which all intersections are transversal 
(see Figure \ref{fi:nontransversal}). The validity of
this definition and the correctness of B\'ezout's theorem 
now follows from our previous proof for the transversal case.
\end{proof}

The statement of B\'ezout's Theorem is also valid for the intersection of
$n-1$ tropical hypersurfaces of degrees $d_1,d_2,\ldots,d_{n-1}$ in $\tp^{n-1}$.
If they intersect in finitely many points, then the number of these points
(counting multiplicities) is always $\,d_1 d_2 \cdots d_{n-1}$.
Moreover, also \emph{Bernstein's Theorem} for sparse
systems of polynomial equations remains valid in the
tropical setting. This theorem states that the number of
intersection points always equals the mixed volume
of the Newton polytopes. For a discussion of the
tropical Bernstein Theorem see \cite[Section 9.1]{SSPE}.

\smallskip

Families of tropical complete intersections have an important feature
which is not familiar from the classical situation, namely, intersections can be 
continued across the entire parameter space of coefficients.
We explain this for the intersection of two curves $C$ and $D$ of
degrees $c$ and $d$ in $\tp^2$. Suppose the (geometric) intersection of
$C$ and $D$ is not finite. Pick {\bf any} nearby curves $C_\epsilon$
and $D_\epsilon$ such that $C_\epsilon$ and $D_\epsilon$
intersect in finitely many points. Then
$\,C_\epsilon \, \cap \, D_\epsilon \,$ has cardinality $cd$. 

\begin{thm} \label{stable}
The limit of the point configuration $\,C_\epsilon \, \cap \, D_\epsilon \,$
is independent of the choice of perturbations. It is a well-defined
subset of $cd$ points in $\,C \,\cap \, D$.
\end{thm}

Of course, as always, we are counting multiplicities in the intersection
$\,C_\epsilon \, \cap \, D_\epsilon \,$ and hence also in its limit as $\epsilon $
tends to $ 0$. This limit is a configuration of points with multiplicities,
where the sum of all multiplicities is $cd$.
We call this limit the \emph{stable intersection} of the curves $C$ and $D$, 
and we denote this multiset of points by
$$ C \,\cap_{st} D \quad = \quad \lim_{\epsilon \rightarrow 0}
(C_\epsilon \, \cap \, D_\epsilon ). $$
Hence we can strengthen the statement of B\'ezout's Theorem as follows:

\begin{cor}
Any two curves of degrees $c$ and $d$ in the
tropical projective plane $\tp^2$ intersect stably in a
well-defined set of $cd$ points, counting multiplicities.
\end{cor}

\ifpictures
\begin{figure}[h]
\vspace*{0cm}

\[
  \includegraphics[width=4.0 cm]{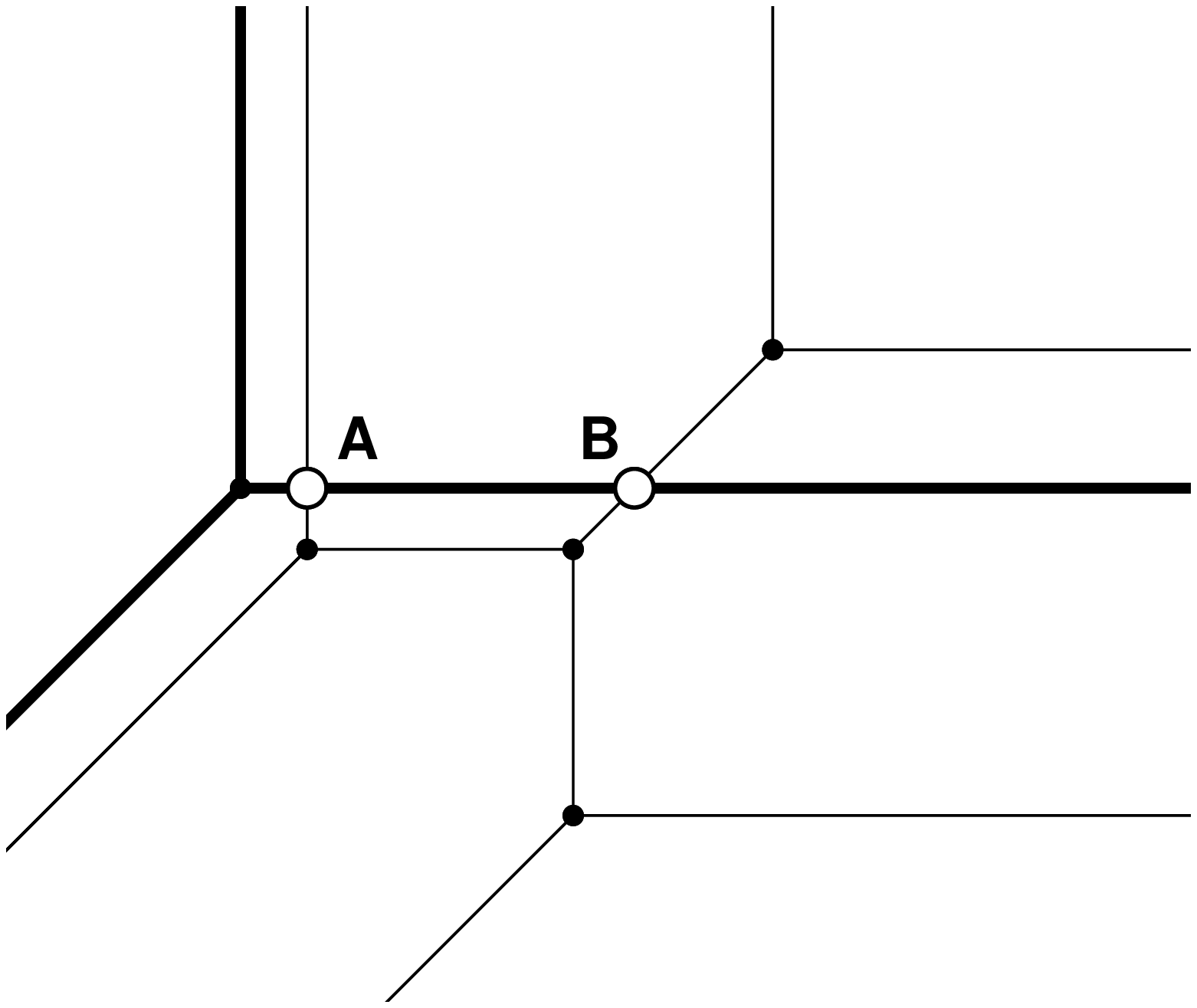}\quad
  \includegraphics[width=4.0 cm]{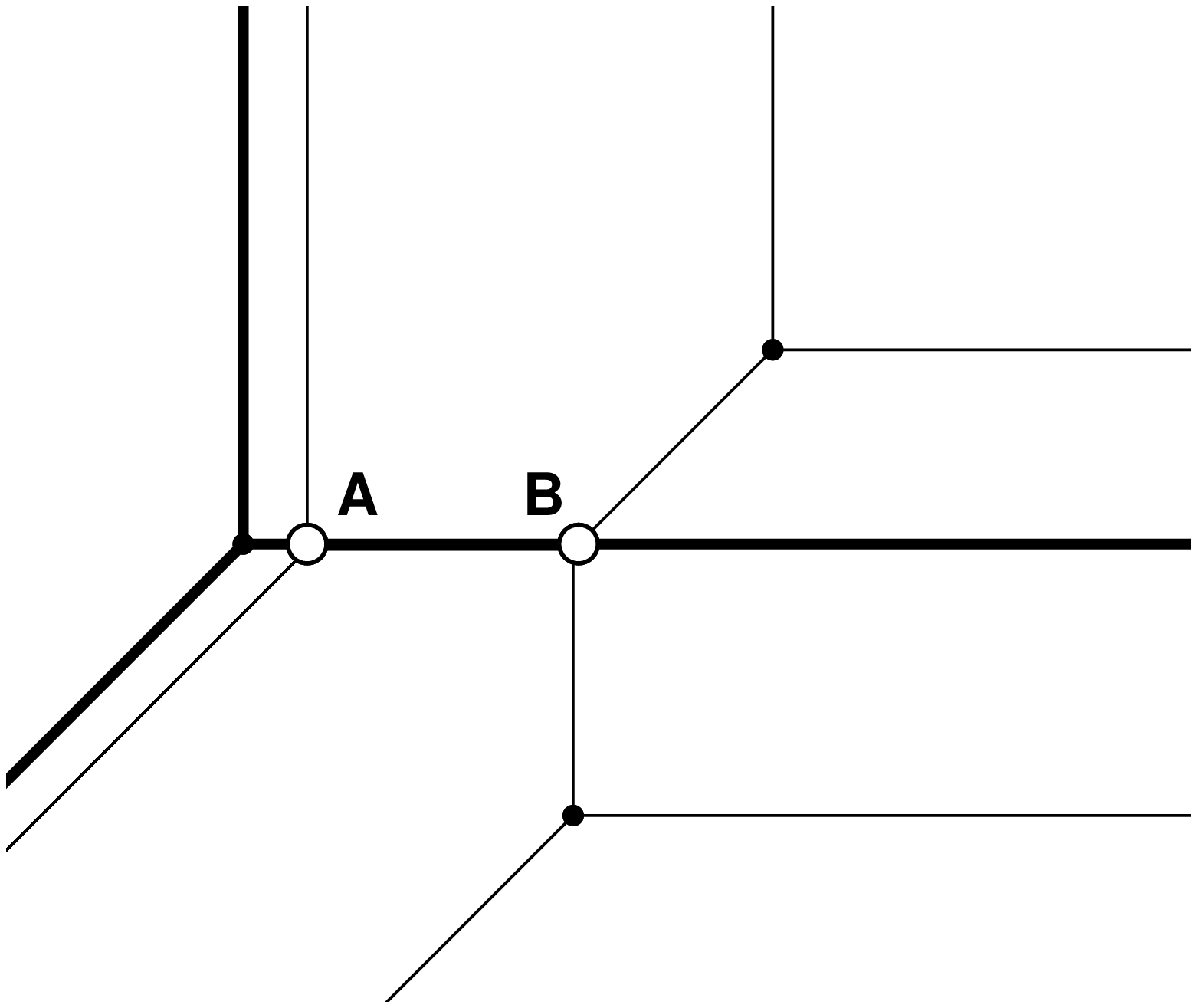}\quad
  \includegraphics[width=4.0 cm]{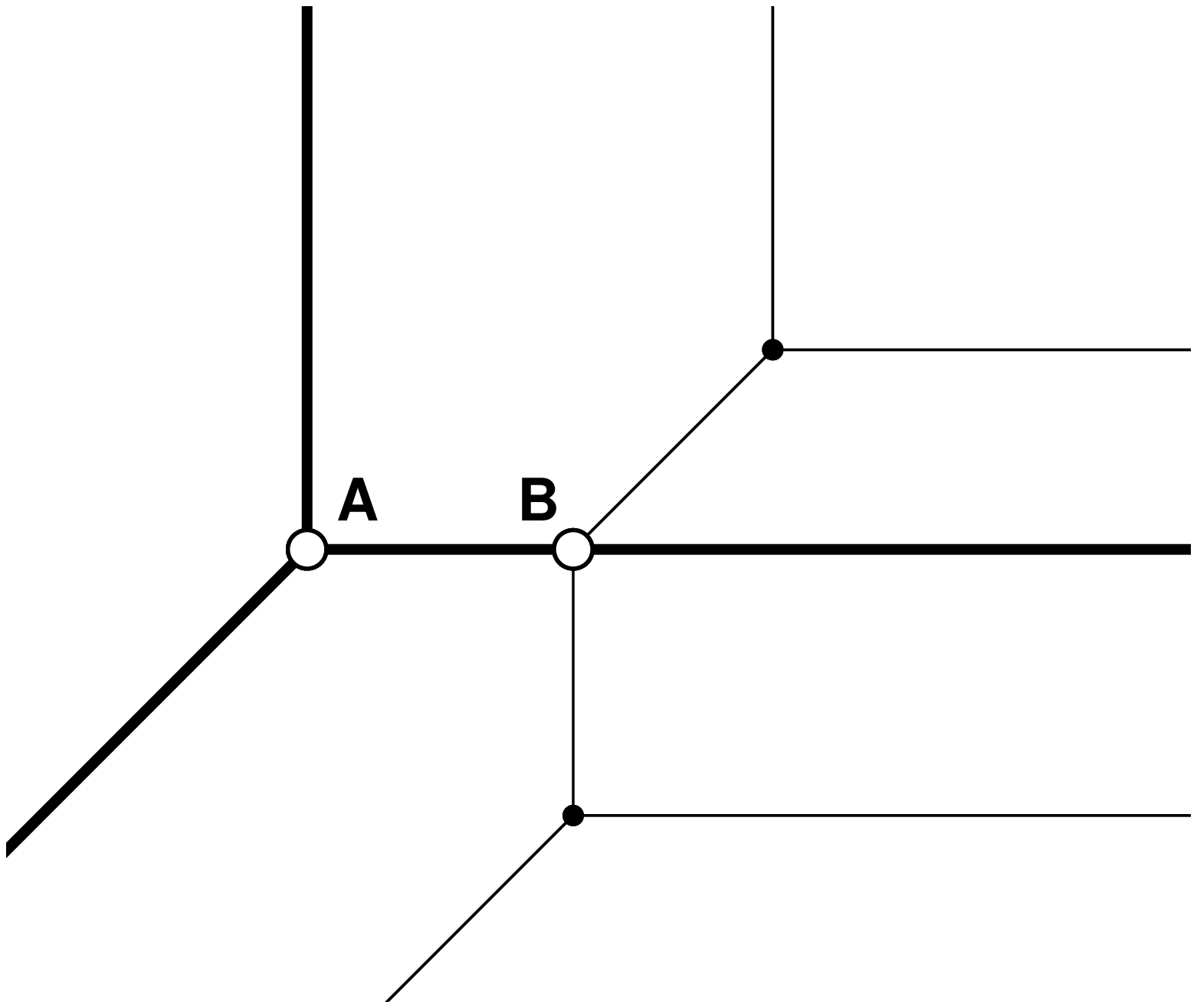}
\]

\vspace*{0cm}

\caption{Stable intersections of a line and a conic.}
\label{fi:stableline}
\end{figure}
\fi

The proof of Theorem \ref{stable} follows from our proof of
the tropical B\'ezout's Theorem.
We shall illustrate
the statement by two examples.
Figure~\ref{fi:stableline} shows the stable intersections of a line and a conic. In the first picture they
intersect transversally in the  points $A$ and $B$. In the second picture the line is moved to a position where the intersection is no longer transversal. The situation in the third picture is even more special. 
However, observe that for {\bf any} nearby transversal situation the intersection points will be close to $A$ and $B$. In all three pictures, the pair of points $A$ and $B$ is the stable intersection of the
line and the conic. In this manner we can construct a 
continuous piecewise linear map which maps any pair of conics
to their four intersection points.

\ifpictures
\begin{figure}[h]
\vspace*{0cm}

\[
  \includegraphics[width=4.0 cm]{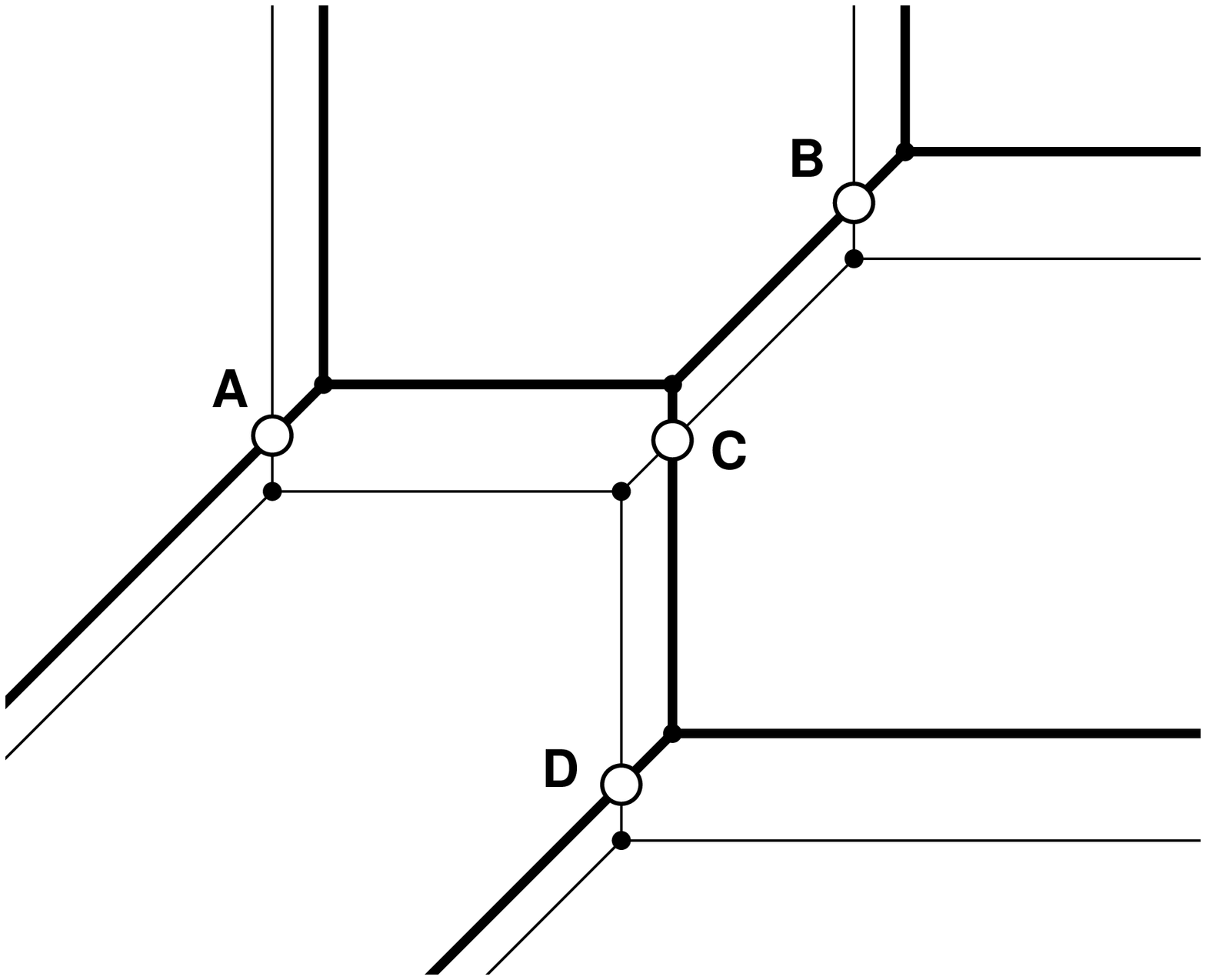}\quad
  \includegraphics[width=4.0 cm]{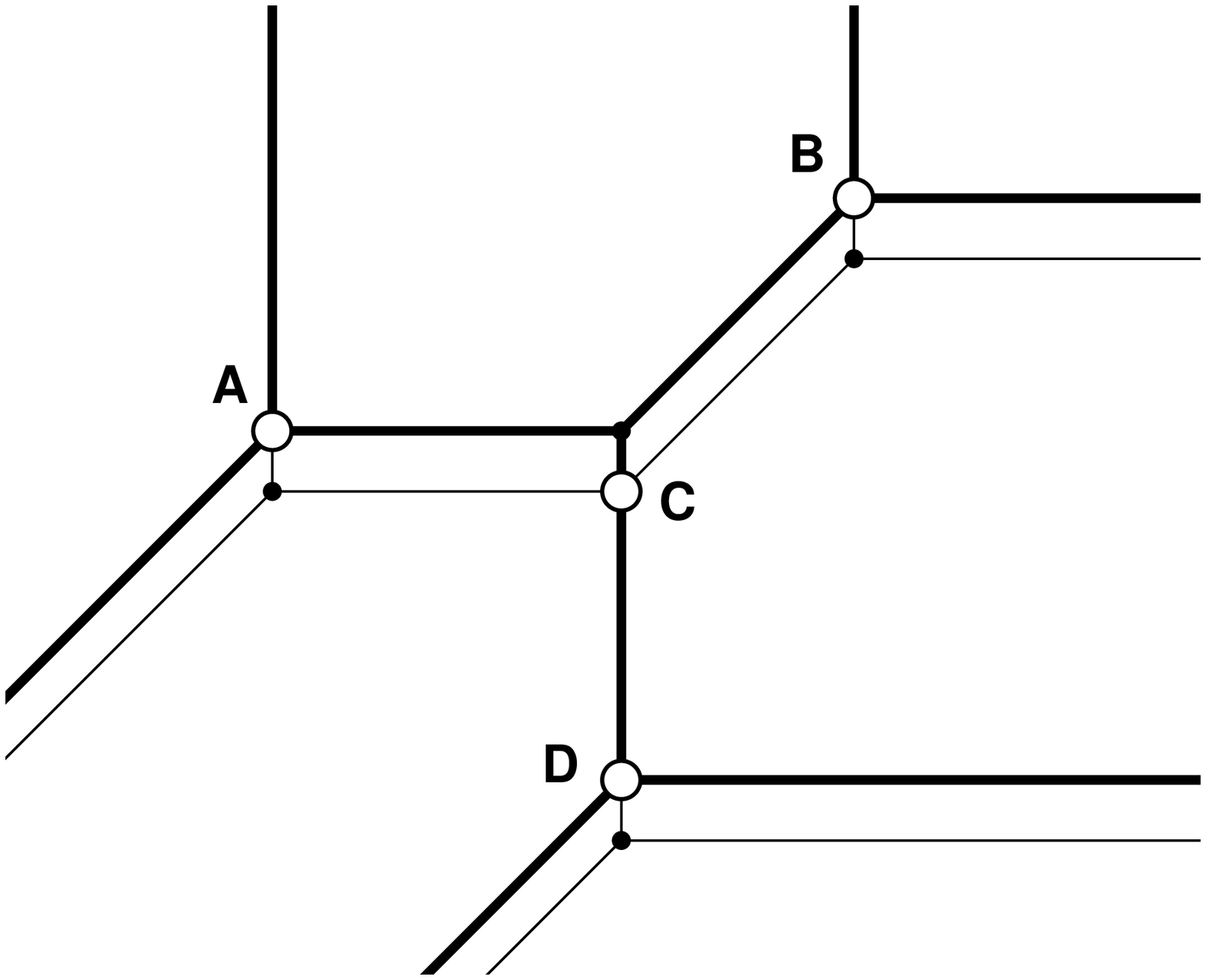}\quad
  \includegraphics[width=4.0 cm]{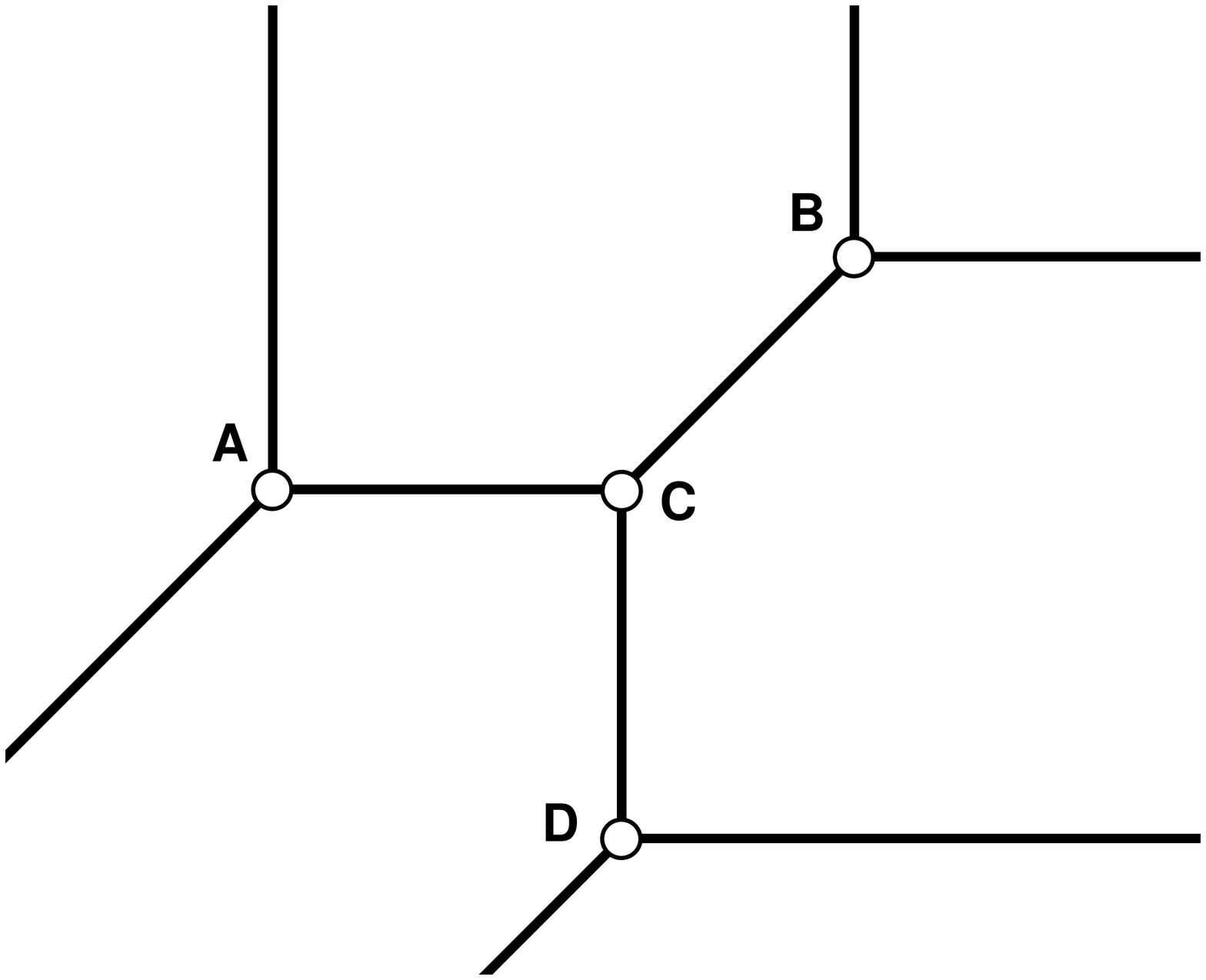}
\]

\vspace*{0cm}

\caption{The stable intersection of a conic with itself.}

\label{fi:stableconic}
\end{figure}
\fi

Figure~\ref{fi:stableconic}
illustrates another fascinating feature of stable intersections. It shows the intersection
of a conic with a translate of itself in a sequence of three pictures. The points
in the stable intersection are labeled $A,B,C,D$. Observe that in the 
third picture, where the conic is intersected with itself, the stable intersections coincide with the 
four vertices of the conic. The same works for all tropical hypersurfaces in all dimensions.
The stable self-intersection of a tropical hypersurface
 in $\tp^{n-1}$ is its set of vertices, each
counted with an appropriate multiplicity.

Our discussion suggests a general result to the effect
that {\sl there is no monodromy in tropical geometry}. It would be worthwhile
to make this statement precise.

\section{Solving linear equations using Cramer's rule\label{se:cramer}}

We now consider the problem of intersecting
$k$ tropical hyperplanes in $\tp^{n-1}$. If these
hyperplanes are in general position then their intersection
is  a tropical linear space of dimension $n-k-1$. If they
are in special position then their intersection is a tropical
prevariety of dimension larger than $n-k-1$ but it is usually not
a tropical variety. However, just like in the previous section,
this  prevariety always contains a well-defined stable intersection
which is a tropical linear space of dimension $n-k-1$.
The map which computes this stable intersection is nothing but
\emph{Cramer's Rule}. The aim of this section is to make
these statements precise and to outline their proofs.

Let $A = (a_{ij})$ be a $k \times k$-matrix with entries in $\rr \cup \{+ \infty \}$.
We define the \emph{tropical determinant} of $A$ by evaluating the
expansion formula tropically:
\[
  \text{det}_{\text{trop}}(A) \,\, = \,\,
  \bigoplus_{\sigma \in S_k} (a_{1,\sigma_1} \odot \cdots \odot a_{k,\sigma_k}) \,\, = \,\,
  \min_{\sigma \in S_k} (a_{1,\sigma_1} + \cdots + a_{k,\sigma_k}) \, .
\]
Here $S_k$ is the group of permutations of $\{1,2,\ldots,k\}$.
The matrix $A$ is said to be \emph{tropically singular} if this minimum
is attained more than once. This is equivalent to saying that $A$ is in
 the tropical variety defined by the ordinary $k \times k$-determinant.

\begin{lemma} \label{bemerkung} The matrix $A$ is tropically singular if and only if
the $k$ points whose coordinates are the column vectors
of $A$ lie on a tropical hyperplane in $\tp^{k-1}$.
\end{lemma}

\begin{proof}
If $A$ is tropically singular then
we can choose a $k \times k$-matrix $U(t)$ with entries in the field
$K = \overline{\qq(t)}$ such that 
${\rm det}(U(t)) = 0$ and ${\rm order}(U(t)) = A$.
There exists a non-zero vector $v(t) \in K^n$ in the
kernel of $\,U(t)$. The identity
$\,U(t) \cdot v(t) = 0 \,$ implies that the column vectors of $A$
lie on the tropical hyperplane 
$\,\mathcal{T}\bigl(v_1(t) \cdot x_1 + \cdots + v_n(t) \cdot x_n \bigr)$.
The converse direction follows analogously.
\end{proof}

The tropical determinant of a matrix is also known as
the \emph{min-plus permanent}. In tropical geometry,
just like in algebraic geometry over a field of characteristic~$2$,
the determinant and the permanent are indistinguishable.
In practice, one does not compute the tropical determinant of a $k \times k$-matrix
by first computing all $k \,! $ sums $\,a_{1,\sigma_1} + \cdots + a_{k,\sigma_k} \,$
and then taking the minimum. This process would take exponential time in $k$.
Instead, one recognizes this task as an
\emph{assignment problem}. A well-known result from
combinatorial optimization \cite[Corollary 17.4b]{Schrijver} implies

\begin{rmk}
The tropical determinant of a $k \times k$-matrix
can be computed in $O(k^3)$ arithmetic operations.
\end{rmk}

Fix a  $k \times n$-matrix $C = (c_{ij})$ with $k < n$.
Each row gives a tropical linear form 
$$ F_i \quad = \quad 
c_{i1}  \odot x_1 \, \oplus \, c_{i2}  \odot x_2 \, \oplus \,
\cdots \, \oplus c_{in}  \odot x_n .$$
For any $k$-subset $I = \{i_1 , \ldots , i_k \}$ of 
$\{1,2,\ldots,n\}$, we let $C_I$ denote the
$k \times k$-submatrix of $C$ with column indices $I$,
and we abbreviate its tropical determinant
$$ w_I \quad = \quad  \text{det}_{\text{trop}}(C_I) . $$
Let $J = \{  j_0 , j_1 ,\cdots , j_{n-k} \}$ be any $(n-k+1)$-subset
of $\{1,\ldots,n\}$ and $J^c$ its complement.
 The following tropical linear form
is a \emph{circuit} of the $k \times n$-matrix $C$:
$$ G_J \quad = \quad 
w_{J^c \cup \{j_0\}} \odot x_{j_0} \,\,\,\oplus \,\,\,
w_{J^c \cup \{j_1\}} \odot x_{j_1} \,\,\,\oplus \,\,\, \cdots \,\,\,\oplus \,\,\,
w_{J^c \cup \{j_{n-k}\}} \odot x_{j_{n-k}} . $$
We form the intersection of the $\binom{n}{k-1}$ tropical hyperplanes
defined by the circuits:
\begin{equation} \label{stableplane}
\bigcap_{1 \leq  j_0 < j_1 < \cdots < j_k \leq n } \,\,\,
\mathcal{T}\bigl(G_J \bigr) \qquad \subset \quad \tp^{n-1}.
\end{equation}
Similarly, we form the intersection of the $k$ tropical hyperplanes
given by the rows:
\begin{equation} \label{preplane}
\mathcal{T}(F_1) \,\,\cap \,\,
\mathcal{T}(F_2) \,\,\cap \,\,\cdots \,\,\cap \,\,\mathcal{T}(F_k)
 \qquad \subset \quad \tp^{n-1}.
 \end{equation}
 The problem of solving a system of $k$ tropical  linear equations in $n$ unknowns
 is the same as computing the intersection (\ref{preplane}). We saw
  in Section 3 that (\ref{preplane}) is in general only a 
prevariety, even for $k = 2$, $n=3$.
 In higher dimensions this prevariety can have maximal faces of different dimensions.
 The intersection (\ref{stableplane}) is much nicer. We show that it is the stable
 version of the poorly behaved intersection (\ref{preplane}):
  
 \begin{thm} \label{stablecramer}
 The intersection (\ref{stableplane}) is a tropical linear space
 of codimension $k$ in $\tp^{n-1}$. It is always contained in
 the prevariety (\ref{preplane}). The two intersections are equal if
 and only if none of the $k \times k$-submatrices $\,C_I\,$ of $\,C\,$
 is tropically singular.
 \end{thm}
 
  \medskip
 
 \noindent {\sc Proof. }
 Consider the vector $w \in \rr^{\binom{n}{k}}$ whose coordinates
 are the tropical $k \times k$-determinants $w_I$. Then $w$ is a point
 on the tropical Grassmannian  \cite{speyer-sturmfels-2003}.
 By \cite[Theorem 3.8]{speyer-sturmfels-2003},
the point $w$ represents a tropical linear space $L_w$ in $\tp^{n-1}$.
 The tropical linear space $L_w$ has codimension $k$, and it is precisely the
 set (\ref{stableplane}). 
 
  The second assertion follows from the if-direction in the third assertion.
 Indeed, if $C$ is any $k \times n$-matrix which has a tropically
 singular $k \times k$-submatrix then we can find a family of matrices
 $C^{(\epsilon)}$, $\epsilon > 0$, with $\, 
 {\rm lim}_{\epsilon \rightarrow 0} C^{(\epsilon)} \, = \, C \,$ such that
 each $C^{(\epsilon)}$ has all $k \times k$-matrices tropically
 non-singular. Let $w^{(\epsilon)} \in \rr^{\binom{n}{k}}$ be the vector
 of tropical $k \times k$-subdeterminants of $C^{(\epsilon)} $.
 The tropical linear space $L_{w^{(\epsilon)}}$ depends continuously on 
 the parameter $\epsilon$. If $x$ is any point in the plane $L_w$ 
 then there exists a sequence of
 points $x^{(\epsilon)} \in L_{w^{(\epsilon)}}$ such that
  $\,{\rm lim}_{\epsilon \rightarrow 0} x^{(\epsilon)} \, = \, x$.
  Now consider the tropical linear form $\,F_i^{(\epsilon)} \,$ given by
  the $i$-th row of the matrix $C_\epsilon$.
 By the if-direction of the third assertion, the point $x^{(\epsilon)}$ lies in
 the tropical hypersurface $\mathcal{T}(F_i^{(\epsilon)})$. This
 hypersurface also depends continuously on $\epsilon$, and therefore
 we get the desired conclusion $x \in \mathcal{T}(F_i)$ for all $i$.
 
 We next prove the  only if-direction in the third assertion.
Suppose that the $k \times k$-matrix $C_I$
is tropically singular. By Lemma \ref{bemerkung},
there exists a vector $p \in \rr^k$ in the
tropical kernel of $C_I$. We can augment
$p$ to a vector in the
tropical kernel (\ref{preplane}) of $C$
 by placing any sufficiently large positive
 reals in the other $n-k$ coordinates. 
  Hence the prevariety (\ref{preplane})
 contains a polyhedral cone of dimension
 $n-k$ in $\tp^{n-1}$.
 
 What is left to prove at this point is the
 most difficult part of Theorem \ref {stablecramer},
 namely,  assuming that the 
 $k \times k$-submatrices of $C$ are tropically
 non-singular then we must show that the tropical prevariety
 (\ref{preplane}) is contained in the tropical linear space
 (\ref{stableplane}). We will now interrupt the general proof
 to give a detailed discussion of the
 special case $n=k+1$,  when 
 (\ref{stableplane}) and (\ref{preplane}) 
  consist of  a single point in $\tp^{n-1}$.
 Thereafter we shall return to the third assertion in
 Theorem \ref {stablecramer} for $n \geq k+2$.
 \qed
 
 \medskip

Let $C$ be a rational $(n-1) \times n$-matrix whose
$n$ maximal square submatrices are tropically
non-singular. We claim that
the associated system of
$n-1$ tropical linear equations in $n$ unknowns
has a unique solution point $p$
 in tropical projective space $\tp^{n-1}$.
 In the next two paragraphs, we shall explain how
 the results on \emph{linkage trees} 
 in \cite[Theorem 2.4]{zelevinsky}
 can be used to prove both this result and to 
devise a polynomial-time algorithm for computing the
point $p$ from the matrix $C$.

Let $Y$ be an $(n-1) \times n$-matrix of indeterminates,
and consider the problem of minimizing the dot product
$C \cdot Y = \sum_{ij} c_{ij} y_{ij} $ subject to the
constraints that $Y$ is non-negative, all its
row sums are $n$ and all its column sums are $n-1$.
This is a \emph{transportation problem} which can be
solved in polynomial time in the binary encoding of $C$.
Our hypothesis that the maximal
minors of $C$ are tropically non-singular means that
$C$ specifies a \emph{coherent matching field}.
Theorem 2.8 in \cite{zelevinsky} implies that
the above transportation problem has a unique optimal solution
$Y^*$.  Each row of $Y^*$ has precisely two non-zero
entries: this data specifies the \emph{linkage tree} $T$
as in Theorem 2.4 of \cite{zelevinsky}. Thus $T$ is a tree
on the nodes $1,2,\ldots,n$ whose edges are labeled by
$1,2,\ldots,n-1$. If the two non-zero entries of the
$i$-th row of $Y^*$ are in columns $j_i$ and $k_i$ then
$\{j_i,k_i\}$ is the edge labeled by $i$.
We claim that
our desired point $p = (p_1,\ldots,p_n)$ satisfies the equations
\begin{equation}
\label{reallyeasy}
 c_{i j_i} + p_{j_i}   \,\, = \,\, c_{i k_i} + p_{k_i}  \qquad i = 1,2,\ldots,n-1. 
\end{equation}
Since $T$ is a tree, these equations have a unique solution $p$
which can be computed in polynomial time from $Y^*$ and hence from $C$.
It remains to prove the claim that $p$ solves the tropical
equations and is unique with this property.

We may assume without loss of generality
that the zero vector $(0,0,\ldots,0)$ is a solution
of the given tropical equations. We can change the given matrix $C$
by scalar addition in each column until each column is
non-negative and has at least one zero entry. Then 
each row of $C$ is non-negative and has at least two zero entries.
The objective function value of our transportation problem is zero,
and the set of zero entries of $C$ supports a unique linkage tree $T$.
Now our tropical linear system consists of the $n-1$ equations
$$  p_{j_i}   \,\, = \,\,  p_{k_i}  \qquad i = 1,2,\ldots,n-1 $$
and $(n-1)(n-2)$ inequalities  $ \,p_{j_i} \leq c_{ik} $,
one for each position $(i,k)$ which is not in the matching field.
The zero vector $p = (0,0,\ldots,0)$ is the unique
solution to these equations and inequalities.
This argument completes the proof of Theorem \ref {stablecramer} 
for the case $n= k+1$. We summarize our discussion as follows.

\begin{cor} Solving a system of $n \! - \! 1$ tropical linear
equations in $n$ unknowns amounts to computing the
linkage tree of the coefficient matrix $C$. 
This can be done in polynomial time by solving the
transportation problem with cost matrix $C$, 
row sums $n$ and column sums $n-1$.
The solution $p$ is given by the system (\ref{reallyeasy}).
\end{cor}

We note that the solution to the inhomogeneous
square system of tropical linear equations
$A \odot x = b$ presented in \cite{olsder-roos-88}
corresponds to a special linkage tree. This tree is
a star with center indexed by $b$ and leaves
indexed by the $n$ columns of $A$.

 \medskip
 
 \noindent {\sc Proof of Theorem \ref {stablecramer} (continued). }
 We now sketch the proof of the if-direction in the third assertion
for $n \geq k+2 $.
 Suppose $C$ has no tropically singular $k \times k$-submatrix.
 We wish to show that if a point $x$ lies in (\ref{preplane}) then it
 also lies in~(\ref{stableplane}).  Clearly, it suffices to show
 this for the zero vector $x = 0$. Thus we will show:
 If $0$ is in  (\ref{preplane}) then it is in 
(\ref{stableplane}). After tropical scalar multiplication, we
may assume that the coefficients $c_{i1}, \ldots,c_{in}$
of $F_i$ are non-negative and their minimum is~$0$.
Since $0 \in \mathcal{T}(F)$, this minimum is attained twice.
At this stage, the $k \times n$-matrix 
 $C$ is non-negative and each row contains
at least two zero entries. 

We next apply the combinatorial theory developed in
\cite{zelevinsky}. Our  hypothesis states that the minimum in
the definition of a tropical $k \times k$-subdeterminant
of $C$ is uniquely attained.  Thus it specifies a \emph{coherent matching field}.
Let $\Sigma $ be
the \emph{support set} of this matching field. It contains
all the locations $(i,j)$ of zero entries of $C$.  This
means that the prevariety (\ref{preplane}) remains
unchanged if we replace all entries outside of $\Sigma$
by $+ \infty$. Now, using the results
in \cite[\S 4]{zelevinsky}, we can transform our matrix $C$ by
tropical row operations to an equivalent matrix $C'$
whose rows are tropical $k \times k$-minors. The
tropical linear forms $F'_i$ given by the rows of that matrix $C'$
are circuits. But they still have the same intersection 
(\ref{preplane}). From the Support Theorem in \cite{zelevinsky} we infer that
this intersection has codimension $k$ locally around $x= 0$.
This holds for any point $x$ in the relative interior of
any maximal face of  the polyhedral complex (\ref{preplane}),
and therefore the complexes (\ref{preplane}) and
(\ref{stableplane}) are equal. \qed
  
  \bigskip

 We apply the Theorem \ref {stablecramer} (with $n=6$) to
 study families of conics in the tropical projective
 plane $\tp^2$.  The support set for conics is
 $$ \mathcal{A} \quad = \quad
 \bigl\{ \,(2,0,0), 
 \, (1,1,0),
 \, (0,2,0),
 \, (0,1,1),
 \, (0,0,2),
 \, (1,0,1) \,\bigr\} . $$
 We identify points $\,a = (a_1,a_2,a_3,a_4,a_5,a_6) \,$
 with tropical conics $\mathcal{T}(F)$ as in Example \ref{ex:conictypes}.
 Fix a configuration of points  $P_i = (x_i,y_i,z_i) \in \tp^2$
 for $i = 1,2,\ldots,k$.
 Let $C$ be the $k \times 6$-matrix whose
 row vectors are
 \[
  (2 x_i, \, x_i + y_i, \, 2 y_i, \, y_i + z_i, \, 2 z_i, \, x_i + z_i) 
  \, , \quad 
  1 \le i \le k \, .
\]

\begin{lemma}
The vector $a$ lies in the tropical kernel of the matrix $C$
if and only if the points $P_1,P_2,\ldots,P_k$ lie on the
 conics $\mathcal{T}(F)$.
 \end{lemma}
  
  The implications of Theorem 
\ref{stablecramer} for $k= 4$ constitute the subject of the next section.
For $k = 5$ we conclude the following results.

\begin{cor}
\label{co:conic5points}
Any five  points in $\tp^2$  lie on
a conic. The conic is unique if and only if the points are not on
a curve whose support is a proper subset of  $\mathcal{A}$.
\end{cor}

\ifpictures
\begin{figure}[h]
\vspace*{0cm}

\[
  \includegraphics[width=6cm]{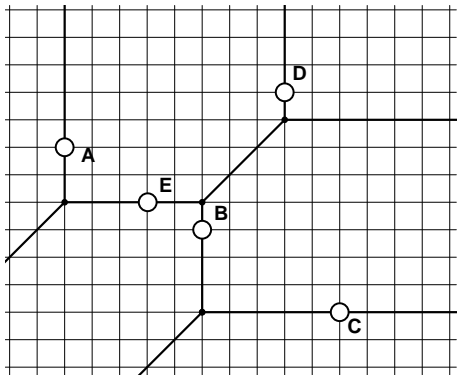}
\]

\vspace*{0cm}

\caption{Conic through five points.}
\label{fi:conic5points}
\end{figure}
\fi

The unique conic through five general points
is computed by Cramer's rule from the matrix $C$,
that is, the coefficient $a_i$ is the tropical
$5 \times 5$-determinant gotten from $C$ by deleting the $i$-th column.
The remarkable fact is that the conic
given by Cramer's rule is stable, i.e., the unique conic
through {\bf any} perturbation of the given points will
converge to that conic, in analogy to the discussion at
the end of Section 4.  It turns out that the stable conics
given by quintuples of points in $\tp^2$
are always proper, that is, they belong to types (d) and (e) in
 Example~\ref{ex:conictypes}.

\begin{thm}
The unique stable conic through any five given (not necessarily distinct)
points is proper.
\end{thm}

\begin{proof}
For $\tau \in \{x^2,xy,y^2,xz,yz,z^2\}$ let $a_\tau$ denote the tropical 
$5 \times 5$-sub\-determinant of the $5 \times 6$ coefficient matrix
obtained by omitting the column associated with the term $\tau$.
Since for any $x$ and $y$ we have
\[
  2 \cdot \min \{ 2x, 2y \} \ \le \ \min \{2x, x+y\} + \min \{x+y,2y\} \, ,
\]
the definition of the tropical determinant implies
\[
  2 a_{xy} \le a_{x^2} + a_{y^2} \, ,
\] 
and similarly
\[
  2 a_{xz} \le a_{x^2} + a_{z^2} \, , \quad
  2 a_{yz} \le a_{y^2} + a_{z^2} \, .\]
This shows that the conic with these coefficients is proper. This proves our claim.
\end{proof}

We have implemented the computation of the stable conic 
through five given points into the geometry software {\tt Cinderella}.
The user clicks any five points with the mouse onto the computer screen.
The program then sets up the corresponding $5 \times 6$-matrix $C$,
it computes the six tropical $5 \times 5$-minors of  $C$, and
it then draws the curve onto the
screen. This is done with the lifting method shown in 
Figure \ref{fi:3dattainedtwice}.

The presence of stability (i.e. the absence of monodromy) ensures
that the software behaves smoothly and always produces the
correct picture of the
 stable proper conic. For instance, if the user inputs the
 same point five times then {\tt Cinderella} will correctly draw the
 double line (from case (a) of Example~\ref{ex:conictypes})
 with vertex at that point.
 Can you guess what happens if the user clicks one point twice
 and another point three times ?

\section{Quadratic curves through four given points}

In this section we study the set of tropical conics passing through
four given points $P_i = (x_i,y_i,z_i) \in \tp^2$, $1 \le i \le 4$.
By Theorem~\ref{stablecramer}, if all $4\times 4$-submatrices of the 
$4 \times 6$-matrix $C$ with rows
\[
  (2 x_i, \, x_i + y_i, \, 2 y_i, \, y_i + z_i, \, 2 z_i, \, x_i + z_i) 
  \, , \quad 
  1 \le i \le 4 \, ,
\]
are tropically nonsingular, then this set of conics is a tropical line in
$\tp^5$. By Example~\ref{linesinPn} these are trees with six leaves.
In the following we only consider quadruples of points
which satisfy this genericity condition. We note that the
pencil of conics contains several distinguished conics.

\begin{enumerate}
\item Three degenerate conics which are the pairs of two lines 
(see Figure~\ref{fi:uniontwolines}).

\ifpictures
\begin{figure}[t]
\vspace*{0cm}

\[
  \includegraphics[width=4cm]{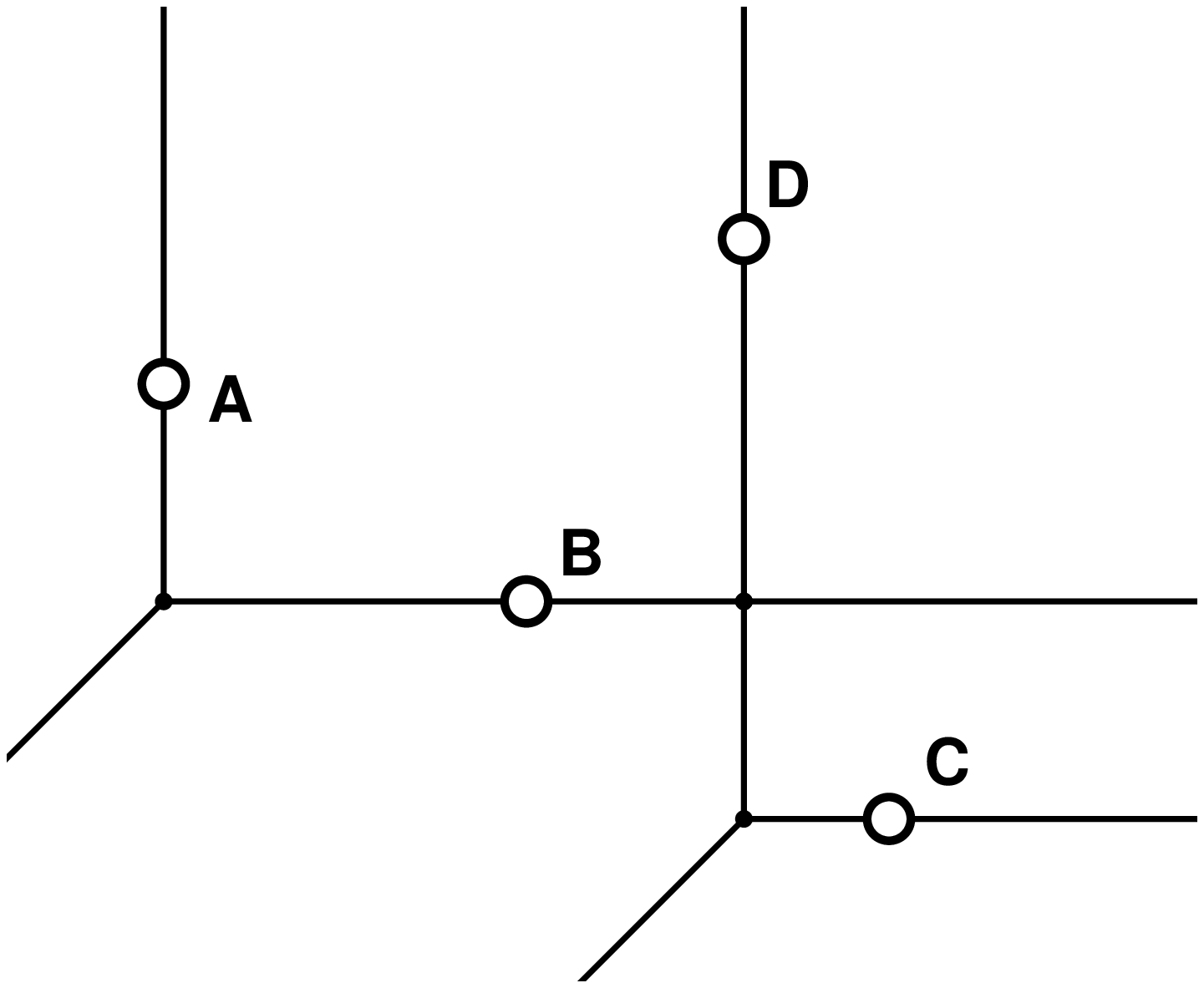} \quad
  \includegraphics[width=4cm]{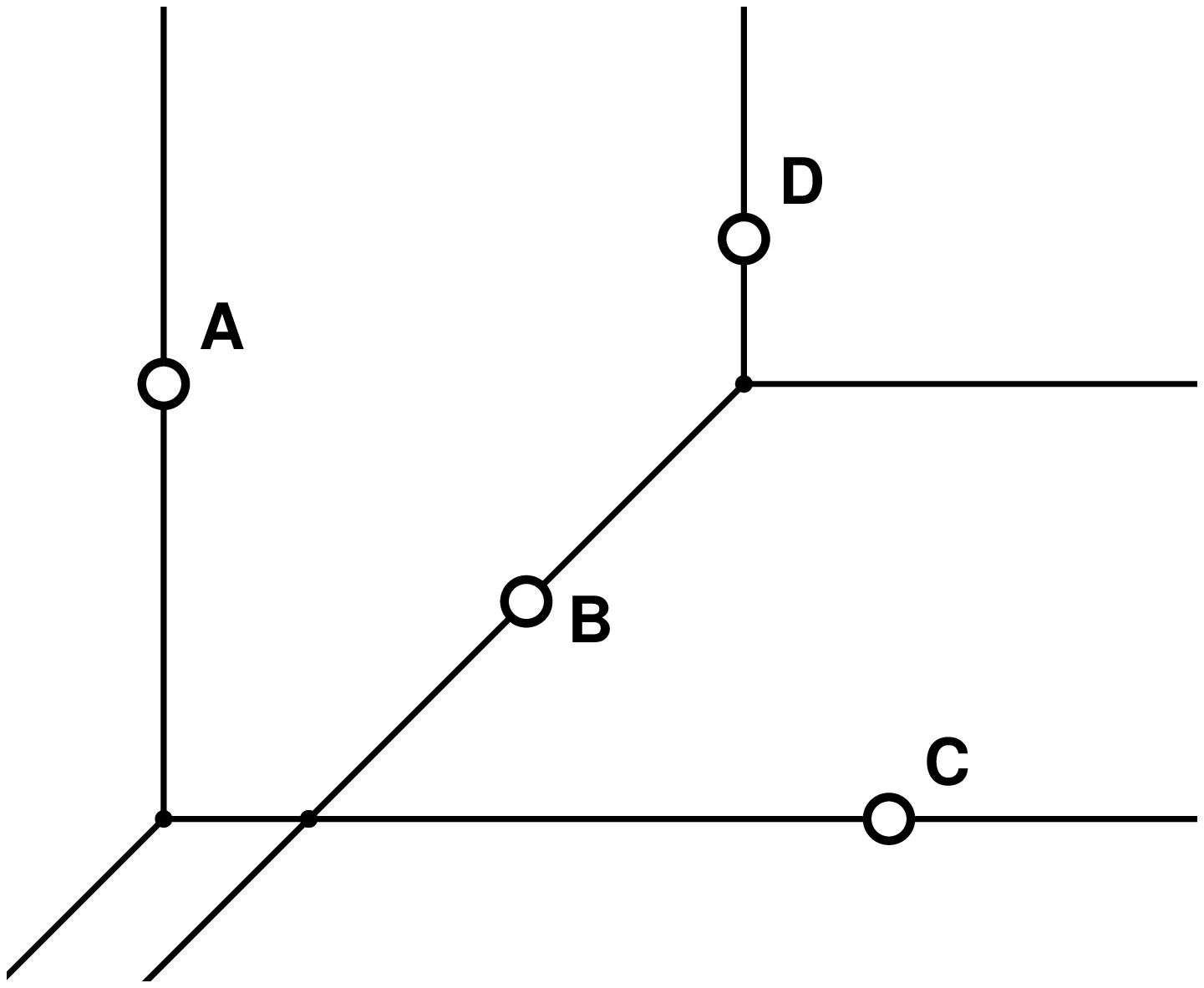} \quad
  \includegraphics[width=4cm]{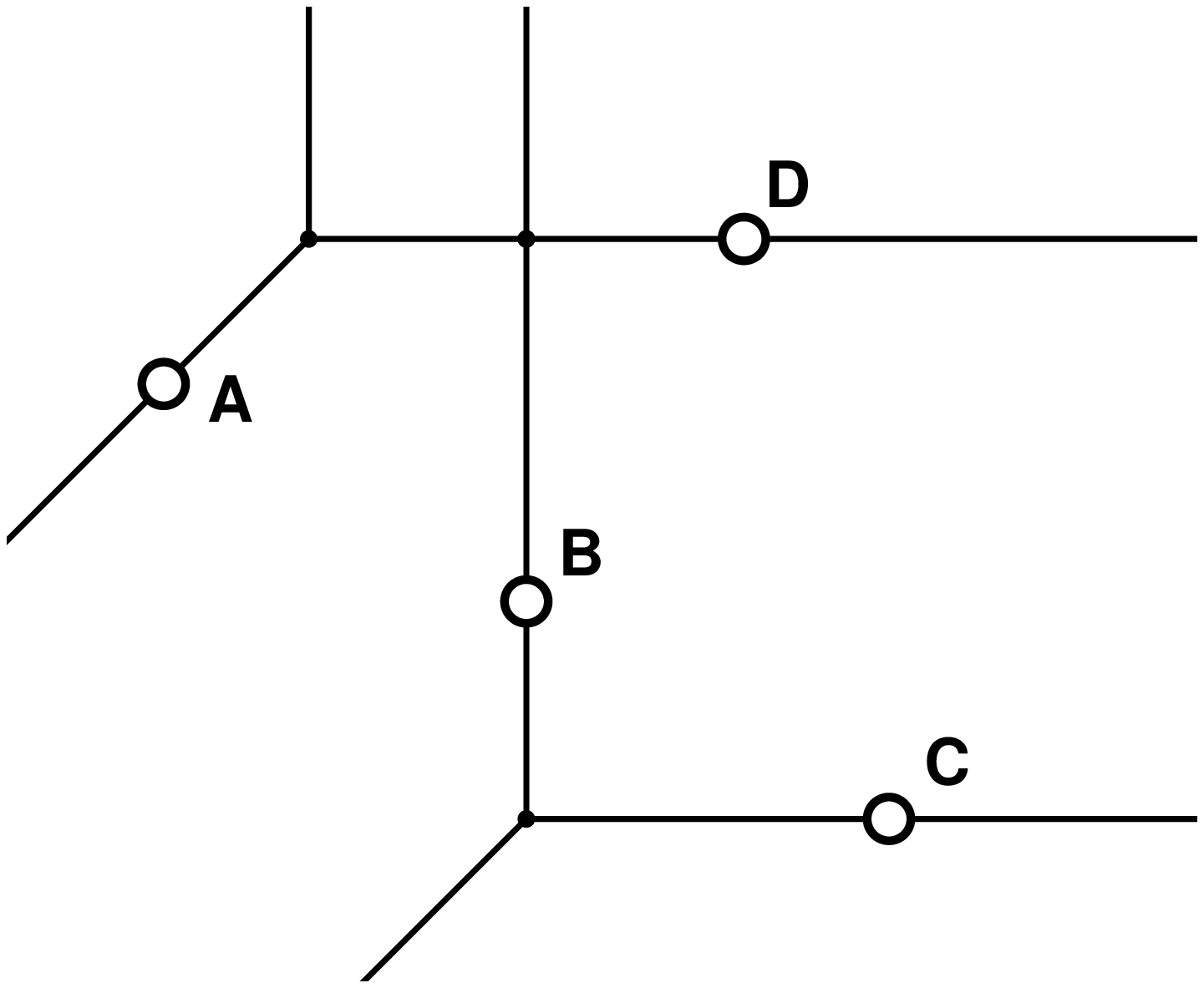}
\]

\vspace*{0cm}

\caption{Pairs of two lines through four given points}
\label{fi:uniontwolines}
\end{figure}
\fi

\item Conics in which one of the four points is a vertex of the conic
(see Figure~\ref{fi:vertexconic}).

\ifpictures
\begin{figure}[t]
\vspace*{0cm}

\[
  \begin{array}{c@{\quad}c@{\quad}c@{\quad}c}
  \includegraphics[width=4cm]{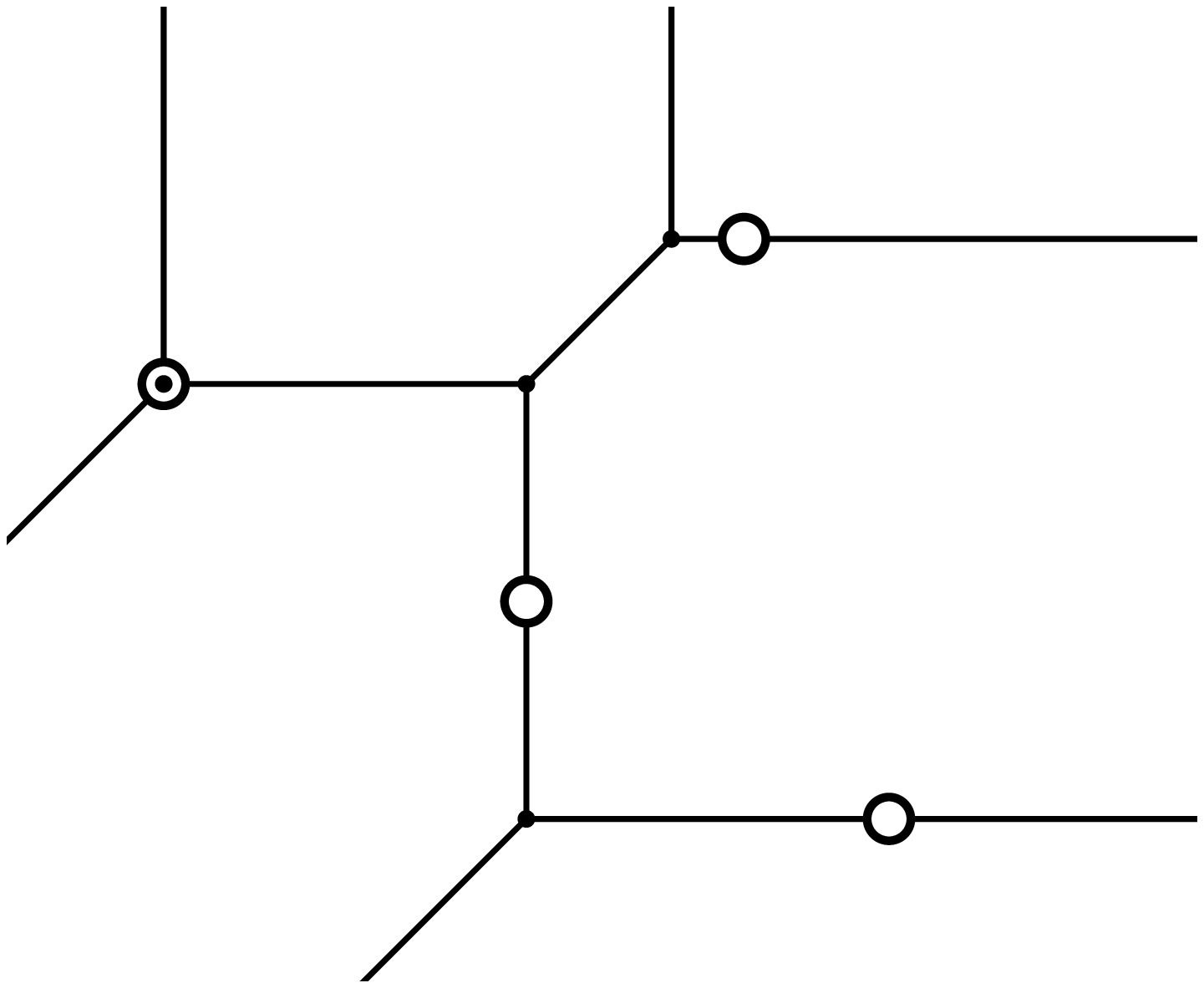} &
  \includegraphics[width=4cm]{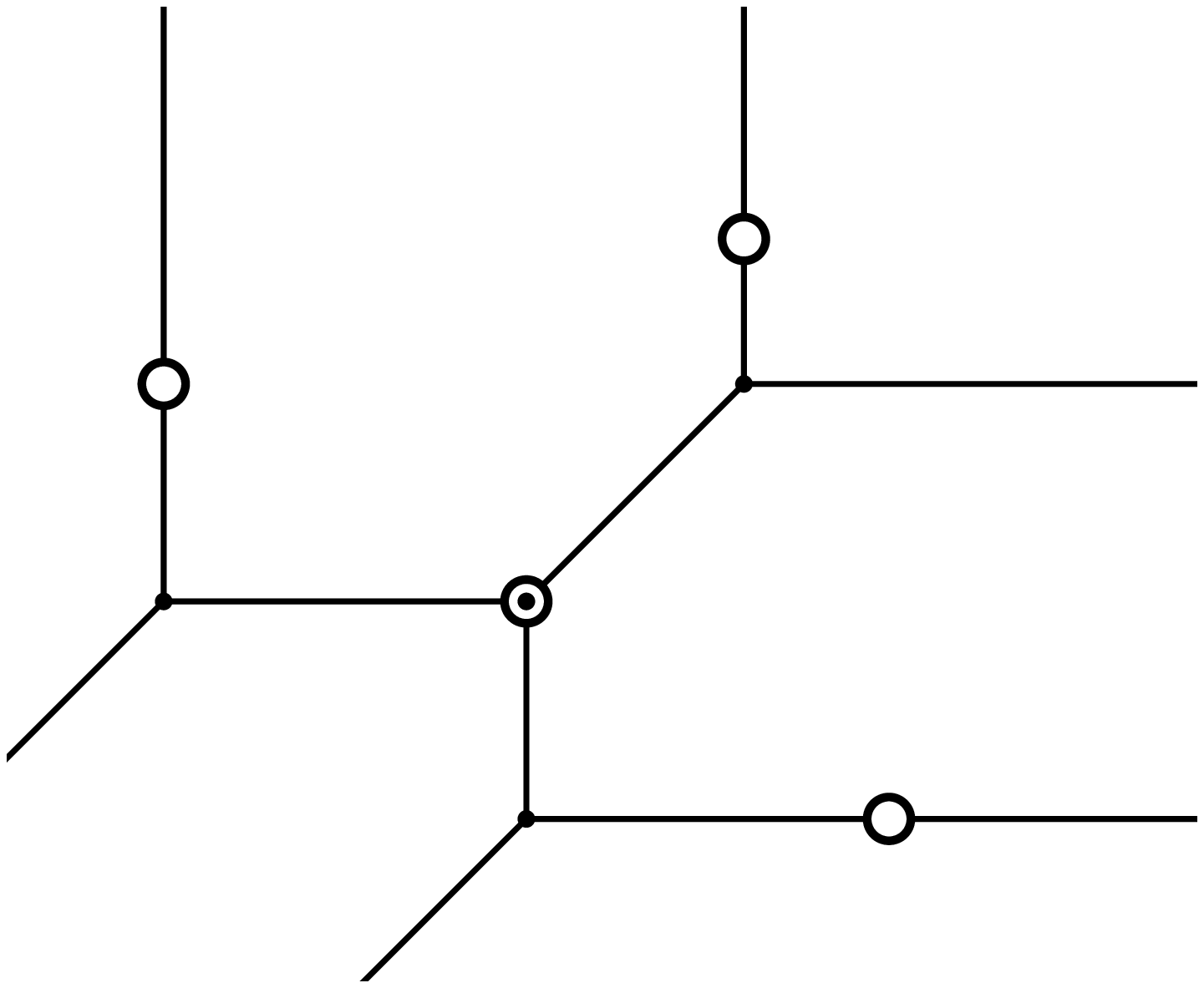} \\
  \includegraphics[width=4cm]{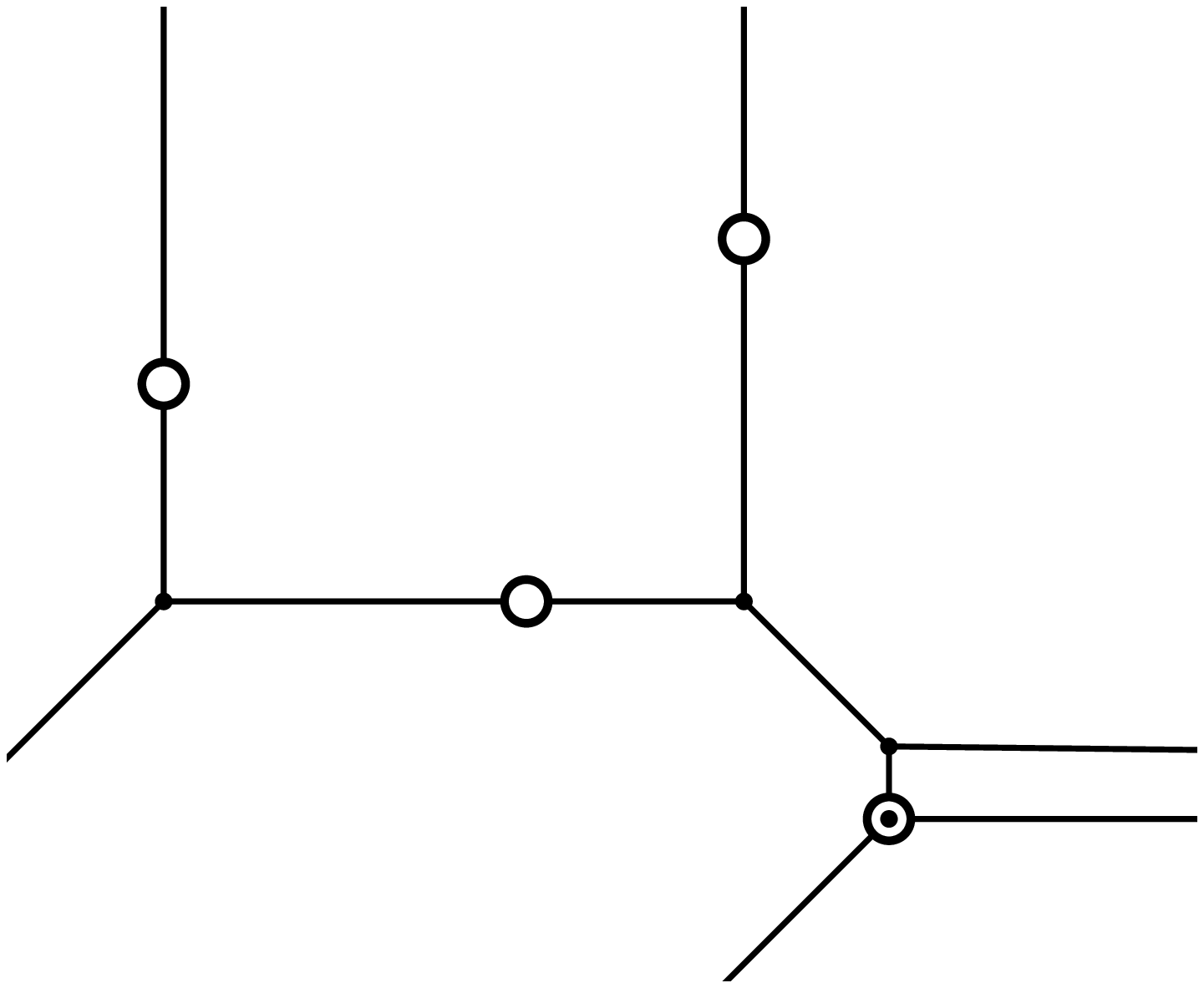} &
  \includegraphics[width=4cm]{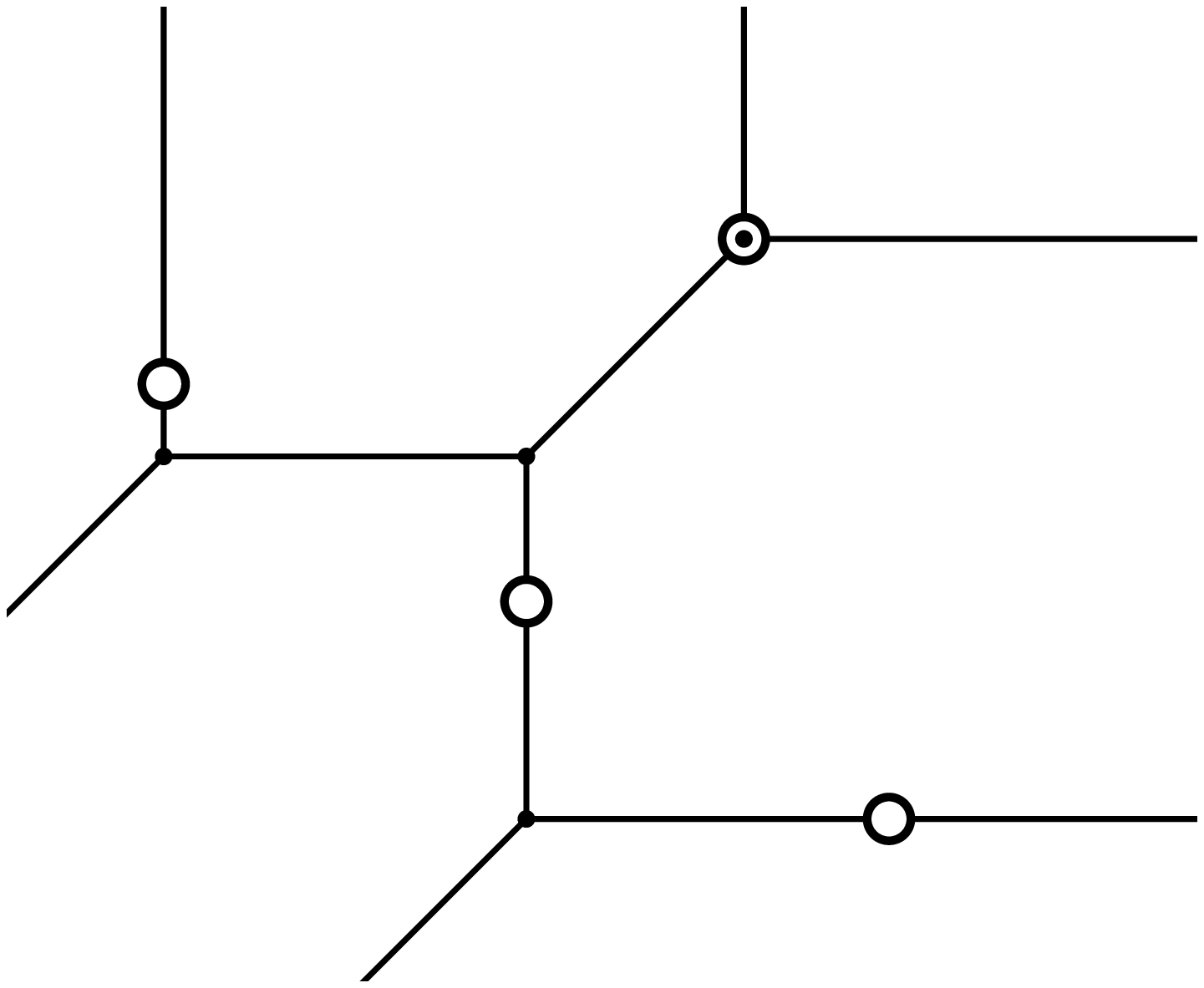}
  \end{array} 
\]

\vspace*{0cm}

\caption{One of the given points is a vertex of a conic.}
\label{fi:vertexconic}
\end{figure}
\fi

\item Conics in which the coefficient of one term converges to $+\infty$.
Figure~\ref{fi:extremeconics} shows the limits of these conics.
Geometrically, the tentacle of the conic associated with the distinguished
term is missing. Note that conics of that type have only three internal vertices.

\ifpictures
\begin{figure}[t]
\vspace*{0cm}

\[
  \begin{array}{c@{\qquad}c@{\qquad}c}
  \includegraphics[width=3.3cm]{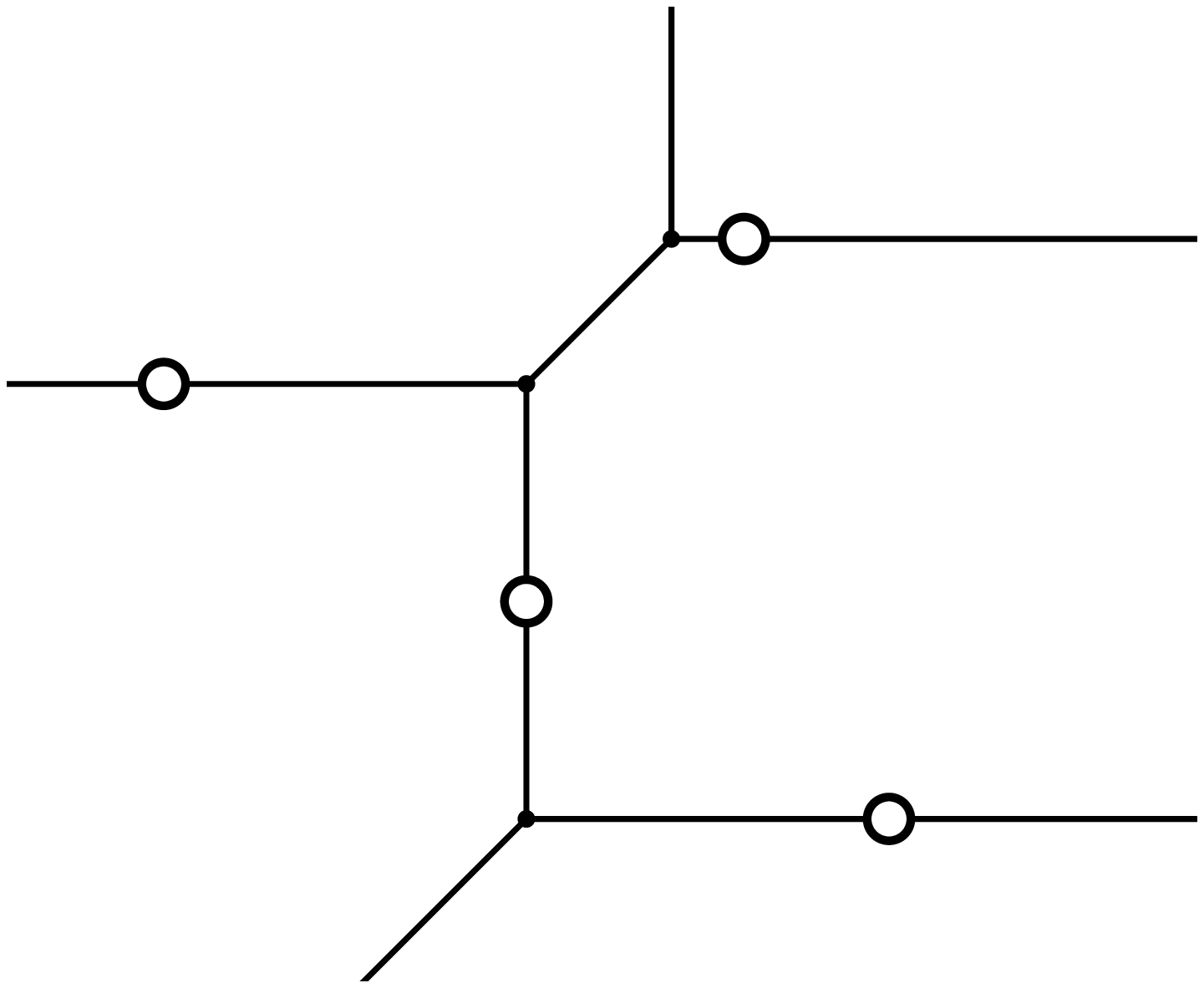} &
  \includegraphics[width=3.3cm]{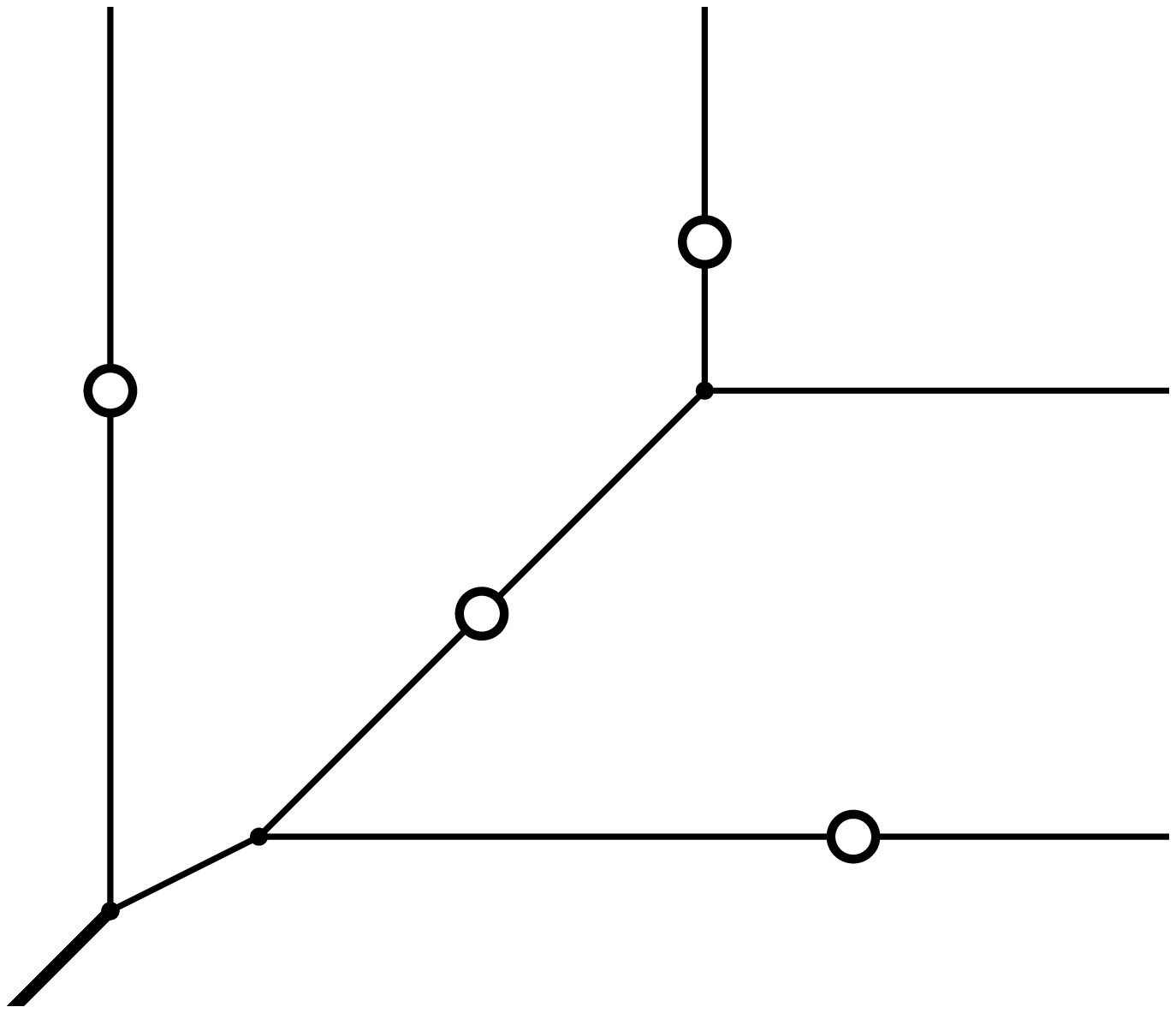} &
  \includegraphics[width=3.3cm]{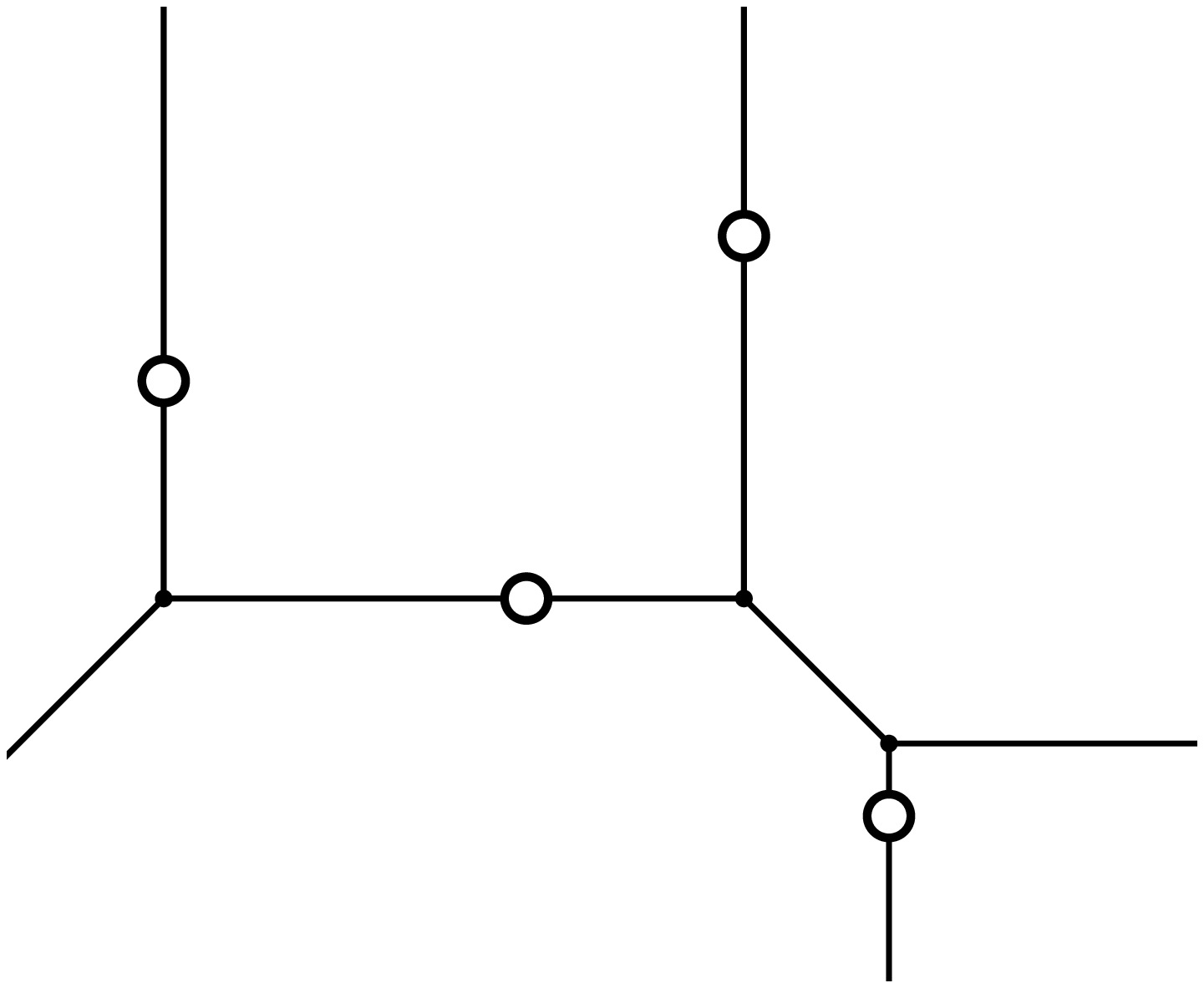} \\ [1ex]
  \includegraphics[width=3.3cm]{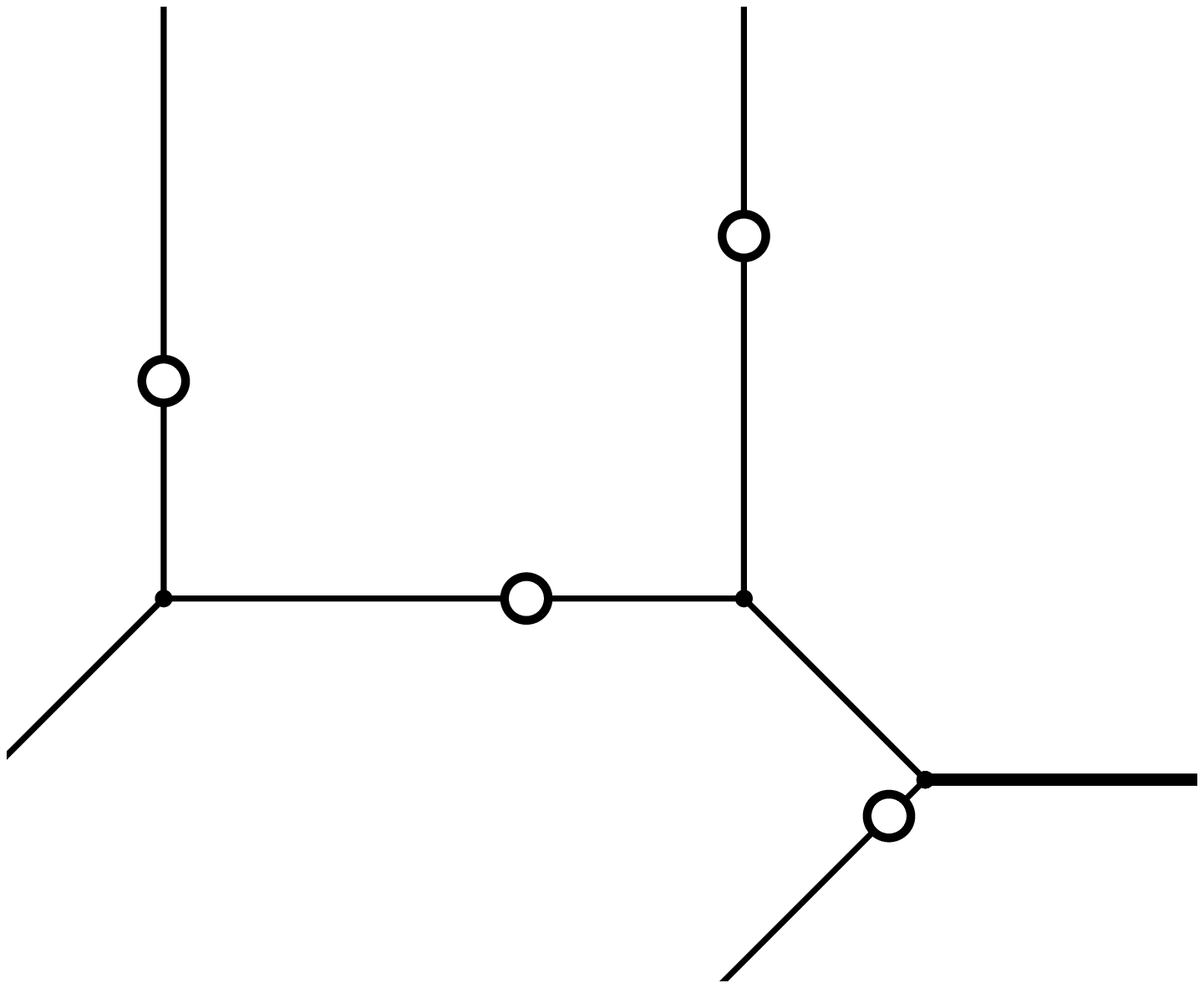} &
  \includegraphics[width=3.3cm]{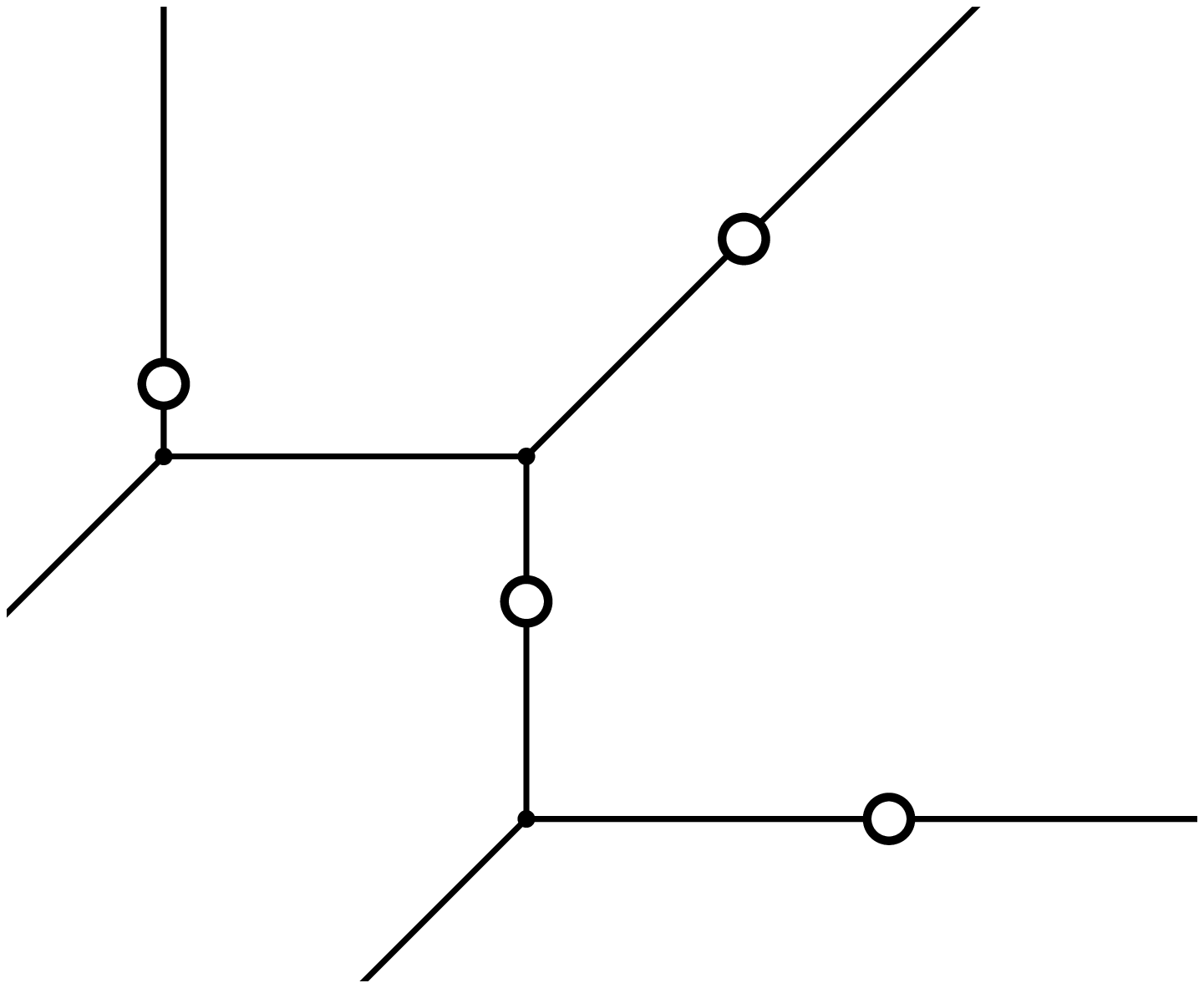} &
  \includegraphics[width=3.3cm]{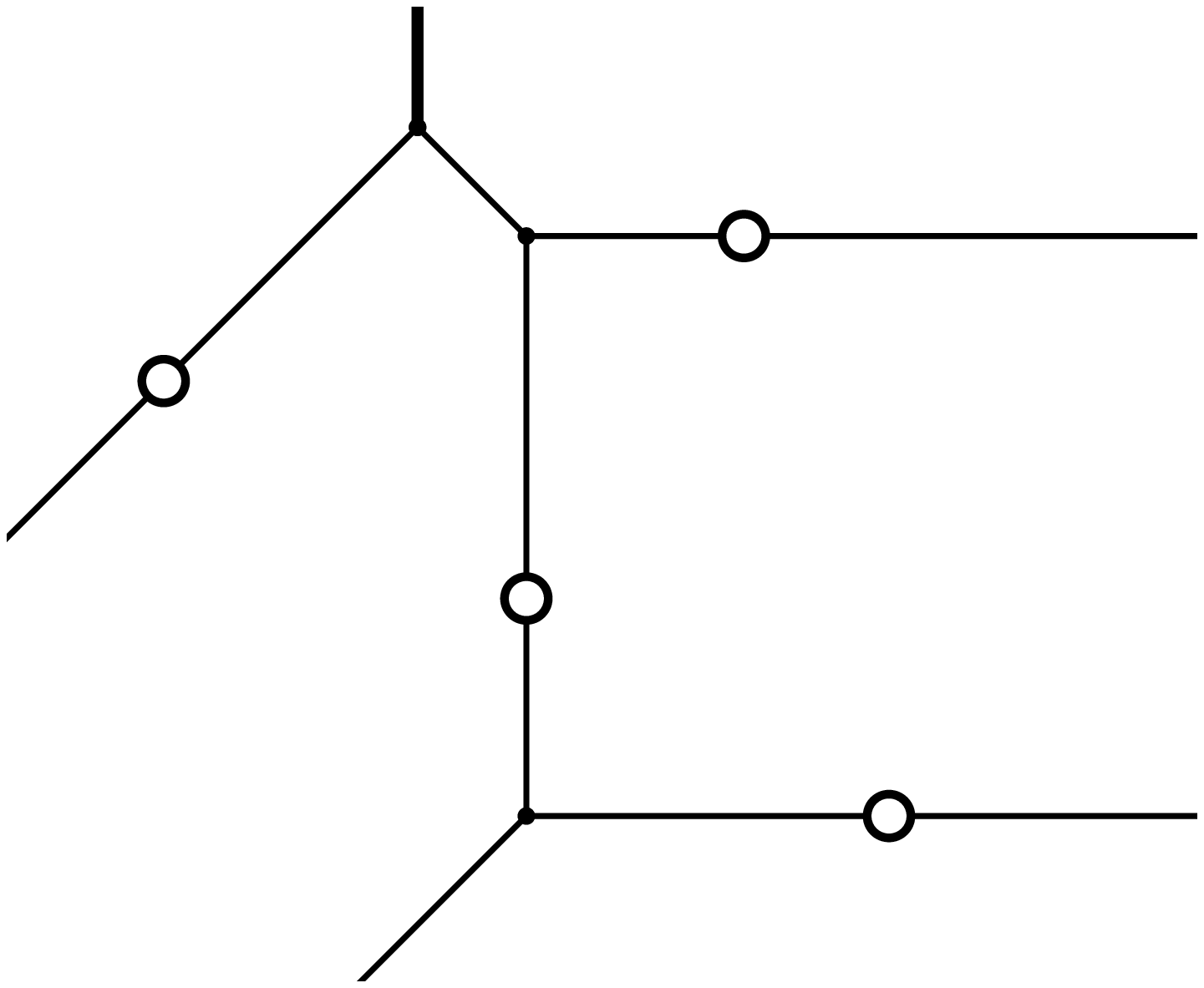}
  \end{array} 
\]

\vspace*{0cm}

\caption{The six limits conics through four given points}
\label{fi:extremeconics}
\end{figure}
\fi

\end{enumerate}

Figure~\ref{fi:pencilconics} depicts the tree 
of conics through the four points $(0,6,0)$, $(5,3,0)$, $(10,0,0)$, $(8,8,0)$.
The six leaves are the limit conics in 
Figure \ref{fi:extremeconics}.

\ifpictures
\begin{figure}[h]
\vspace*{0cm}

\[
  \includegraphics[width=13.2cm]{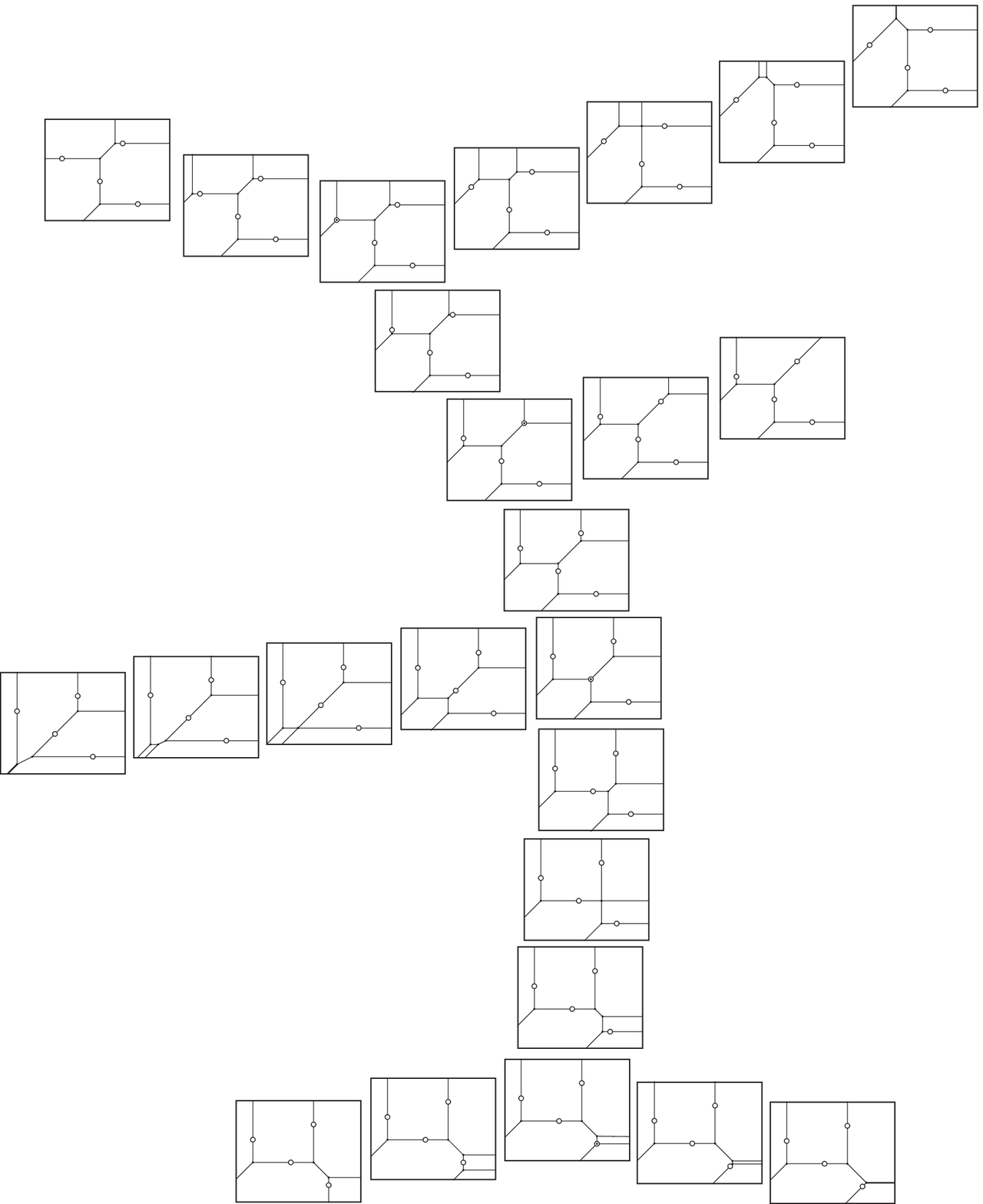}
\]

\vspace*{0cm}

\caption{The tree of conics through four given points}
\label{fi:pencilconics}
\end{figure}
\fi

Hence, the \emph{combinatorial type} of a pencil of conics 
is described by a tree with six labeled leaves
$x^2,xy,y^2,yz,z^2,xz$
and four trivalent (unlabeled) internal vertices.
Since the number of trivalent trees with $n$ labeled leaves is the Schr\"oder number
\[
  (2n-5) !! = (2n-5) \cdot (2n-7) \cdots 5 \cdot 3 \cdot 1 \, ,
\]
there are at most 105 different combinatorial types of pencil of conics.
These 105 trees come in two symmetry classes:

\begin{description}
\item[Caterpillar trees] There are 90 trees in which each of the
four vertices is adjacent to at least one half ray
(such as the tree in Figure~\ref{fi:pencilconics}).
\item[Snowflake trees] There are 15 trees in which there exists
one vertex which is not adjacent to any half ray.
\end{description}

We say that a tree $\Gamma$ {is \emph{realizable by a pencil of
conics} if there exists a configuration of four points in $\tp^2$
whose pencil of conics gives the tree $\Gamma$. 
The following statements characterize which of the 105 trees 
are realizable.

\begin{thm}
\label{th:penciltypes}
A tree $\Gamma$ is realizable if and only if $\Gamma$ can be
embedded as a planar graph into the unit disc such that the 
six labeled vertices are located 
in the cyclic order $x^2, xy, y^2, yz, z^2, zx$
on the boundary of the disc.
\end{thm}

\begin{cor}
\label{co:penciltypes}
Exactly 14 of the 105 trees are realizable. Twelve of them
are caterpillar trees and two of them are snowflake trees.
\end{cor}

Figure~\ref{fi:penciltypes} shows one of the realizable caterpillar trees
and one of the realizable snowflake trees. In particular, the statement 
shows that a snowflake tree can be obtained as the complete intersection 
of four tropical hyperplanes in $\tp^5$. By 
Example 6.2 in \cite{speyer-sturmfels-2003}), the
tropical $3$-plane  in $\tp^5$
which is dual to the snowflake tree is 
not a complete intersection.

\ifpictures
\begin{figure}[h]
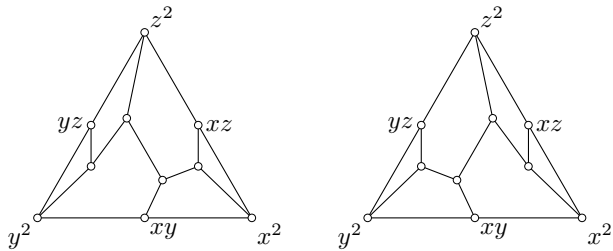

\vspace*{-0.2cm}

\[
  \includegraphics[scale=0.95]{picturesmpost/pictropical.1} \qquad
  \includegraphics[scale=0.95]{picturesmpost/pictropical.2}
\]

\vspace*{-0.2cm}

\caption{A realizable caterpillar and a realizable snowflake}
\label{fi:penciltypes}
\end{figure}
\fi

In order to prove Theorem~\ref{th:penciltypes}, we consider the following
more general setting.
Let $\mathcal{A} = \{a_1,a_2,\ldots,a_n\}$ be an $n$-element subset of
$\{(i,j,k) \in \nn_0^3 : i+j+k = d \}$ for some $d$.
Consider a tree $\Gamma$ with $n$ leaves which are
labeled by the elements of $\mathcal{A}$.  The combinatorial type
of the tree $\Gamma$ is specified by the induced subtrees
on four leaves $\{a_i,a_j,a_k,a_l\}$. The four possible
subtrees are denoted as follows:
$$ (ij \vert kl) \, , \,\,
(ik \vert jl) \, , \,\,
(il \vert jk) \, , \,\, (ijkl).$$
The first three trees are the trivalent trees on $i,j,k,l$
and the last tree is the tree with one $4$-valent  node.
We say that $\Gamma$ is \emph{compatible with} $\mathcal{A}$
if the following condition holds: if $(ij \vert kl)$ is  a 
trivalent subtree of $\Gamma$, then the convex hull
of $a_i, a_j, a_k$, and $a_l$ has at least one of the
segments ${\rm conv}(a_i,a_j)$ or
${\rm conv}(a_k,a_l)$ as an edge.

The space of curves with support $\mathcal{A}$ is
identified with the tropical projective space $\tp^{n-1}$.
Now consider a configuration $C$ of $n-2$ points
in $\tp^2$ which does not lie on any tropical curve
with support $\mathcal{A} \backslash \{a_i,a_j\}$ for any pair $i,j$.
Then the set of all curves with support $\mathcal{A}$ which pass
through $C$ is a tropical line  $\Gamma_C$ in $\tp^{n-1}$.
Combinatorially, the line $\Gamma_C$ is a tree whose leaves
are labeled by $\mathcal{A}$.

\begin{thm} For every configuration $C$,
the tree $\Gamma_C$ is compatible with $\mathcal{A}$.
\end{thm}

\begin{proof}
The theorem is trivial for $n \leq 3$. In the case $n=4$, it
is proved by examining all combinatorial
types of four-element-configurations of $\mathcal{A}$ and
two-point configurations  
in the plane. This involves an exhaustive case analysis
which we omit here. For the general
case $n \geq 5$, we assume by induction that the result
is already known for $n-1$.

The Pl\"ucker coordinate $p_{ij}$ of the line $\Gamma_C$
is the tropical determinant of the $(n-2)\times (n-2)$-matrix
whose rows are labeled by $C$, whose columns
are labeled by  $\mathcal{A} \backslash \{a_i,a_j\}$ and whose
entries are the dot products of the row labels
with the column labels. Fix any quadruple $ijkl$
and consider the restricted Pl\"ucker vector 
$$ P \quad = \quad (p_{ij}, p_{ik}, p_{il}, p_{jk}, p_{jl}, p_{kl}). $$
Supposing that $n \not\in \{i,j,k,l\}$, we can compute these six
tropical $(n-2) \times (n-2)$-determinants by tropical Laplace expansion
with respect to the last row $n$:
$$ P \quad = \quad \bigoplus_{c \in C} \, (a_n \cdot c) \odot  
(p_{ij}^{(c)}, p_{ik}^{(c)}, p_{il}^{(c)}, p_{jk}^{(c)}, p_{jl}^{(c)}, p_{kl}^{(c)}). $$
Here $p_{ij}^{(c)}$ is the tropical $(n-3)\times (n-3)$-minor
gotten from the matrix for $p_{ij}$ by deleting column $n$
and row $c$.  Each of the Pl\"ucker vectors in this sum
defines a tree which is compatible with
the set $\{a_i,a_j,a_k,a_l\}$, since
 $\Gamma_{C \backslash \{c\}}$ is compatible with
 $\mathcal{A} \backslash \{a_n\}$ by induction.
 The proof now follows from a lemma to the effect
 a tropical linear combination
 of compatible Pl\"ucker vectors is always compatible.
\end{proof}

\begin{cor} \label{co:compatible}
If ${\rm conv}(\mathcal{A})$ is a convex polygon 
which has $a_1,a_2,\ldots,a_n$ on its boundary, then
the compatible trees are precisely the planar trees 
whose leaves form an $n$-gon.  The number of compatible
trees is the Catalan number $\,\frac{1}{n-1} \binom{2n-4}{n-2}$.
\end{cor}

We do not know whether every tree $\Gamma$ which is compatible
with $\mathcal{A}$ can be realized as $\Gamma = \Gamma_C$
by some configuration $C$ of $n-2$ points in $\tp^2$. For the special case
$$ \mathcal{A} \quad = \quad
\bigl\{ \, (2,0,0), \,(1,1,0), \, (0,2,0), \, (0,1,1) , \, (0,0,2), \, (1,0,1) \, \bigr\} ,$$
there are $14$ compatible trees, by Corollary~\ref{co:compatible}, and we checked
that each of them is realizable. This proves Theorem~\ref{th:penciltypes}
and Corollary~\ref{co:penciltypes} on quadratic curves through four given points.

\section{Incidence Theorems and Tropical Cinderella}

The previous sections showed that central concepts of classical projective
geometry (such as intersection multiplicities, B\'ezout's and Bernstein's Theorem,
as well as determinants) can be transferred to tropical geometry. In this section, 
we report on first insights on the question 
in how far projective incidence theorems 
carry over to the tropical world. By pointing out some pitfalls, we
would like to argue that these generalizations require great care. 

Our investigations were supported by a prototype of a {\it tropical version} of 
the dynamic geometry software {\tt Cinderella} 
\cite{CINDERELLA}
(most pictures in this article were generated with this tool). This software supports the interactive manipulation of elementary geometric constructions. 
Here, an elementary geometric construction is considered as a {\it construction sequence} 
that starts with a set of free elements (e.g., points) and proceeds by constructively 
adding new dependent elements (e.g., the line passing through two points, 
intersection points of two lines, conics through five points).
Once a construction is finished, one can explore its dynamic behavior by simply dragging 
the free elements. The dependent elements move according to 
the construction sequence. 
Our experimental version of this software provides basic operations for {\it join}
(i.e., the line passing through two points), 
{\it meet} (i.e., the intersection point of two lines)
and the {\it stable conic through five points} in $\tp^2$, as discussed
in Section 5.

The possibilities and limitations of dynamic geometry are tightly connected to the
degenerate situations which can occur. E.g., a real user might choose 
two points $a, b$ in the plane which are identical, in which case
the line through $a$ and $b$ is neither unique in classical projective geometry nor
unique in tropical geometry. Thus, for the purposes of dynamic geometry software
it is very desirable to have as few degeneracies as possible.
In contrast to classical geometry, tropical geometry offers the 
distinguished feature to have absolutely no degeneracies in our basic
operations. Namely, as discussed in Sections~\ref{se:bezout} and Section~\ref{se:cramer},
the concept of stable intersections always defines a distinguished solution, 
say, for the line passing through two given points $a$ and $b$. 
This holds true even if the two points coincide or if there is 
an infinite number of tropical lines passing through $a$ and $b$.

Let $a \otimes b$ denote the tropical cross product
\[ a\otimes b := (
    a_2\odot b_3 \, \oplus \, a_3 \odot b_2, \,
    a_3\odot b_1 \, \oplus \, a_1 \odot b_3, \,
    a_1\odot b_2 \, \oplus \, a_2 \odot b_1
) \, .
\]
The \emph{stable join} of two points $a,b \in \tp^2$ is the tropical
line 
\[
  \mathcal{T}(u_1 \odot x \oplus u_2 \odot y \oplus u_3 \odot z)
\] 
defined by $u:=a \otimes b$.
Similarly, the \emph{stable meet} of two lines
\[
  \mathcal{T}(u_1 \odot x \oplus u_2 \odot y \oplus u_3 \odot z) \quad \mbox{and} \quad
  \mathcal{T}(v_1 \odot x \oplus v_2 \odot y \oplus v_3 \odot z) 
\]
is the point $u \otimes v$.

Three (tropical or projective) lines $a,b,c$ are said to be 
\emph{concurrent} if $a$, $b$, and $c$ have a point in common.
Pappus' Theorem is concerned with certain concurrencies among nine lines.
In classical projective geometry, it is well-known that one has to be careful
in stating the right non-degeneracy assumptions.
However, we will see that we have to be even more careful in the tropical world.

\ifpictures
\begin{figure}[h]
\vspace*{-0.2cm}

\[
  \includegraphics[width=6.3cm]{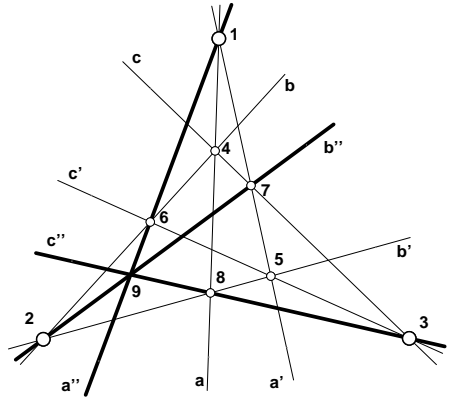}
\]

\vspace*{-0.3cm}

\caption{Pappus' theorem in classical projective geometry. The lines $a''$, $b''$ and
$c''$ are drawn in bold.}
\label{fi:pappos}
\end{figure}
\fi

\medskip
The following statement expresses one version of Pappus' Theorem that holds in the usual projective plane 
over an arbitrary field.

\medskip
\noindent
{\sc Pappus' Theorem, incidence version:}
{\sl Let $a,a',a'',b,b',b'',c,c',c''$ be nine distinct lines in the projective plane.
If the following triples of lines 
\[
[a,a',a''],[b,b',b''],[c,c',c''],[a,b,c],[a',b',c'],[a'',b,c'],[a',b'',c],[a,b',c'']
\]
are concurrent then also $[a'',b'',c'']$ are concurrent.}
\medskip

In Figure~\ref{fi:pappos}, the points in which the lines meet are labeled by 
$1, \ldots, 8$. 
The final common intersection point of the three lines $a'',b'',c''$ is labeled by $9$. 

Experimentally, it turns out that a tropical analogue of this statement holds for
many instances.
However, Figure~\ref{fi:pappos} shows a counterexample, which proves that 
the above version of Pappus' Theorem does not generally hold in the tropical plane.

\ifpictures
\begin{figure}[h]
\vspace*{0cm}

\[
  \includegraphics[width=10cm]{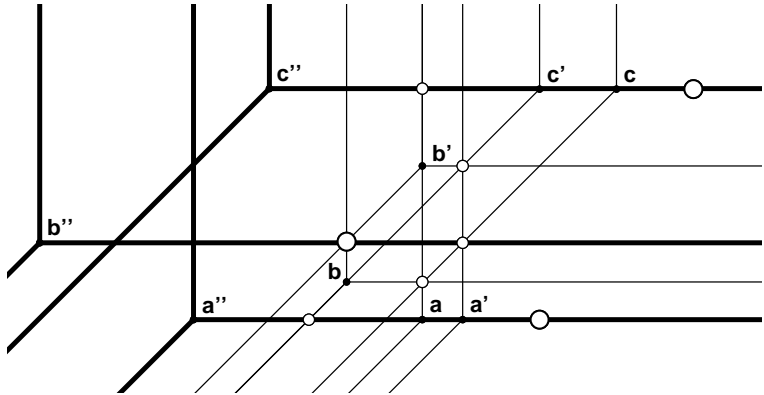}
\]

\vspace*{0cm}

\caption{A tropical non-Pappus configuration: 
the triples $[a,a',a'']$, $[b,b',b'']$, $[c,c',c'']$, $[a,b,c]$, $[a',b',c']$, $[a'',b,c']$, $[a',b'',c]$, $[a,b',c'']$ are concurrent, but 
$[a'',b'',c'']$ is not.
}
\label{fi:nonpappos}
\end{figure}
\fi

In this picture all concurrencies of the hypotheses are satisfied, but the conclusion is violated. Concrete coordinates for the lines in this counterexample are given by the following matrix:
\[
\bordermatrix{
&a&b&c&a'&b'&c'&a''&b''&c''\cr
&-4 &-2 &-9 &-5  &-4  &-7  &2    &6   &0  &\cr
&6&5&0 &6 &2  &0    &6 &4  &0  &\cr
&0&0&0&0&0&0&0&0&0\cr
}
\]

The main reason why the conclusion of the tropical incidence version of 
Pappus' Theorem holds for many examples is that there 
is also a constructive version in the projective plane that we conjecture to 
hold also in the tropical projective plane.

\medskip

\begin{samepage}
\noindent
{\sc Pappus' Theorem, constructive version:}
{\sl Let $1,2,3,4,5$ be five freely chosen points in the projective plane given by homogeneous coordinates.
Define the following additional three points and nine lines by a sequence of 
(stable) join and (stable) meet operations (carried out by cross-products):
\[
\begin{array}{llllll}
a:=1\otimes 4,&
b:=2\otimes 4,&
c:=3\otimes 4,&
a':=1\otimes 5,&
b':=2\otimes 5,&
c':=3\otimes 5,\cr
6:=b\otimes c',&
7:=a'\otimes c,&
8:=a\otimes b',&
a'':=1\otimes 6,&
b'':=2\otimes 7,&
c'':=3\otimes 8.\cr
\end{array}
\]
Then the three tropical lines $a'',b''$ and $c''$ are concurrent.}
\end{samepage}

\medskip
The construction is organized in a way such that the eight hypotheses of the 
incidence-theoretic version of the theorem are satisfied
automatically by the construction. 
For instance, $a,a',a''$ meet in a point (namely point $1$) since they all arise from a
join operation in which $1$ is involved.
For some choices of the initial points it may happen that during the construction 
the cross product of two linearly dependent vectors is calculated. In this case the final conclusion is automatic.

\ifpictures
\begin{figure}[h]
\vspace*{0cm}

\[
  \includegraphics[width=6.0 cm]{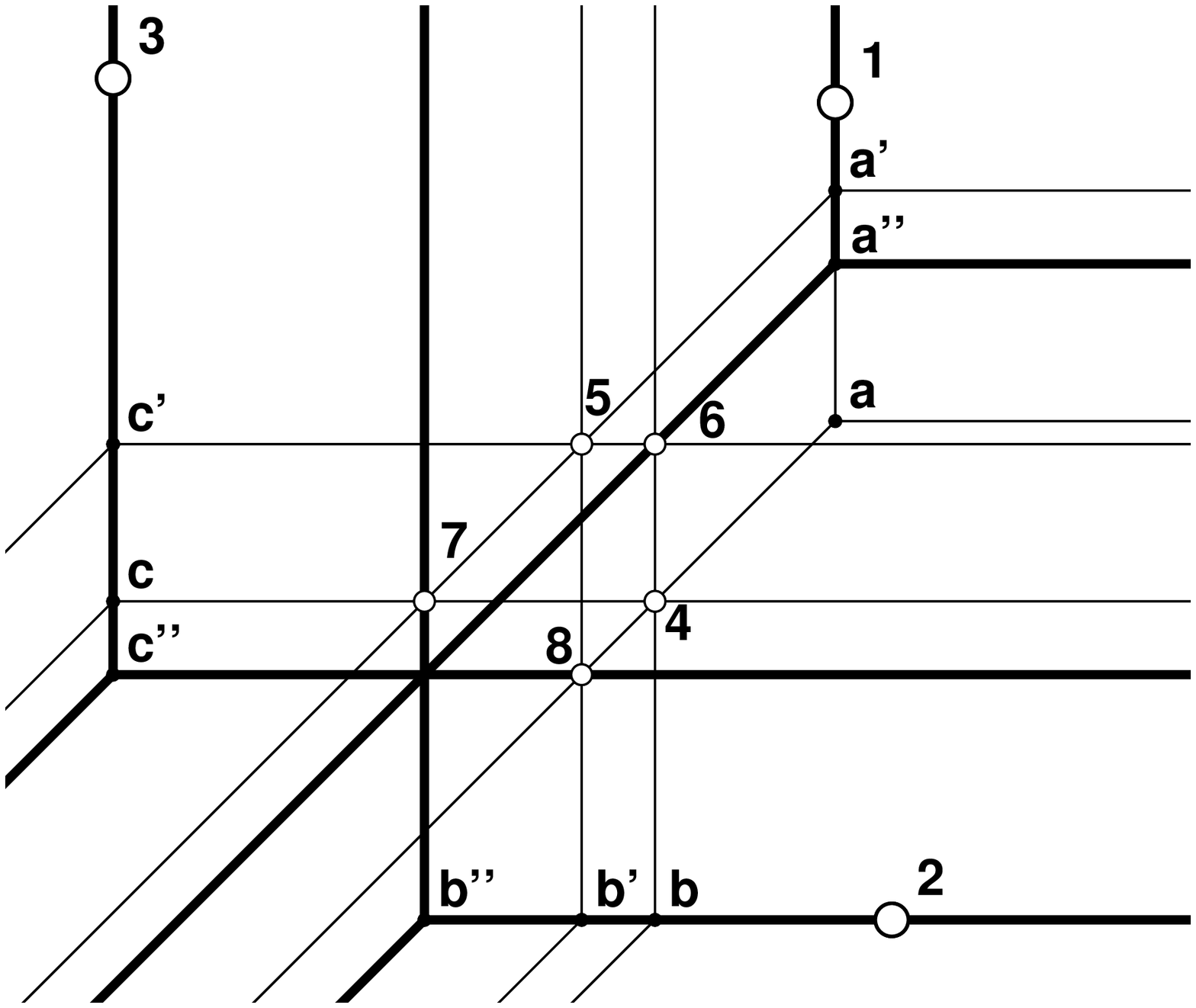}\quad
  \includegraphics[width=6.0 cm]{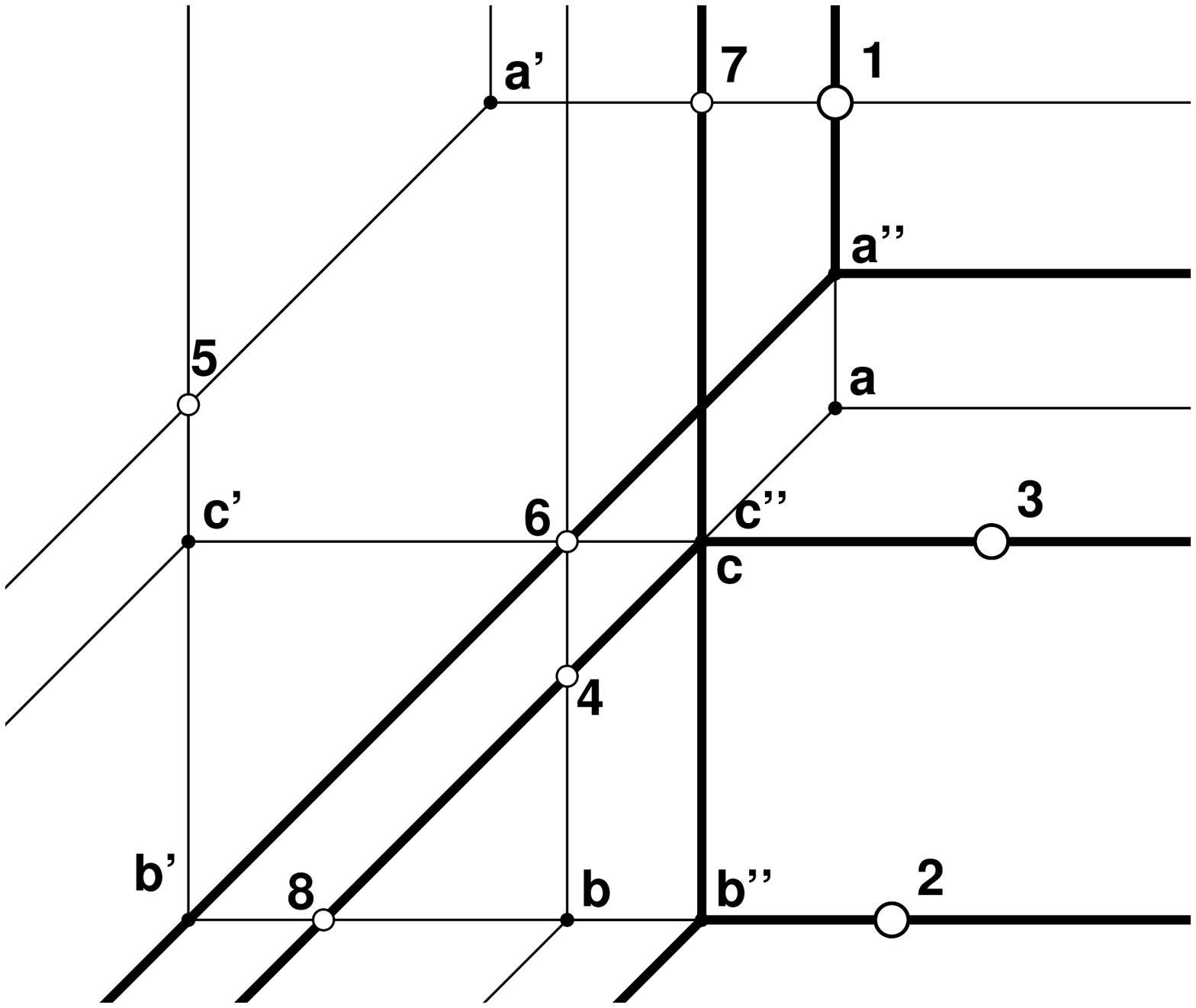}
\]

\vspace*{0cm}

\caption{The constructive tropical Pappus' Theorem.}
\label{fi:troppappos}
\end{figure}
\fi

There is strong experimental evidence that this constructive version of 
Pappus' Theorem also holds in the tropical projective plane. Here, the vanishing of the
determinant is replaced by the property that
the matrix with rows $a''$, $b''$, and $c''$ is tropically singular.
The construction sequence is carried out by using the tropical cross product.
As mentioned before, degenerate cross products cannot occur. 
 
Figure~\ref{fi:troppappos} shows two possible situations of how the final conclusion of the theorem can hold. 
Either the three lines are in a degenerate position (right picture) or they meet properly (left picture). 
In the latter case, the final coincidence arises since there exists an interesting subconfiguration in the picture that is an incidence theorem of classical affine geometry.
This subconfiguration is formed by the points $4,5,6,7,8$ and the three {\it straight} lines passing through them (which are rays of the nine tropical lines). The rays form three bundles of parallel lines. If the incidences at $4,5,6,7,8$ are satisfied, then the final 
coincidence
of the bold lines are satisfied as well. So the tropical constructive Pappus' Theorem
seems to hold since the final concurrence arises either from degeneracy or from an affine incidence theorem.

So far, we do not have a proof for the tropical constructive version of Pappus' Theorem.
However, let us point out that in principle,
one can decide the validity of the tropical version by an exhaustive 
enumeration of all combinatorial equivalence classes of realizations of the hypotheses 
(which correspond to cones in the hypotheses-space of the construction). 
However, this requires to check many different situations, since 
already a complete quadrilateral (i.e., joining four points by all six possible lines)  can 
be tropically realized in 3141 different combinatorial ways. 
\penalty -10000

\end{document}